\documentclass[12pt,letter,oneside]{amsbook}
\usepackage{amssymb}

\newlength{\originalbase}
\newcommand{\spacing}[1]{\setlength{\baselineskip}{#1\originalbase}}

\newcommand{\mb}[1]{\mbox{\rm {#1}}}
\newcommand{\mr}[1]{\mathfrak{#1}}
\newcommand{\mc}[1]{\mathcal{#1}}

\newcommand{\C}{\mathbb{C}}
\newcommand{\Q}{\mathbb{Q}}
\newcommand{\R}{\mathbb{R}}
\newcommand{\Z}{\mathbb{Z}}
\newcommand{\F}{\mathbb{F}}

\newcommand{\gen}{\mathrm{gen}}
\newcommand{\triv}{[0]}

\newtheorem{lm}{Lemma}[section]
\newtheorem{thq}[lm]{Theorem}
\theoremstyle{remark}
\newtheorem{reme}[lm]{Remark}
\newtheorem{df}[lm]{Definition}
\newtheorem{co}[lm]{Corollary}
\newtheorem{con}[lm]{Conjecture}
\newtheorem{prop}[lm]{Property}
\newtheorem{ex}[lm]{Example}

\DeclareMathOperator{\ad}{ad}
\DeclareMathOperator{\AD}{A\!D}
\newcommand{\ADD}[2]{\sideset{_{#1}}{_{#2}}\AD}
\DeclareMathOperator{\Del}{\Delta}
\newcommand{\DEL}[2]{\sideset{_{#1}}{_{#2}}\Del}

\DeclareMathOperator{\Sym}{Sym}
\DeclareMathOperator{\mdeg}{mdeg}
\DeclareMathOperator{\res}{res}
\DeclareMathOperator{\rank}{rank}

\DeclareMathOperator{\tr}{tr}
\DeclareMathOperator{\im}{Im}
\DeclareMathOperator{\id}{id}
\DeclareMathOperator{\spann}{span}

\DeclareMathOperator{\Aut}{Aut}

\DeclareMathOperator{\Inn}{Inn}

\DeclareMathOperator{\IAut}{I\!Aut}
\DeclareMathOperator{\CIAut}{CI\!Aut}
\DeclareMathOperator{\Hom}{Hom}
\DeclareMathOperator{\D}{Der}

\DeclareMathOperator{\IGamma}{I\Gamma}

\DeclareMathOperator{\Out}{Out}

\DeclareMathOperator{\Cen}{Cen}
\DeclareMathOperator{\End}{End}

\DeclareMathOperator{\GL}{GL}
\DeclareMathOperator{\SL}{SL}

\newcommand{\gll}{\mathfrak{gl}}
\newcommand{\sll}{\mathfrak{sl}}

\newcommand{\la}{\langle}
\newcommand{\ra}{\rangle}

\newcommand{\ifnewpage}{}

\begin{document}						

\spacing{1.70}  
\thispagestyle{empty}


\begin{center}

{\large \bf Abstract}
\bigskip

{ \Large  \bf
On the Automorphism Tower \\
of Free Nilpotent Groups
}

\bigskip

{\Large
Martin Dimitrov Kassabov\\
2003
}
\end{center}

\vspace{2cm}
\begin{center}
\begin{minipage}{4.5in}
\rm

\mbox{\hspace{5mm}}
In this thesis I study the automorphism tower of free nilpotent
groups. Our main tool in studying the automorphism tower is to
embed every group as a lattice in some Lie group. Using known
rigidity results the automorphism group of the discrete group can
be embedded into the automorphism group of the Lie group.
\vspace{2mm}

\mbox{\hspace{5mm}}
This allows me to lift the description of the derivation tower of
the free nilpotent Lie algebra to obtain information about the
automorphism tower of the free nilpotent group.
\vspace{2mm}

\mbox{\hspace{5mm}}
The main result in this thesis states that the automorphism 
tower of the free nilpotent group $\Gamma(n,d)$ on $n$ generators 
and nilpotency class $d$, stabilizes after finitely many steps.
If the nilpotency class is small compared to the number of 
generators we have that the height of the automorphism tower is
at most 3.

\end{minipage}
\end{center}


\newpage


\begin{titlepage}
\vspace*{1cm}
\begin{centering}
\begin{LARGE} \bf                         
	On the Automorphism Tower \\
of Free Nilpotent Groups\\
\end{LARGE}                           
\vfill
	A Dissertation\\
	Presented to the Faculty of the Graduate School\\
	of\\
	Yale University\\
	in Candidacy for the Degree of\\
	Doctor of Philosophy\\
\vfill
	by\\
	Martin Dimitrov Kassabov\\
\vspace{1cm}
	Dissertation Director: Professor Efim Zelmanov\\
\vspace{1cm}
	May 2003\\
\end{centering}
\end{titlepage}


\thispagestyle{empty}

\begin{titlepage}
	\begin{centering}
	\vspace*{3in}
		\copyright 2003 by Martin Dimitrov Kassabov.\\
		All rights reserved.\\
	\end{centering}
\end{titlepage}

\pagenumbering{roman}

\thispagestyle{empty}
\tableofcontents
\thispagestyle{empty}
\newpage



\newpage
\addcontentsline{toc}{chapter}{Notations}
\thispagestyle{empty}
\vspace{5mm}
\begin{center}
\bf \LARGE Notations
\end{center}
\vspace{5mm}
	
Here we list the notations, used in this diseration. 

\begin{tabular}{ll}
$\Z$ & set of integres \\
$\Q$ & field of rational numbers \\
$\R$ & field of real numbers \\
$\F_p$ & finite field with $p$ elements\\
$\gll_n$ & Lie lagebra of $n\times n$ matrices \\ 
$\sll_n$ & Lie lagebra of traceless $n\times n$ matrices \\ 

$\R \la V \ra$ & free associative algebra over $\R$ generated by a spcase $V$ \\
$L[V]$ & free Lie algebra generated by a spcase $V$ \\
$L[V,d]$ & free nilpotent of class $d$ Lie algebra generated by a spcase $V$ \\
$\D L$ & derivation algebra of the algebra $L$ \\
$\D_0 L$ & algebra restricted derivations of the algebra $L$ \\
$\D^k L$ & $k$-th algebra in the derivation tower of $L$ \\
$\D^k_0 L$ & $k$-th algebra in the tower of restricted derivations of $L$ \\
$\Cen(L)$ & center of the algebra $L$ \\
$\Cen_L(I)$ & centralizer of the ideal $I$ in the algebra $L$ \\
$[A,B]$ & commutator of the idales $A$ and $B$ \\

$UL[V]$ & nilpotent part of the derivation algebra of $L[V]$\\
$UL[V,d]$ & nilpotent part of the derivation algebra of $L[V,d]$\\

$\GL_n(R)$ & group of invertable $n\times n$ matrices over a ring $R$ \\
$\SL_n(R)$ & group of $n\times n$ matrices with determinant $1$ over a ring $R$ \\
$\SL_n^{\pm 1}(R)$ & group of $n\times n$ matrices with determinant $\pm 1$ over a ring $R$ \\

$\mc{F}_n$ & free group on $n$ generators \\
$\Gamma(n,d)$ & free nilpotent group on $n$ generators of nilpotency class $d$\\
$G(n,d)$ & free unilpotent group on $n$ generators of nilpotency class $d$\\

$\Aut \Gamma$ & automorphism group of the group $\Gamma$ \\
$\Aut^k \Gamma$ & $k$-th group in the automorphism tower of the group $\Gamma$ \\ 
$\Cen(\Gamma)$ & center of the group $\Gamma$ \\
$[A,B]$ & commutator of the subgroup $A$ and $B$ \\

$M^\mr{g}$ & elemnts in the module $M$ fixed by the algebra $\mr{g}$ \\
$M^G$ & elemnts in the module $M$ fixed by the group  $G$ \\

\end{tabular}

\newpage

\addcontentsline{toc}{chapter}{Introduction}
\thispagestyle{empty}
\vspace{5mm}
\begin{center}
\bf \LARGE Introduction
\end{center}
\vspace{5mm}

	\begin{section}{General}

A classical result of Wielandt~\cite{Wielandt} from the 1930-es
asserts that the automorphism tower stabilizes for any finite
group $\Gamma$ with trivial center. The original proof 
was simplified by Schenkman in~\cite{Schenkman2}. This proof 
can be modified to
yield an effective bound on the height (see~\cite{Pettet}).

There are two possible directions to generalize this result, one
is to start with an infinite group $\Gamma$ and the other is to
start with a group with non trivial center. The straightforward
generalization of Wielandt's result does not hold in any of these
two cases.

If the original group $\Gamma$ has trivial center the automorphism
tower consist of an increasing sequence of group and this sequence
can be extended to a transfinite sequence. Recently
Hamkins~\cite{Hamkins} proved that the transfinite automorphism tower
stabilizes, at some (transfinite) level, for any group $\Gamma$
(even not finitely generated).

To the best of our knowledge there is no generalization of
Wielandt's theorem for a large class of groups with centers. It is
not known if the analog of this theorem holds for finite
$p$-groups. One of the few cases where the generalization of
Wielandt's result is known is the case of finitely generated
abelian groups.

The most general version of Wielandt's theorem for infinite groups
with center does not hold --- there are known examples of
polycyclic groups with infinite automorphism tower
(see~\cite{Hulse}). Knowing that the automorphism tower stabilizes
for abelian groups and does not stabilize for polycyclic groups,
Baumslag conjectured that for any finitely generated nilpotent
group $\Gamma$, the automorphism tower of $\Gamma$ stabilizes
after finitely many steps.

There has been almost no progress in proving Baumslag's conjecture
during the last 50 years. In~\cite{DyerForm2}, Formanek and Dyer
proved that the automorphism tower of the free nilpotent group of
class $2$ stabilizes at the second level, when the number of
generators is different than $3$. For the three generator case,
the tower stabilizes at the third level. Similar results are known
for other nilpotent groups, like the group of $n\times n$ upper
triangular matrices, and in all these case the automorphism tower
is very short.

The description of the automorphism tower is also know for the case
of free groups. It was proved by Dayer and formanek (see~\cite{DyerForm1})
in the case of finitely generated groups and by Tolstykh (see~\cite{Tolstykh})
in the case of infinitely many generators.

Baumslag's conjecture was verified for several classes of groups
and no example of a nilpotent group with infinite automorphism
tower has been found. I started this project hoping to find a
counter example among the free nilpotent groups. Some computations
in the free Lie algebras and their derivations suggested that it
may be possible to find such a counterexample among the free $2$
generated nilpotent groups of odd nilpotency class.

It turns out (see Theorems~\ref{auttower:st}, \ref{grtower:stab}
and~\ref{auttowerrst2gen}) that the free nilpotent group is too
big and its automorphism tower stabilizes after finitely many
steps. I still have some hope that it is possible to construct an
example of a relatively free group in some nilpotent variety,
whose automorphism tower does not stabilize. However it seems that
in this case the tower is going to stabilize weakly, i.e., for
some $i$ the groups $\Aut^i \Gamma$ and $\Aut^{i+1} \Gamma$ will
be isomorphic but the natural homomorphism between them, sending
every element to the inner automorphism obtained by conjugation
with this element will not be an isomorphism.

\bigskip

The main result of this dissertation is the following theorem

{\renewcommand{\thelm}{\!\!}
\begin{thq}
The automorphism tower of the free nilpotent group $\Gamma(n,d)$
on $n$ generators and nilpotency class $d$ stabilizes after
finitely many steps.
\end{thq}
}

The idea of the proof of the above theorem is to embed every group
which appears in the automorphism tower of $\Gamma$ into a
suitable Lie group over the reals and first study the automorphism
group of this Lie group using the Lie theory. Using the connection
between that automorphism of a Lie group and the derivations of
its Lie algebra I first study the derivation tower of the free
nilpotent Lie algebra. It turns out that derivation tower is very
short and that it is not a good approximation of the automorphism
tower of the free nilpotent group (although it is a very good
approximation of the automorphism tower of the free unipotent
group). In chapter $3$, I define the notions of `restricted
derivation algebra' and `tower of restricted derivations' which is
a better approximation to the automorphism tower of the free
nilpotent group. After obtaining a description of this tower I
used rigidity results to obtain a similar description of the
automorphism tower of the free nilpotent group.
\end{section}

\ifnewpage

\begin{section}{Structure of the manuscript}

This dissertation consist of $4$ chapters. The first chapter
describes free nilpotent Lie algebras and their derivations. The
main result in this chapter is Theorem~\ref{gen:maincom}, which
describes the generators of the derivation algebra.

In the second chapter I study the automorphism group of the
free nilpotent group, using the derivation algebra of the free
nilpotent Lie algebra. Theorem~\ref{whenT} gives a necessary and
sufficient condition for this group to have Kazhdan property $T$.

The third chapter is devoted to the derivation tower of the free
nilpotent Lie algebras. Since the derivation tower stabilizes
after $3$ steps and it has no connection with the automorphism
tower of the free nilpotent groups, I define the notion of
restricted derivations and study the tower of restricted
derivations. The main results in this chapter are
Theorems~\ref{dertower:nocenter} and~\ref{dertower:center}, which
show that the derivation tower stabilizes and give bound of its
height.

In the last chapter I study the automorphism tower of the free
nilpotent group.  Each group of this tower is embedded in certain
Lie group (over the reals). Using rigidity results I can (in most
cases) embed the next group in the gower into the automorphism
group of the Lie group. This allows me to generalize the results
about the stabilization of the derivation tower to
Theorems~\ref{auttower:st} and~\ref{grtower:stab}, which give the
stabilization of the derivation tower of the free nilpotent groups
in different cases if the number of generators is more then $2$.
The case of $2$ generated free nilpotent groups is a bit different
because not all groups from the automorphism tower can be embedded
into Lie groups. In section~\ref{sec:2gen}, I show how this
problem can be avoided and Theorem~\ref{auttowerrst2gen} gives the
stabilization of the automorphism tower in this case.

\end{section}

\ifnewpage

\begin{section}{Acknowledgments}

Firstly, I would like to thank my adviser Efim Zelmanov for his
support during my studies. I also would like to thank Alex
Lubotzky, for proposing the problem and suggesting the main steps,
and pointing many useful references.

Many fellow graduate students at the department of Mathematics at
Yale University helped me with useful discussions and suggestions
during my work on this thesis. I am especially grateful to Victor
Potsak, Roman Muchnik, Tal Poznansky, Sergei Krutelevich, Misha
Ershov, Keivan Mallahi--Karai, and Alireza Salehi--Golsefidy. I
would also like to thank Jordan Milev, for being a great friend
and roommate during the last 4 years.

Finally, I want to thank my wife Irina and my daughter Dora for
their love and support during the years I've spent in New Haven.

\end{section}

\newpage

\renewcommand{\thesection}{\arabic{chapter}.\arabic{section}}

\pagenumbering{arabic}

\begin{chapter}{Free Nilpotent Lie Algebras}
\thispagestyle{empty}
	
In this chapter we discuss free nilpotent Lie algebras and their
derivations. The main result is Theorem~\ref{gen:maincom} which
describes the generating set of the derivation algebra of the free
nilpotent Lie  algebra if the nilpotency class is small compared
to the number of generators.

\begin{section}{Free Associative and Lie Algebras}
\label{subsec:lie:liealgebra}

Let $V$ be the linear space spanned by (formal letters)
$x_1,\dots,x_n$.
Let $\R\la V\ra = \R \la x_1, \dots, x_n \ra$ be the free
associative algebra over $\R$ generated by $V$, i.e., the algebra
with a basis consisting of all words (including the empty one) on
the letters $x_i$, where the product of two words is their
concatenation.

This algebra is called the free associative algebra because it
satisfies the following universal property:
\begin{prop}
For any associative algebra $A$ (over $\R$) and any elements $y_i
\in A$ there exists a unique homomorphism $\phi: \R\la x_i \ra \to
A$ of associative algebras, such that $\phi(x_i) = y_i$.
\end{prop}

The algebra $\R\la V\ra$ has two natural gradings - the first one
is a $\Z$-grading and the degree of each word in the basis is
given by its length. We will denote the degree of a homogeneous
element $a\in \R\la V\ra$ by $\deg a$. Using the grading the
algebra $\R\la V\ra$ can be written as
$$
\R\la V\ra = \bigoplus_{i=0}^\infty \R\la V\ra^{(i)},
$$
where $\R\la V\ra^{(0)} = \R.1$ and $\R\la V\ra^{(1)} =  V$.

By the universal property of $\R\la V\ra$ it can be seen that the
automorphism group $GL(V)$ of the space $V$ acts naturally on the
algebra and this action preserves the grading. Each homogeneous
component $\R\la V\ra^{(i)}$ is a polynomial representation of
$GL(V)$ which is isomorphic as representation to the $i$-th tensor
power of the standard module. Therefore, it decomposes as a direct
sum of irreducible representations
$$
\R\la V\ra^{(i)} \simeq \bigoplus_\lambda m_\lambda V^\lambda,
$$
where the sum is over all partitions $\lambda$ of $i$ with no more
than $n$ parts and $m_\lambda$ are some positive integers. Here
$V^\lambda$ denotes the irreducible $\GL_n$ module corresponding
to the partition $\lambda$.

The second grading of the algebra $\R\la V\ra$ is given by the
multi degree --- let us denote by $\deg_{x_i} \omega$ the number
of occurrences of the letter $x_i$ in the word $\omega$. The multi
degree of a basis word $\omega$ is defined as
$$
\mdeg \omega  = (\deg_{x_1} \omega, \dots, \deg_{x_n} \omega).
$$
Now we can define the second grading (this grading is $\Z^n$ valued)
of the algebra $\R\la V\ra$. Note that the action of $\GL(V)$ does not
preserve this grading.

Let $I_A$ be the ideal in $\R\la V\ra$ generated by the space $V$.
This is an ideal of codimension $1$ and it is called the
augmentation ideal. The powers of this ideal form a filtration of
the algebra, which coincide with the filtration coming from the
degree grading.

\begin{df}
The quotient $\R\la V,d\ra = I_A/I_A^{d+1}$ is called free
nilpotent algebra of class $d$ because they satisfy a similar
universal property as $\R\la V\ra$ but with respect to all
nilpotent algebras of class $d$ (note that these algebras do not
have a unit). The basis of these algebras consists of all nonempty
words of length less or equal to $d$ and they can be decomposed
into homogenous components as follows:
$$
\R\la V,d\ra = \bigoplus_{i=1}^d \R\la V,d\ra^{(i)} = \bigoplus_{i=1}^d \R\la V\ra^{(i)}.
$$
\end{df}

Let us consider $\R\la V\ra$ as a Lie algebra with respect to the
bracket $[f,g] = fg-gf$, and let $L[V] = L[x_1,\dots,x_n]$ denote
the Lie subalgebra generated by $x_i$-es. The algebra $L[V]$ is
called the free Lie algebra because it satisfies the universal
property:
\begin{prop}
For any Lie algebra $L$ and any elements $y_i\in L$ there exists
a unique Lie algebra homomorphism $\phi: L[x_i] \to L$ such that
$\phi(x_i) =y_i$.
\end{prop}

Since the algebra $L[V]$ is embedded into $\R \la V \ra$ it
inherits from it the two gradings and the $\GL(V)$ module
structure. With respect to the $\Z$ grading the algebra $L[V]$
decomposes as
$$
L[V] = \bigoplus_{i=1}^{\infty} L[V]^{(i)},
$$
where $L[V]^{(i)}$ is the homogenous component of degree $i$ and
we have that $L[V]^{(1)} = V$.

Let $I_L = I_A \cap L[V]$ be the augmentation ideal of the algebra
$L[V]$ - this is the ideal generated by the spaces $V$ which by
definition coincide with the algebra $L[V]$. The powers of this
ideal form a filtration of the algebra, which coincides with the
filtration coming from the degree grading, i.e., we have
$$
L[V]^{(i)} \simeq I_L^i/I_L^{i+1}.
$$

As in the associative case, each homogenous component $L[V]^{(i)}$
is a polynomial representation of $\GL(V)$, which decomposes as a
direct sum of irreducible representations
$$
L[V]^{(i)} \simeq \bigoplus_\lambda l_\lambda V^\lambda,
$$
where the sum is over all partitions $\lambda$ of $i$ with no more
than $n$ parts and $l_\lambda$ are some nonnegative integers. The
multiplicities $l_\lambda$ are positives for all $\lambda$ except
the following: $\lambda = [i]$, for $i >1$; $\lambda = [1^k]$ for $k\leq n$;
and $\lambda = [2^2]$ and $ \lambda = [2^3]$. 

\begin{df}
Finally, we can define the free nilpotent Lie algebra $L[V,d]$ of
nilpotency class $d$ as $L[V]/I_L^{d+1}$ or equivalently as the
Lie subalgebra of $R\la V,d\ra$ generated by the space $V$. These
algebras decompose as
$$
L[V,d] = \bigoplus_{i=1}^d L[V,d]^{(i)} = \bigoplus_{i=1}^d L[V]^{(i)}.
$$
\end{df}

\end{section}

\ifnewpage

\begin{section}{Derivation of free algebras}
\label{subsec:lie:der}

We start with the definition of the derivation of an algebra.

\begin{df}
Let $A$ be an algebra.  An endomorphism $\delta$ of $A$ (as a
linear space) is called a derivation, if for any elements $x,y\in
A$ we have
$$
\delta(xy) =  x\delta(y) + \delta(x)y.
$$
Direct computation shows that the commutator of two derivations is
again a derivation, which gives that the set of all derivations
forms a Lie subalgebra of $\End (A)$, which is denoted as $\D(A)$.
\end{df}

If the algebra $A$ has a grading then this grading induces a
grading on the algebra $\D A$. We say that the derivation $\delta$
is homogenous of degree $k$, if and only if, we have that $\deg
\delta(x) = k + \deg x$ for any homogenous element $x$.

Also, note that if a group $G$ acts on the algebra $A$ by
automorphisms, then this action induces an action on $\End(A)$. It
is easy to see that $\D(A)$ is a $G$ invariant subspace, therefore
we have a natural action of the group $G$ onto the derivation
algebra of $A$.

Finally, note that a derivation $d$ is determined by its values on
the generators of the algebra, because if we know them, using the
Leibnitz rule we can reconstruct the values of $d$ on all
elements. This shows that the map $res_{\gen}:\D A \to \Hom
(\mbox{Gen } A , A)$ is an injection. In general this map is not a
surjection because not all maps can be extended to derivations of
the whole algebra $A$.

There is one very important case where this map is a surjection
--- if the algebra $A$ is the relatively free algebra in certain
variety generated by the space $V$ and the ground field is
infinite, then the map $\res_{\gen}: \D A \to \Hom(V, A)$ is a
surjection. We will use this result only in the cases of free (and
free nilpotent) associative and Lie algebras where it is well
known.

\vspace{3mm}

{\bf Derivations of the free associative algebra}

Let us consider the algebra $\R\la V \ra$ and its derivations ---
since this algebra is free in the class of all associative
algebras, we have that
$$
\D \R\la V \ra \simeq \Hom(V,\R\la V \ra).
$$
The algebra $\R \la V \ra$ has two gradings, which induces the
gradings on the derivation algebra. For a derivation $d$ we will
denote by $\deg d$ and $\mdeg d$ its degree and multi degree,
respectively.

The degree grading gives the following decomposition into
homogenous components
$$
\D \R\la V \ra = \bigoplus_{i=-1}^\infty \D \R \la V \ra^{(i)}
\simeq  \bigoplus_{i=-1}^\infty \Hom (V,\R \la V \ra^{(i+1)}),
$$
where the isomorphism is as linear spaces. Actually this is also a
$\GL(V)$ module isomorphism because the action of $\GL(V)$
preserves the grading by the degree.

\vspace{3mm}

{\bf Derivations of the free Lie algebra}

The algebra $L[V]$ is a Lie subalgebra of $\R \la V \ra$, which
generates it as an associative algebra. This implies that we can
embed the derivation algebra of $L[V]$ into the derivation algebra
of $\R \la V \ra$. Using the decomposition of $\D \R \la V \ra$
into homogenous components we can see that

\begin{lm}
The algebra $\D L[V]$ has the following decomposition into
homogenous components
$$
\D L[V]= \bigoplus_{i=0}^{\infty} \D L[V]^{(i)},
$$
where $\D L[V]^{(0)} = \gll (V)$ and $\D L[V]^{(i)} \simeq
\Hom(V,L[V]^{(i+1)})$ (as a $\GL(V)$ module). There is no
homogenous component of degree $-1$ as in the case of free
associative algebras, because the Lie algebra $L[V]$ does not have
homogenous component of degree $0$.
\end{lm}

Let us look closely at the Lie algebra structure of $\D L[V]$ with respect
to this decomposition.
\begin{reme}
The Lie algebra structure on $\D L[V]$ is given by:

The adjoint action of $\gll (V)$ on
$\D L[V]^{(k)} = \Hom(V,L[V]^{(k+1)})$,
coincides with the adjoint action of the
Lie algebra of the group $GL(V)$ and it comes from the natural
action of $\gll (V)$ on $V$ and $L[V]^{(k+1)}$;

The element $1\in \gll (V)$ acts as the
degree derivation of $\D L[V]$;

The Lie bracket of the elements $f$ and $g$, such that
$f \in \Hom(V,L[V]^{(k+1)})$ and
$g \in \Hom(V,L[V]^{(l+1)})$ is defined to satisfy
$ D_{[f,g]} = [ D_f, D_g]$.
This is equivalent to
$$
[f,g](x) = D_f (g(x)) - D_g(f(x)),
$$
for any $x\in V$.
\end{reme}

The action of $\GL(V)$ on the algebra $\D L[V]$ is not polynomial
(in general), because the action of $\GL(V)$ on the dual space
$V^*$ is not. However if we restrict the action to the subgroup
$\SL(V)$, it becomes polynomial. 
\begin{df}
A poly homogeneous derivation $\delta$ of the algebra $L[V]$ is
called positive
$\mdeg \delta \geq 0$, 
 i.e., we have that
$\deg_{x_i} \delta \geq 0$ for all $i$
.

We will call a derivation $\delta$ totaly positive if the $\GL(V)$ module
generated by $\delta$ contains a basis of positive derivations.
\end{df}

The following lemma gives a more natural classification of the
totally positive derivations.

\begin{lm}
The derivation $\delta$ is totally positive if and only if the
action of $\GL(V)$ on the submodule generated by $\delta$ is
polynomial.
\end{lm}

\begin{proof}
We will not prove this lemma directly here. The statement of
the lemma follows easily from Theorem~\ref{com:split}.
\end{proof}

Let $UL[V]$ denote the `nilpotent' radical
of $\D L[V]$, i.e.,
$$
UL[V] = \bigoplus_{k=1}^{\infty} \D L[V]^{(k)}.
$$

\vspace{3mm}

{\bf Derivations of the free nilpotent Lie algebras}

The derivation algebra of $L[V,n]$ has a similar description:

\begin{lm}
The derivation algebra of $L[V,n]$ has the following decomposition
$$
\D (L[V,d])= \bigoplus_{k=0}^{d-1} \D L[V]^{(k)}.
$$
\end{lm}

Let $UL[V,d]$ denote the nilpotent radical of $\D (L[V,d])$, i.e.
$$
UL[V,d] = \bigoplus_{k=1}^{d-1} \D L[V]^{(k)}.
$$

\begin{reme}
The center of the Lie algebra $UL[V,d]$ is equal to $\D L[V]^{d-1}$.
$UL[V,d]$ is also a nilpotent Lie algebra of class $d-1$.
\end{reme}
\begin{proof}
For any $f \in L[V]$ we can define the inner derivation $\ad(f)$
as $\ad (f)(g) = [f,g]$. Direct computations gives that for any
derivation $\delta$ of the algebra $L[V]$ or $L[V,d]$,
 we have $[\delta ,\ad x ] = \ad(\delta(x))$.
This computation shows that, if a derivation $\delta$ lies in the
center of the algebra $UL[V,d]$ then we have $\delta(x) \in \Cen
(L[V,d])$ for any $x \in L[V,d]$. By taking for $x$ the generators
of the algebra $L[V,d]$ we see that $\deg \delta =d-1$, i.e. $\Cen
UL[V,d] \subset  \D L[V]^{d-1}$. The other inclusion follows from
trivial degree arguments. The nilpotent class is at least $d-1$
because the image of the map $\ad$ is isomorphic to
$L[v,d]/\Cen(L[V,d]) = L[V,d-1]$, which is a nilpotent algebra of
class $d-1$. The nilpotency class is at most $d-1$ by degree
arguments.
\end{proof}

In the next chapters we will see that many properties of the
automorphism tower of the free nilpotent group depend on the
generating set of the algebra $UL[V,d]$.

\end{section}

\ifnewpage

\begin{section}{Some computations in Lie algebra $UL[V]$}
\label{subsec:lie:commutators}

In this section we define notations for some elements in the
algebra $UL[V]$ which will be used later. We also compute the Lie
bracket between these elements.

\begin{df}
Let $h$ be an element in $\R \la V\ra$. Let us define
a derivation $\AD_h \in \D \R \la V \ra$ by
$$
\AD_h(x) = (h_{|x_i \to \ad (x_i)})(x)
$$
for any $x\in V$ and extend it to the whole algebra $R\la V
\ra$ by the Leibnitz rule. By construction we have that
$\AD_h(x) \in L[V]$ for all $x\in V$, therefore we can consider $\AD_h$
also as a derivation of the Lie
algebra $L[V]$. Note that if $h\in L[V] \subset \R \la V \ra$,
then $\AD_h = \ad(h)$
is the inner derivation obtained by the element $h$.

Finally, we can view $\AD$ as a linear map from $\R\la V \ra$ to
 $UL[V] \subset \D L[V]$. We will call the derivations in
the image of the map $\AD$ associative derivations.
\end{df}

\begin{reme}
The map $\AD$ is not injective. For example if $n=2$ then
$$
\AD_{[x_1,x_2]x_1}=0,
$$
and if $n=3$, then
$$
\AD_{[[x_1,x_2],[x_1,x_3]][x_1,x_2]x_1}=0.
$$
Note that if we restrict the map $\AD$ to the elements of degree at most
$n$ it becomes injective, but there are homogeneous elements of degree
$n+1$ in the kernel of $\AD$. In the case of infinitely many variables
($n=\infty$), the map $\AD$ is injective.
\end{reme}
\begin{proof}
Using Shirshov's basis of the free Lie algebra it can be seen that
if $f$ does not depend on $x$ then $\AD_f(x)\not = 0$. This shows
that the kernel of the map $\AD$ is trivial in the case of
infinitely many variables and in the case of $n$ variables there
are no elements in the kernel of degree less than $n$. A more
complicated argument based on the fact that in $L[V]$ there are no
$\GL(V)$ modules corresponding to the partition $[1^n]$, gives
that there are also no elements of degree $n$ in $\ker \AD$.

We can construct an element in $\ker \AD$ of degree $n+1$ using
the following idea --- define a non trivial element $f \in L[x_i,
x_{n+1}]$ by
$$
f(x_i;x_{n+1}) = \sum_{\sigma\in S_{n+1}}
[[x_1,x_{\sigma_1},\dots, x_{\sigma{n+1}}]
$$
Using Shirshov's basis we can rewrite $f$ in such a way that all
commutators start with $x_{n+1}$. Therefore, there exists $h\in
R\la V \ra$ of degree $n+1$ such that $f(x_i;x_{n+1}) = \AD_h
x_{n+1}$. By construction, $f(x_i;x_j)=0$ for $j=1,\dots,n$.
Therefore, $\AD_h$ acts trivially on the space $V$.
\end{proof}
\begin{reme}
The map $\AD$ is surjective only if $n=2$. For example, the
derivation $\delta$ defined by $\delta(x_1)=[x_2,x_3]$ and
$\delta(x_i)=0$ for $i\not=1$ does not lie in $\im(\AD)$.
\end{reme}

\begin{lm}
\label{com:com}
In the Lie algebra $UL[V]$ we have the following equalities:
\begin{enumerate}
\item
$[g,\AD_f] = \AD_{g(f)}$ for any $g \in \gll$ and
$f\in \R \la V \ra$, i.e. $\AD$ is a map of $\gll$ modules.

\item
$[\AD_f,\AD_g] = \AD_{[g,f]}$
for any $f,g \in L[V]$.

\item
$[\delta,\AD_f] = \AD_{\delta(f)}$
for any $f \in L[V]$ and $\delta\in UL[V]$.

\item
$[\AD_{h_1},\AD_{h_2}] = \AD_h$
where $h=-[h_1,h_2] + \AD_{h_1}(h_2) - \AD_{h_2}(h_1)$.
\end{enumerate}
\end{lm}
\begin{proof}
The proofs of all these equalities are similar,
so here we only prove the third and fourth.
$$
[\delta,\AD_f](x) = \delta(\AD_f(x)) - \AD_f(\delta(x))=
\delta([x,f]) - [\delta(x),f] =
$$
$$
=[\delta(x),f] + [x,\delta(f)] - [\delta(x),f] = [x,\delta(f)]=
\AD_{\delta(f)}(x).
$$
Similarly,
$$
[\AD_{h_1},\AD_{h_2}](x) =
\AD_{h_1}(\AD_{h_2}(x)) - \AD_{h_2}(\AD_{h_1}(x))=
$$
$$
\AD_{h_2h_1}(x) + \AD_{\AD_{h_1}(h_2)}(x) -
\AD_{h_1h_2}(x) - \AD_{\AD_{h_2}(h_1)}(x) =
\AD_H(x),
$$
because
$\AD_{f}(\AD_{g}(x))= \AD_{gf}(x) + \AD_{\AD_{f}(g)}(x)$.
\end{proof}

\begin{co}
The space $I_L=\{\AD_f | f\in L[V] \}$ is an ideal in the algebra
$UL[V]$ (the ideal $I_L$ coincides with the space of all inner derivations of
the algebra $L[V]$) and $I_A =\im(\AD)$ 
is a subalgebra. However, by $(4)$ from the previous lemma it can
be seen that it is not a Lie algebra homomorphism.
\end{co}

\begin{df}
Let us define a derivation
$\ADD{x_i}{h}$ by
$$
\ADD{x_i}{h}(x_i) = \AD_h(x_i)
\qquad \ADD{x_i}{h}(x_j) = 0, \mbox{ for } i \not= j.
$$
We will also need the derivation $\DEL{x_i}{f}$
for $f\in L[V]$ defined by
$$
\DEL{x_i}{f}(x_i) = f
\qquad \DEL{x_i}{f}(x_j) = 0, \mbox{ for } i \not= j.
$$
\end{df}

The next lemma gives the commutators in $UL$ of the elements
$\AD_f$ and $\DEL{x}{f}$.

\begin{lm}
\label{com:com1}
In the Lie algebra $UL[V]$ we have the following equalities:
\begin{enumerate}
\item
$[\delta,\ADD{x_i}{h}] = \ADD{x_i}{\delta(h)}$
for any $\delta(x_i)=0$ and
$\deg_{x_i} \delta(x_j)=0$.

\item
$[\ADD{x_i}{h},\ADD{x_j}{g}] =
\ADD{x_i}{\ADD{x_j}{g}(h)} - \ADD{x_j}{\ADD{x_i}{h}(g)}$
if $i\not = j$.

\item
$[\DEL{x_i}{\ad^l(x_k)(x_j)},\DEL{x_j}{\ad^m(x_k)(x_i)}]=
\ADD{x_i}{x_k^{l+m}} - \ADD{x_j}{x_k^{l+m}}$
if $i,j,k$ are pairwise different indexes.
\end{enumerate}
\end{lm}
\begin{proof}
(1) Let us see how $[\delta,\ADD{x_i}{h}]$ acts on an element
$x_j$. There are two possibilities $j=i$
$$
[\delta,\ADD{x_i}{h}](x_i) = \delta(\ADD{x_i}{h}(x_i)) - \ADD{x_i}{h}(\delta(x_i)) =
\delta(h) - 0 = \ADD{x_i}{\delta(h)}(x_i)
$$
and $j \not = i$
$$
[\delta,\ADD{x_i}{h}](x_j) = \delta(\ADD{x_i}{h}(x_j)) - \ADD{x_i}{h}(\delta(x_j)) =
0 - \ADD{x_i}{h}(\delta(x_j)) = 0 = \ADD{x_i}{\delta(h)}(x_j).
$$
This shows that the left and the right side act in the same way on
the generators of the algebra $L[V]$, therefore they coincide as
an element in $\D L[V]$.

The proofs of (2) and (3) are similar.
\end{proof}

The next theorem gives that $UL[V]$ is generated as a linear space
by the image of the map $\AD$ and by the $\GL(V)$ module generated
by poly homogeneous derivations which are not positive.
\begin{thq}
\label{com:split}
The algebra $UL(n,d)$ is generated as an $\GL$ (or $\gll$) module
by the image of the map
$\AD$ and elements $\DEL{x_n}{f}$, where $\deg_{x_n} f =0$.
\end{thq}
\begin{proof}
Since the action of the group $\GL$ on $UL$ is semi simple, it is enough to
prove that any simple sub module is either in the image of the map $\AD$ or
generated by $\DEL{x_n}{f}$, where $\deg_{x_n} f =0$. Consider a new action of
$\GL$ (or $\gll$) on $UL$ defined by
$$
\begin{array}{ll}
g\circ\delta = g(\delta) \det g & \mbox{ for } g\in \GL \\
g\circ\delta = g(\delta) + \tr g . \delta  & \mbox{ for } g\in \gll.
\end{array}
$$
This makes the representation polynomial and if we have a simple
sub module $W$ of $UL$ we can find a highest weight vector
$\delta$ in $W$ (with respect to standard Cartan sub algebra).

By the choice of the element $\delta$, there exists a nonnegative
integer $s$ such that  $e_{n,n}\circ \delta = s \delta$, here
$e_{i,j}$ denote the matrix which have zeros everywhere except for
the $i,j$ place where it has a one.

Now there are two possibilities --- either $s=0$ or $s\geq 0$.

If $s=0$ then $e_{n,n}$ acts trivially on $\delta(x_n)$, therefore $\delta(x_n)$ is
a Lie polynomial on the letters $x_1,\dots, x_{n-1}$, i.e.,
$\deg_{x_n} \delta(x_n) = 0$.
The element $\delta$ is highest weight vector,
therefore $e_{n,i} \delta =0$ for $i \leq n$,
which implies that
$$
\delta(x_i) = e_{n,i} (\delta(x_n)) =0.
$$
The last equation is the same as $ \delta = \DEL{x_n}{\delta(x_n)}$.

If $s>0$, then the standard action of $\GL$ on $W$ is also
polynomial. We will use induction on the maximal $l$ such that
$\delta(x_{n-k})=0$ for all $k \geq l$, to prove that $\delta$
lies in the image of the map $\AD$. The base case is trivial since
if $l=0$ then $\delta =0$.


Let us linearize $\delta(x_k)$ with respect to the variable $x_k$,
i.e., let us substitute $x_k + y$ in the place of $x_k$ in
$\delta(x_k)$ and expand the result as a function on the variable
$y$:
$$
\delta(x_k)|_{x_k\to x_k + y} = \delta(x_k) + \sum F(\delta)_k^i(x_1,\dots,x_n;y),
$$
where $F_k^i$ is a Lie polynomial with degree $i$ with respect to $y$.
It follows from the definition of $F_k^1$ that
$$
F(\delta)_k^1(x_1,\dots,x_n;x_i) = e_{k,i} (\delta(x_k)),
$$
and that $F_k^1$ lies in the image of the map $\AD$.
Now let us define the derivation $F$
by
$$
F(\delta)(y) = \sum_k F(\delta)_k^1(x_1,\dots,x_n;y),\,\,\mbox{for }y \in V.
$$

First, we want to show that the element $F$ is also a highest
weight vector for some $\gll$ module isomorphic to $W$, which
follows from the next lemma.

\begin{lm}
Let the element $F(\delta)$ be defined in the way above. Then, for
$g\in \gll$ we have $g F(\delta) = F(g\delta)$.
\end{lm}
\begin{proof}
In order to prove the lemma it is enough to show it for $g=e_{p,q}$. Let us
compute how $e_{p,q}$ acts on $F(\delta)(x_i)$.
$$
\begin{array}{rl}
e_{p,q}(F(\delta)(x_i)) &= \sum_k e_{p,q}(e_{k,i}(\delta(x_k))) = \\
& = \sum_k [e_{p,q}e_{k,i}](\delta(x_k)) + \sum_k e_{k,i}(e_{p,q}(\delta(x_k))) =\\
& = \delta_{p,i} \sum_k e_{k,q} (\delta(x_k)) - e_{p,i} \delta (x_q) +\\
& \mbox{\hspace{1cm}}\sum_k e_{k,i}((e_{p,q} \delta)(x_k)) + e_{p,i} \delta(x_q) =\\
& = F(\delta)(e_{p,q}x_i) + F(e_{p,q}\delta),
\end{array}
$$
where $\delta_{p,i}$ is the Kronecker delta. Therefore,
$e_{p,q}(F(\delta)) = F(e_{p,q}\delta)$, which proves the lemma.
\end{proof}

Let us compute $F(x_i)$
$$
F(x_i) = \sum e_{k,i}\delta(x_k) = \sum_{k\leq i} e_{k,i}\delta(x_k) + \sum_{k>i} e_{k,i}\delta(x_k)=
\sum_{k\leq i} e_{k,i}\delta(x_k) + (n-i) \delta(x_i),
$$
i.e., $F(\delta) (x_{n-k}) = 0$ for $k \geq l$ and
$$
F(\delta)(x_m) = ((n-m) \delta + e_{m,m} \delta)(x_m) = s_l \delta(x_m)
$$
for $m= n-l-1$ and some positive integer $s_l\geq s>0$. Therefore,
we may use the induction assumption of the derivation $\delta -
s_l^{-1} F$, i.e., we may assume that it lies in the image of
$\AD$. The derivation $s_l^{-1}F$ also lies in the image of $\AD$,
therefore the same is true for $\delta$. This finishes the proof
of the induction step and completes the proof of the theorem.

\end{proof}

\end{section}

\ifnewpage

\begin{section}{Generators of the algebra $UL[V]$}
\label{subsec:lie:red}

The $\sll$-invariant component of the set of generators of the
algebra $UL[V]$ has direct connections with Kazhdan property $T$
of the automorphism group of the free nilpotent group. In this
section our aim is to prove Theorem~\ref{gen:maincom}, which fully
describes the generators of $UL[V,d]$ if $d\leq n(n-1)$.

\begin{lm}
The algebra $\D L[V]$ is generated as a Lie algebra by $\gll$ and
the elements which generate $UL[V]/\big[UL[V],UL[V]\big]$ as
$\gll$ module.
\end{lm}
\begin{proof}
For any $d$, the algebra $UL[V,d]$ is a nilpotent algebra.
Therefore, it is generated as a Lie algebra by the elements in
$UL[V,d]/\big[UL[V,d],UL[V,d]\big]$. Using limit arguments and the
fact that
$$
\bigcap_d \ker (UL[V] \to UL[V,d]) =0
$$
we see that $UL[V]$ behaves like a nilpotent Lie algebra.
Therefore, it is generated by the elements which generate
$UL[V]/\big[UL[V],UL[V]\big]$ as $\gll$ module.
\end{proof}

The algebra $\sll$ is simple. Therefore, the $\sll$ module
$\big[UL[V],UL[V]\big]$ has a (non invariant) complement
$U_{\gen}[V]$, which is also a graded $\sll$ module. Thus, we have
the isomorphism
$$
UL[V]/\big[UL[V],UL[V]\big] \simeq U_{\gen}[V].
$$
\begin{reme}
This module is not uniquely defined in general but its isomorphism
type does not depend on the choice of $U_{\gen}$ (if $d$ is small
compared to $n$ it follows from Theorem~\ref{gen:maincom} that the
homogenous components of $U_\gen [V]$ of degree less than $d$ are
uniquely defined). In this section we want to treat the elements
of $U_{\gen}$ as elements in $UL$, that is why we define
$U_{\gen}$ as a submodule. In general it is better to use the
above isomorphism to define $U_{\gen}$ as an abstract $\gll$
module (see Chapter 3).
\end{reme}

Similarly, we can define $U_{\gen}[V,d]$, which satisfies
$$
U_{\gen}[V,d] = U_{\gen}[V] \cap UL[V,d].
$$

\begin{lm}
\label{gen:negdeg}
Let $\delta$ be a derivation such that
$\deg_{x_i} \delta= -1$ for some index $i$ and $\deg \delta >1$.
Then $\delta$ lies in
the commutator subalgebra of $UL[V]$.
\end{lm}
\begin{proof}
If $\deg_{x_i} \delta= -1$, then $\delta= \DEL{x_i}{f}$, for some $f$ with $\deg_{x_i} f =0$.
Using Lemma~\ref{com:com1} we can see that
$$
[\DEL{y_i}{h},\DEL{x}{[y_i,y_j]}] = \DEL{x}{[h,y_j]},
$$
for any $i\not=j$ and $h\in L[V]$ such that $\deg h >1$, but $\deg_x h =0$. Any
element $f\in L[V]$, which does not depend on $x$, can be expressed
as a linear combination of
elements $[h,y_j]$. This, together with the observation above,
proves the lemma. This proof does nor work when $n=2$ (we can not find two different
$y$-es), but in that case there are no derivations $\delta$ such that
$\deg_{x_i} \delta = -1$, since the map $\AD$ is surjective.
\end{proof}

This result together with lemma~\ref{com:split} gives that
the module $U_{\gen}$ `almost' lies in the image of the map $\AD$:

\begin{thq}
We have that
$$
UL[V]^{(1)} \subset U_{\gen}[V] \subset \im (\AD) + UL[V]^{(1)}.
$$
\end{thq}

\begin{reme}
The elements $\AD_{x_i^k}$ for $k \geq 1$ do not lie in
$\big[UL[V],UL[V]\big]$, i.e., they lie in $U_{\gen}[V]$.
\end{reme}
\begin{proof}
Using multi degree arguments we can see that the only way to
express a derivation of multi degree $(k,0,\dots,0)$ as a
commutator of two derivations is if they have multi degrees
$(l,0,\dots,0)$ and $(k-l,0,\dots,0)$ or $(l,1,-1,0,\dots,0)$ and
$(k-l,-1,1,0,\dots,0)$. Using the computations from
Lemmas~\ref{com:com} and~\ref{com:com1}, we can see that in both
cases the commutator of these derivations has trivial projection
in the submodule generated by $\AD_{x_i^k}$.
\end{proof}

In order to show that these elements
generate $U_{\gen}$ for small values of $d$, we need to prove several
lemmas which allow us to `simplify' elements in $UL$ modulo
$\big[UL[V],UL[V]\big]$.


From now on, we need to distinguish one element of the set
$X=\{x_1,\dots,x_n\}$ from the others. Therefore, we will split
$X$ as $X = \{x\} \cup Y$, where $Y = \{ y_1,\dots, y_l\}$, $l=n-1$, and $x
\not\in Y$.

First we will prove a lemma about free associative
algebras, which will be used later in lemma~\ref{gen:onevariable} to
simplify the elements in $UL$ modulo the commutator subalgebra.

\begin{lm}
\label{gen:simplifyass}
Let $f\in \R \la x, Y\ra$ be an element which lies in the
commutator ideal of $\R\la x,Y\ra$ (i.e., $f$ has a trivial
projection in the polynomial ring $\R[x,Y]$). Then $f$ can be
expressed as a linear
combination of elements of the type $\DEL{x}{g} h$, where
$h \in \R \la x, Y\ra$ and $\DEL{x}{g}$ is a derivation of $\R \la x,Y\ra$,
which kills all $y_i$-es and sends $x$ to a homogeneous element $g$
in the free Lie algebra $L[x, Y]$ of degree at least $2$.
\end{lm}
\begin{proof}
The free associative algebra $\R \la x, Y\ra$ is the universal
enveloping algebra of the free Lie algebra $L[x, Y]$. Therefore,
any element $f$ in $\R \la x,Y \ra$ can be expressed as a linear
combination of the elements
$$
\Sym_k(h_1,\dots, h_k)=
\sum_{\sigma \in S_k} h_{\sigma(1)} \dots h_{\sigma(k)} ,
$$
where the elements $h_i \in L[x, Y]$ are homogeneous, and $k \leq
p=\deg f$.

The condition that $f$ has a trivial projection in the polynomial
algebra $\R[x,Y]$, implies that in the above decomposition there
are no elements with $k=p$. This gives that at least one of the
elements $h_i$ has degree at least $2$. Therefore, it suffices to
prove the lemma only for $f=\Sym_k(h_1,\dots, h_k)$. Without loss
of generality, we may assume that if $\deg h_i=1$, then either
$h_i=x$ or $h_i \in Y$.

We will prove the lemma for these $f$'s
(assuming that $k$ is fixed)
 using induction on the
number $q$ of elements  $h_i$ of degree greater then $1$. The base
case $q=0$ is impossible, because $k < p$.

Let $q>0$.
We can permute the elements $h_i$, in order to have
$$
\deg h_1 \geq \dots \geq \deg h_q > 1
$$
and $ h_{q+1}= \dots =h_{q+t}=x$ and the
remaining $h$'s lie in $Y$.
In this case we have the following equality:
$$
\DEL{x}{h_1}(\Sym_p(x,h_2,\dots,h_p))=
$$
$$
=(t+1) \Sym_p(h_1,\dots,h_p) +
\sum_{i=2}^q \Sym_p (x,h_2,\dots,\DEL{x}{h_1}(h_i),\dots h_p).
$$
In each term of the sum we have that the number of Lie elements of
degree at least $2$ is $q-1$. Then, by the induction hypothesis,
$\Sym_k(h_1,\dots, h_k)$ is expressible by a linear combination of the
elements $D_g(h)$. This, together with the fact that the left hand
side is of the same type, proves the induction step.
\end{proof}

\begin{lm}
\label{gen:onevariable}
Let $f\in L[x,Y]$ be an element such that
$1\leq\deg_x f \leq l = |Y|$.
Then the derivation $\DEL{x}{f}$ is equivalent
modulo $\big[UL[x,Y],UL[x,Y]\big]$
to a derivation
lying in the $\sll$ module generated by
$\Delta_{s;n_1,\dots, n_s}$, for $s\leq l$. Here we have denoted
$$
\Delta_{s;n_1,\dots, n_s}=
\ADD{x}{y_1^{n_1} x y_2^{n_2} x \dots x y_s^{n_s}}.
$$
\end{lm}
\begin{proof}
Any such element $f$ can be expressed as a linear combination of
elements
$$
\left({h_1}_{|y_i\to \ad(y_i)} \circ \ad(x) \circ
h_2 \circ \ad(x) \circ \dots \circ h_s \right)(x),
$$
where $s=\deg_x f $ and $h_i$ are associative polynomials on $y_i$-es,
and we have substituted $y_i$ with $\ad(y_i)$ in each of them.
Denote the derivation which kills all $y_i$-es
and sends $x$ to the element above by $\DEL{x}{s;h_1,\dots,h_s}$.

Using Lemma~\ref{com:com} part 5, we can show that
$$
\sum_{j=1}^s \DEL{x}{s;h_1,\dots,\nu(h_j),\dots,h_s} \in \big[UL[x,Y],UL[x,Y]\big],
$$
for any derivation $\nu$ such that $\nu(x)=0$ and
$\nu(y_i)\in L[Y]$.

Let us fix $s\leq l$.
We will show that the derivation $\DEL{x}{s;h_1,\dots,h_s}$ is
equivalent to a derivation which
lies in the $\sll$ module generated by $\Delta_{s;n_1,\dots, n_s}$. To
show this we will use induction on $(\deg h_1, \dots,\deg h_s)$ with respect to
lexicographical order.

Note that the base case is trivial because
$\DEL{x}{s;1,\dots,1}=0$.
On each step we will use induction on $1\leq k \leq s$
to show that such $\DEL{x}{s;h_1,\dots,h_s}$ lies in
$\sll$ the module generated by
$$
\DEL{x}{s;y_1^{n_1},\dots,y_k^{n_k},h_{k+1},\dots ,h_s}.
$$
Thus, the second induction allows us to make the induction step
for the first one.

Therefore, the only thing left is to show the induction step for
the second induction (the second induction also has a trivial base
case).

Let us first consider the case where $h_{k+1}$ lies in the
commutator algebra of $\R\la Y\ra$. After rearranging the variables
as follows: $x:=y_{k+1}$; $y_i:=y_i$ for $i\leq k$, and $y_i:=y_{i+1}$ for
$i > k$ we can apply Lemma~\ref{gen:simplifyass} and, without loss
of generality, we can assume that
$h_{k+1} = D_{|y_{k+1}\to g} t$ for some Lie polynomial $g$ and an
associative polynomial $t$. Using Lemma~\ref{com:com} we
notice that
$$
[D_{|y_{k+1}\to g},
\Delta_{s;y_1^{n_1},\dots,y_k^{n_k},t,h_{k+2},\dots ,h_s}](x)=
$$
$$
=\Delta_{s;y_1^{n_1},\dots,y_k^{n_k},h_{k+1},\dots ,h_s} +
\sum_{j=k+2}^s \Delta_{s;y_1^{n_1},\dots,y_k^{n_k},t,h_{k+2},\dots,
D_{|y_{k+1}\to g} h_j,\dots ,h_s}.
$$
All terms in the sum have multi-degree less than $(n_1,\dots,n_s)$,
because at position $k+1$ the degree is
$$
\deg t = \deg h_{k+1} + 1 - \deg g < \deg h_{k+1} = n_{k+1},
$$
and for all $j\leq k$ the degree is $n_j$. Therefore, by the
induction hypothesis they are equivalent to elements in the $\sll$
module generated by $\Delta_{s;n_1,\dots,n_s}$. Therefore, the
same is true for
$\DEL{x}{s;y_1^{n_1},\dots,y_k^{n_k},h_{k+1},\dots ,h_s}$.

Using this argument, we can view $h_{k+1}$ as an element (up to
equivalence) in the polynomial algebra $\R[Y]$.  Now let us consider
the Lie algebra
$$
\mathfrak{u} =\{ t \in \sll(Y) \mid t(y_i)=0, \mb{ for }j\not= k+1\}
$$
and its action on the set of derivations of the form
$\DEL{x}{s;y_1^{n_1},\dots,y_k^{n_k},h_{k+1},\dots ,h_s}$. This
set is closed under the action of the algebra $\mathfrak{u}$,
because all elements in the algebra $\mathfrak{u}$ kill the
letters $y_i$, for $i\leq k$, and $x$. The elements $t_{k+1}^l$
generate the $\R[Y]$ as a $\mathfrak{u}$ module. Therefore, any
derivation $\delta
=\DEL{x}{s;y_1^{n_1},\dots,y_k^{n_k},h_{k+1},\dots ,h_s}$ up to
equivalence lies in the $\sll$ module generated by the elements
$\DEL{x}{s;y_1^{n_1},\dots,y_{k+1}^{n_{k+1}},h_{k+2},\dots,h_s}$.
This argument completes the induction step.

Thus, we have shown that up to an element in $[UL,UL]$, every
derivation $\DEL{x}{f}$ with $\deg_x f \leq l$ lies in the $\sll$
module generated by $\Delta_{s;n_1,\dots,n_s}$.
\end{proof}

\begin{lm}
If $s\geq 2$ then the derivation $\Delta_{s;n_1,\dots,n_s}$ lies
in the $\sll$ module generated by $\Delta_{1;y_1^{n}}$, for
$n=\sum n_i+s-1$, plus the commutator subalgebra of $UL$.
\end{lm}
\begin{proof}
Let us compute the commutator of derivations
$\Delta_{s-1;y_1^{n_1},\dots,y_{s-2}^{n_{s-2}},y_{s-1}^{n_s}y_s}$
and $\DEL{y_s}{\ad^{n_s}(y_s)(x)}$.
$$
[\Delta_{s-1;y_1^{n_1},\dots,y_{s-2}^{n_{s-2}},y_{s-1}^{n_s}y_s},
\DEL{y_s}{\ad^{n_s}(y_s)(x)}](x) =
$$
$$
=-\DEL{y_s}{\ad^{n_s}(y_s)(x)}
\left(\ad^{n_1}(y_1)\ad(x) \dots \ad(x) \ad^{n_{s-1}}
(y_{s-1})[x,y_s]\right)=
$$
$$
=\left(\ad^{n_1}(y_1)\ad(x) \dots \ad(x)
\ad^{n_{s-1}}(y_{s-1})\ad(x)\ad^{n_s}(y_s)(x) \right)
$$
and
$$
[\Delta_{\dots},\DEL{y_s}{\ad^{n_s}(y_s)(x)}](y_s) =
\Delta_{\dots}\left(\ad^{n_s}(y_s)(x)\right) =
$$
$$
=\left(\ad^{n_s}(y_s)\ad^{n_1}(y_1)
\ad(x) \dots\ad(x)\ad(y_{s-1})^{n_{s-1}}\ad(x)\right)(y_s).
$$

Using a permutation (the one which interchanges $x$ and $y_s$) of
variables we can see that the last expression is equivalent to one
in the $\sll$ module generated by
$$
\Delta_{n_s+1;1,\dots,1,y_s y_1^{n_1}
y_s \dots y_s y_{s-1}^{n_{s-1}}y_s}.
$$
Using the same argument as in the previous lemma we can show that
the last derivation lies in the $\sll$ module generated by
$\Delta_{n_s+1;1,\dots,1,y_1^{n'}}$, where $n'=\sum_{i=1}^{s-1}
n_i + s-1$.

Now using the commutator of derivations $\DEL{x}{\ad^{n_s}(x)(y_2)}$
and $\DEL{y_2}{\ad^{n'}(y_1)(x)}$, we can show that
$\Delta_{n_s+1;1,\dots,1,y_1^{n'}} \sim
\DEL{y_2}{\ad^{n'}(y_1)\ad^{n_s}(x)(y_2)}$.

By the argument from Lemma~\ref{gen:onevariable}, we can see that the
last derivation is equivalent to a derivation in the $\sll$-module
generated by $\Delta_{1;y_1^n}$ for
$n = n' + n_s = \sum n_i + s -1$.
\end{proof}

\begin{lm}
\label{gen:actonone}
Let $\delta$ be a derivation such that
$\deg_{y_i}\delta \geq 0$ for all $i$ and $\sum \deg_{y_i}\delta
>1$. Then $\delta$ is equivalent (modulo the commutator subalgebra
of $UL$) to a derivation which acts nontrivially only on the
variable $x$.
\end{lm}
\begin{proof}
Let us first consider the case when $\delta =\DEL{y_i}{[f,y_j]}$ for
some indices $i, j$ and some Lie polynomial $f \in L[x,Y]$. Using
the derivations $\DEL{x}{f}$ and $\DEL{y_i}{[x,y_j]}$,
it can be seen
that $\delta$ is equivalent to $\DEL{x}{\DEL{y_i}{[x,y_j]}(f)}$,
which is a derivation acting nontrivially only on $x$.

Now, let us consider the case when $\delta$ is such that
$\delta=\DEL{y_i}{\ad^n(x)[f,y_j]}$
for some $i,j$ and some Lie polynomial $f \in L[x,Y]$.

This time, using derivations $\DEL{y_i}{\ad^n(x)(y_i)}$ and
$\DEL{y_i}{[f,y_j]}$, we can show that
$$
\delta \sim \DEL{y_i}{[\DEL{y_i}{\ad^n(x)(y_i)}(f),y_j]}.
$$
But the last derivation is of the type considered before.
Therefore, it is equivalent to one acting nontrivially only on
$x$.

Since every derivation satisfying the conditions of the
lemma can be
expressed as a linear combination of
the derivations considered above,
the lemma is proved.
\end{proof}

%

\begin{lm}
\label{gen:smalldeg}
Let $\delta$ be a poly-homogeneous derivation in $UL[x,Y]$,
and $\deg_x \delta <l=|Y|$. Then, if $\deg \delta >1$, then
$\delta$ is equivalent modulo $[UL,UL]$ to a derivation
in  the $\sll$ module
generated by $\AD_{x^n}$, where $\deg \delta =n$.
\end{lm}
\begin{proof}
Suppose that there exists $i$ such that $\deg_{y_i} \delta = -1$.
Then by lemma~\ref{gen:negdeg} $\delta$ lies in the commutator subalgebra
of $UL$.

If we have that $\deg_{y_i} \delta > 0$ for all $i$, from
Lemma~\ref{gen:actonone} it follows easily that
$\delta \sim \delta'$,
where $\delta' = \DEL{x}{f}$ and $\deg_x f \leq l$.
Now by Lemma~\ref{gen:onevariable}, $\delta'$ is equivalent to a
derivation
in a $\sll$ module generated by $\ADD{x}{y_1^n}$. But from
Lemma~\ref{com:com} we have that
$\ADD{y_i}{x^n}\sim
\ADD{y_j}{x^n}$, i.e.
$l. \ADD{x}{y_1^n} \sim  \AD_{y_1^n}$,
which completes the proof.
\end{proof}

\begin{thq}
\label{gen:maincom}
If $d\leq n(n-1)$, then the Lie algebra
$UL(n,d)$
is generated as a Lie algebra and $\sll$ module by the elements:
$\AD_{x_1^k}$, for $1\leq k<d$, and if $n\geq 3$ the element
$\DEL{x_1}{[x_2,x_3]}$.
\end{thq}
\begin{proof}
If $\delta \in UL(n,d)$ is a poly-homogeneous derivation,
then there exists an
$i$ such that $\deg_{x_i} \delta < n-1$. If $\deg \delta > 1$,
we can apply Lemma~\ref{gen:negdeg} or Lemma~\ref{gen:smalldeg},
to show that
$$
\delta \in \spann_{\sll}\{ \AD_{x_1^n} \}
+ \big[UL[X],UL[X]\big].
$$
Using this and the fact that $\D L[X]^{(1)}$ is generated as
$\sll$ module by $\AD_{x_1}$ and $\DEL{x_1}{[x_2,x_3]}$,
we can easily finish the proof.
\end{proof}
\begin{co}
The same argument shows that $UL[x_1,\dots,x_n]$ is generated (as
a Lie algebra and $\sll$ module) by $\AD_{x_1^k}$, for $k\geq 1$,
$\DEL{x_1}{[x_2,x_3]}$ and elements $\AD_h$, where $h\in \R\la X
\ra$ such that $h$ is a generator of a $\gll$ module corresponding
to a partition $\lambda$ with $\lambda_n\geq n-1$.
\end{co}
\begin{proof}
Let $h\in \R \la X \ra$ be an associative polynomial,
which lies in the $\sll$-module $M$ corresponding to the partition
$\mu=(\mu_1,\dots,\mu_n)$, with $\mu_n<n-1$. Then, there is a
generator $\tilde h$ of the module $M$ such that
$\deg_{x_1} \tilde h = \mu_n < n-1$. By Lemma~\ref{gen:smalldeg}
we have that $\tilde h \in \big[UL[X],UL[X]\big]$.
\end{proof}

\end{section}

\ifnewpage

\begin{section}{Conjecture}

We will see that the degree of the first $\sll_n$ invariant
generator of $UL[V]$ determines when the automorphism group of the
free nilpotent group on $n$ generators and nilpotency class $d$
has Kazhdan property $T$, and it also determines the behavior of
the automorphism tower of that group. Therefore, it is useful to
have a notation of that degree

\begin{df}
\label{gen:hes}
Let $h(n)$ be the minimal degree of an $\sll$-invariant
generator of $UL[x_1,\dots,x_n]$, i.e.,
$$
h(n) = \min \{\deg f \mid
f \in U_{\gen}\mbox{ and } [f,\sll]=0\}.
$$
\end{df}

Theorem~\ref{gen:maincom} states that $h(n) \geq n(n-1)$. We were unable
to show any $\sll$ invariant element which does not lie
in $[UL,UL]$, unless $n=2$. Therefore we want to make a
conjecture that Theorem~\ref{gen:maincom} can be
significantly improved.

\begin{con}
\label{gen:conj}
If $n\geq 3$, then the Lie algebra $UL[x_1,\dots,x_n]$ is
generated as a Lie algebra and $\sll$ module by: $\AD_{x_1^n}$
for $n\geq 1$ and $\DEL{x_1}{[x_2,x_3]}$, i.e., $h(n)=\infty$.
\end{con}

If this conjecture is true then the automorphism tower of the free
nilpotent groups on $n$ generators will have relatively simple
description and it will stabilizes after at most $3$ steps.

\begin{reme}
Let $V$ be an infinite dimensional linear space with basis $x_i$,
for $i=1, \dots$. We can define the free Lie algebra $L[V]$ and
its derivation algebra in the same way as in sections~2.1 and~2.2.
Then, the proof of Theorem~\ref{gen:maincom} with slight
modifications gives that $UL[V]$ is generated as a Lie algebra
(allowing certain infinite sums) and $\Aut(V)$ module by the
elements $\AD_{x_1^k}$, for $k=1,\dots, \infty$, and the element
$\DEL{x_1}{[x_2,x_3]}$, i.e., the conjecture~\ref{gen:conj} is
true in the case $n=\infty$.
\end{reme}

\begin{reme}
\label{gen:conj2} The analog of the above conjecture for $n=2$
does not hold. For example, the element $\AD_{[x_1,x_2]^2}$ does
not lie in the commutator subalgebra of $ UL[x_1,x_2]$. This is
true because $\AD_{[x_1,x_2]^2}\not \in I_L$  and there is only
one $\sll_2$ module of degree less then $3$ which is not in $I_L$,
namely the module $U$ generated by $\AD_{x_1^2}$, and direct
computation shows that there are no trivial $\sll_2$ submodules in
$[U,U]$. This argument gives that $h(2)=4$ and suggests that it is
possible to find a counter example to Baumslag's conjecture among
the $2$ generated relatively free groups in some nilpotent
varieties.
\end{reme}

We tried to describe the generators of the algebra $UL(n,d)$ for
some small values of $n$ and $d$ hoping to disprove
Conjecture~\ref{gen:conj}. Some computer simulations verified this
conjecture for $n=3$ and $d \leq 13$ and for $n=4$ and $d \leq
13$. Unfortunately, we were unable to examine the case $n=5$
because by Theorem~\ref{gen:maincom}, we need to consider $d>20$
and in this case the dimension of $UL$ is too big and the brute
force approach does not work, because we do not have sufficient
computer resources.
\end{section}

\end{chapter}
\newpage

\begin{chapter}{Automorphism Group of Free Nilpotent Group}
\thispagestyle{empty}
	In this chapter we study the automorphism groups of
free nilpotent groups. Our main goal is to prove Theorem~\ref{whenT}
which gives partial answer to a
question posed by A. Lubotzky and I. Pak in~\cite{LubPak}.
In the last section of this chapter we state a recent result by
E. Formanek (see~\cite{Formanek})  describing the center of this
automorphism group. This result is later used in Chapter 4, to obtain description
of the automorphism tower of free nilpotent groups.

Our main method in studying $\Aut \Gamma$ is to embed $\Aut \Gamma$ in some Lie group
$G$ over $\R$. Using this embedding we transfer many questions concerning
the group $\Aut \Gamma$ to similar question about the group $G$. After that
we use methods from Lie theory to transfer these problems to the Lie
algebra of this group, which we studied in Chapter 1.

\begin{section}{Free Nilpotent Groups}
\label{subsec:auto:groups}

In this section we describe the free nilpotent group $\Gamma(n,d)$
of nilpotency class $d$ on $n$ generators. We also describe
the free unipotent group $G(n,d)$, which naturally contains
$\Gamma(n,d)$ as a lattice.

\begin{df}
Let $\mc{F}_n$ be the free group on $n$ generators which we will label
$g_1, \dots, g_n$. The elements in this group are all words on $g_i$ and
$g_i^{-1}$ without cancellations. This group is called a free group
because it satisfy the universal property:
\end{df}
\begin{prop}
For any group $G$ and any elements $a_i \in G$
there exist a unique group homomorphism $\phi: F_n \to G$,
such that $\phi(g_i) = a_i$.
\end{prop}

\begin{df}
Let $G$ be a group. We can define the sequence of normal
subgroups $G^{(i)}$ of $G$ called the lower central series of the group $G$ by
$G^{(0)} = G$ and $G^{(i+1)} = [G,G^{(i)}]$, where $[A,B]$ denotes the
subgroup generated by all commutators $[a,b] = aba^{-1}b^{-1}$, for
$a\in A$ and $b \in B$. A group $G$ is called
nilpotent of class $d$ if and only if $G^{(d+1)} = 1$.
\end{df}

Taking a quotient of the free group by a group in the lower central series we
obtain a free nilpotent group.

\begin{df}
Let $\Gamma(n,d)$ be the free nilpotent group on $n$ generators of
nilpotency class $d$, i.e., $\Gamma(n,d)= \mc{F}_n/\mc{F}_n^{(d+1)}$, where
$\mc{F}_n$ is the free group on $n$ generators and $\mc{F}_n^{(d+1)}$ is the
$d+1$ term of its lower center series.
This group is called a free nilpotent group
because it satisfies the universal property:
\end{df}
\begin{prop}
For any nilpotent of class $d$ group $G$ and any elements $a_i \in G$
there exists a unique group homomorphism $\phi: \Gamma(n,d) \to G$,
such that $\phi(g_i) = a_i$.
\end{prop}

We want to describe the automorphism group $\Aut \Gamma(n,d)$. One
way of doing this is to embed the group $\Aut \Gamma(n,d)$ as
a 
lattice into a Lie group (over the reals) $G(n,d)$.
Then, using rigidity results, the group $\Aut \Gamma(n,d)$
can be embedded into $\Aut G(n,d)$. Finally, we use the
properties of the group $\Aut G(n,d)$ to obtain results
about $\Aut \Gamma(n,d)$.

Our first step is to define a unipotent Lie group
$G(n,d)$ such that $\Gamma(n,d)$ is a lattice in it. The
general construction of embedding a discrete nilpotent group into the
unipotent group was first done by P. Hall (see~\cite{Hall}). This embedding is known
as Malcev embedding --- actually the Malcev embedding gives a unipotent group over
$\Q$. In order to obtain an unipotent  group over $\R$ we need to take its
completion. Many results in this chapter as well as in Chapter 4 can be obtained by
working with Lie groups over $\Q$ instead of Lie groups over $\R$. We feel that
working with Lie group over the real numbers is more natural, and this is the main
reason why we use Lie group over $\R$ in these chapters.

\begin{df}
Let $G(n,d)$ be the free unipotent group over $\R$ on $n$
generators $g_1,\dots,g_n$ and of class $d$.
One way to construct the group $G(n,d)$ is to use the fact that it can
be embedded in the multiplicative
group of the algebra $\R.1 + \R\la V,d\ra$. Its image is
the group generated
by $g_i = 1+x_i$ and operation $g \to g^{\alpha}$ for $\alpha \in \R$,
which is defined formally by
$(1+h)^\alpha = \sum_{k=0}^d {\scriptstyle \left( \!\!\begin{array}{c}k\\n \end{array}\!\! \right) }
h^k$ for
any $h\in I$
Here $I$ is the augmentation ideal in $\R \la x_1,\dots,x_n ;d\ra$.
Using this embedding we will often assume that $G(n,d)$ is a subgroup of
the multiplicative group of the algebra $\R.1 +\R\la V,d\ra$.
\end{df}

The following result is well known.
\begin{lm}
\label{gr:aut0}
The group $G(n,d)$ is a Lie group and its Lie
algebra is isomorphic to the free nilpotent Lie algebra on $n$ generators and
class $d$. Also, the group $\Gamma(n,d)$
is a Zariski dense lattice in it.
\end{lm}
\begin{proof}
Using the universal properties of the groups $\Gamma(n,d)$ and $G(n,d)$ it can
be seen that the map $\rho: \Gamma(n,d) \to G(n,d)$ sending $g_i$ to $1+x_i$ is
an injection, i.e., the group $\Gamma(n,d)$ can be embedded into $G(n,d)$.
The image of this embedding is Zariski dense because the quotient
$G(n,d)/\Gamma(n,d)$ is compact.
\end{proof}

Using this lemma from now on we will assume that $\Gamma(n,d)$ is a
subgroup of $G(n,d)$.

\end{section}

\ifnewpage

\begin{section}{Automorphism group}

Now we are ready to define the automorphism group of $\Gamma$, which is the
main object in this chapter.
\begin{df}
Let $G$ be a group, by
$\End(G)$ we will denote the set of all endomorphisms of the group $G$.
Let us note that every endomorphism is uniquely determined by
its values on the generators of the group $G$, i.e., the restriction map
$\res_{\gen}: \End G \to \Hom( \mbox{Gen }G ,G)$
is an injection (here $\Hom$ denotes just the
set of maps from $\gen G$ to $G$).
As in the case of algebras this map in not a surjection in general
but is a surjection if and only if $G$ is the relatively free group
in some variety.

Let also denote by $\Aut G  \subset \End(G)$ the subset of all invertible
endomorphisms of the group $G$. By definition
$\Aut G$ is a group, which is called the automorphism group of the group $G$.
It contains a normal subgroup consisting of all inner automorphisms
$\Inn (G) = \{ \ad g | g \in G\} < \Aut (G)$, where $\ad g$ denotes
the conjugation by the element $g$, i.e., $(\ad g)(h) = ghg^{-1}$.
\end{df}

Using the universal property of the group $G(n,d)$ it can be seen that
every map from
the generating set of $G$ to $G$ itself, can be extended to an endomorphism of $G$, i.e.,
$\End G(n,d)$ is isomorphic to $G(n,d)^{\times n}$ as a set. The
automorphism group contains all the invertible endomorphisms --- these are the
endomorphisms which act on the abelinization as invertible linear
transformations.

\begin{lm}
The automorphism group of $G(n,d)$
satisfies the exact sequence
$$
1 \to UG(n,d) \to \Aut G(n,d) \to \GL_n(\R) \to 1
$$
where $UG(n,d)$ is a unipotent group of class $d-1$ consisting of
all automorphisms which act trivially on the abelinization
$G/[G,G] \simeq \R^n$. The group $UG(n,d)$ is isomorphic to
$\Hom(G/[G,G],[G,G])$ as a topological space.

In fact the exact sequence splits and we have the isomorphism
$$
\Aut G(n,d) = \GL_n(\R) \ltimes UG(n,d).
$$
\end{lm}
\begin{proof}
The group $UG$ is unipotent because it can be seen that
$UG^{(i+1)}$ acts trivially on $G/G^{(i)}$. Therefore
$UG$ is unipotent group of class at most $d-1$. The class is exactly $d-1$
because $UG$ contains the group of inner automorphism which is isomorphic to
$G/\Cen(G)$ and this is unipotent group of class exactly $d-1$.

The exact sequence splits because $\GL_n$ is reductive group and $UG$ is
an unipotent group, since $\R$ is a field.
\end{proof}

Similar result holds for the automorphism group of the discrete subgroup
$\Gamma(n,d)$ of $G(n,d)$.

\begin{lm}
The automorphism group $\Aut \Gamma(n,d)$ satisfies the following
exact sequence:
$$
1 \to \IAut \Gamma(n,d) \to \Aut \Gamma(n,d) \to \GL_n(\Z) \to 1,
$$
where the group $\IAut\Gamma(n,d)$ consist of all automorphisms,
which act trivially on $\Gamma/[\Gamma,\Gamma]$. This is a nilpotent group of
nilpotency class $d-1$. Note that the
exact sequence does not split into semidirect product unless $d=1$ and $\Gamma$
is abelian.
\end{lm}
\begin{proof}
The group $\Gamma(n,d)$ is a nilpotent group and a lattice in $G(n,d)$.
Therefore, its automorphism group can be embedded into $\Aut G(n,d)$.

Let us define $\IAut\Gamma(n,d) = \Aut \Gamma \cap UG$. Since $UG$
is a normal subgroup of $\Aut G$, $\IAut \Gamma$ is a normal
subgroup of $\Aut \Gamma$. The group $\IAut \Gamma$ is a lattice
in $UG$, since it is isomorphic to $\Hom(\Gamma/[\Gamma,\Gamma],[\Gamma,\Gamma])$
as a set.
The quotient $ \Aut \Gamma / \IAut \Gamma$ acts on
$\Gamma/[\Gamma,\Gamma] \simeq \Z^n$, therefore it is a subgroup
of $\GL_n(\Z)$.
It is the whole group because $\GL_n(\Z)$ is generated by the
elementary and diagonal matrices, and they can be easily lifted to
automorphisms of $\Gamma$.

Unlike in the case of $G(n,d)$ this exact sequence does not split in general.
In fact it splits only when $\Gamma$ is abelian and the group $\IAut \Gamma$ is trivial.
\end{proof}

Since we are interested in whether the group $\Aut \Gamma(n,d)$ has Kazhdan property $T$,
and property $T$ is invariant under passage to finite index subgroups,
in this chapter we will
work mainly with the subgroup of index $2$ in $\Aut \Gamma$ ---
$\Aut_1 \Gamma(n,d) = \pi^{-1}(\SL_n(\Z))$ consisting of all
automorphisms which act on the abelinization as matrices with determinant $1$.

\begin{lm}
The Zariski closure $\Aut_1 G(n,d)$ of the group $\Aut_1 \Gamma(n,d)$ in
$\Aut G$ is
$$
\Aut_1 G(n,d) = \SL_n(\R) \ltimes UG(n,d).
$$
It is a connected Lie group and the
group $\Aut_1 \Gamma$ is a lattice in $Aut^1 G(n,d)$.
\end{lm}
\begin{proof}
The group $\Aut_1 \Gamma(n,d)$ is a lattice in $\Aut_1 G(n,d)$ because
$\IAut \Gamma(n,d)$ is a lattice in $UG(n,d)$ and $\SL_n(\Z)$ is a
lattice in $\SL_n(\R)$.
\end{proof}

\begin{df}
Let us denote by $G^1(n,d)$ the Zariski closure of $\Aut \Gamma$ in
$\Aut G$. This is a Lie group which has two connected components and
is isomorphic to
$$
G^1(n,d) = \GL_n^{\pm1}(\R) \ltimes UG(n,d).
$$
\end{df}
\end{section}

\ifnewpage

\begin{section}{Kazhdan Property $T$}

Kazhdan property $T$ originated from the representation theory of
the Lie groups (see~\cite{Kazhdan}). For a group $G$, let $G^*$ denote the unitary dual of the group
$G$ consisting of all irreducible unitary representations of $G$ up to
isomorphism.
It is possible to put topology (called Felt topology)
on the set $G^*$, saying that two representations
$(\rho_0,\mc{H}_0)$ and $(\rho_1,\mc{H}_1)$ are close if it is
possible to find unitary vectors $v_0 \in \mc{H}_0$ and $v_1 \in \mc{H}_1$
such that the matrix coefficients $\la v_0, g v_0 \ra$ and $\la v_0, g v_0 \ra$
are close (on the compact subsets of $G$) as functions of $g$.
In this terminology it is said that $G$ has Kazhdan property $T$,
if and only the trivial representation is an isolated point in  $G^*$
with respect to this topology.

Here we will use an equivalent definition which is more appropriate for
our goals.

\begin{df}
A topological group $G$, generated by a compact set $Q$, is said
to have Kazhdan property $T$ if there exists a constant $\epsilon$
such that any unitary representation $(\rho,\mc{H})$ of the group
$G$, which contains a unit vector $v$ such that
$||\rho(g)v-v|| \leq \epsilon$
for any $g\in Q$, contains a $G$ invariant vector. The maximal
$\epsilon$ with this property is called the Kazhdan constant of
$G$ with respect to $Q$ and is denoted by $\mc{K}(G,Q)$. This is equivalent
to following expression
$$
\mc{K}(G,Q) = \inf_{(\rho, \mc{H})\in G^*}
\inf_{\begin{array}{c} \scriptscriptstyle v\in \mc{H} \\ \scriptscriptstyle ||v||=1\end{array}}
\sup_{g\in Q} ||\rho(g)v-v||,
$$
where the first infimum is over all unitary irreducible representations of
the group $G$.
\end{df}

\begin{reme}
The above definition says when the pair $(G,Q)$ has Kazhdan property $T$. It
can be shown that having property $T$ depends only on the group $G$ and its topology,
but the value of the Kazhdan constant depends also on the generating set $Q$.
\end{reme}
As an immediately corollary of the remark we have
\begin{co}
Compact (or finite) groups have property $T$, because for the generating set
$Q$ we can take $G$ itself.
In fact it can be shown that for any group $G$ we have
$\mc{K}(G,G) \geq \sqrt{2}$.
\end{co}

\begin{reme}
\label{finiteindex}
Let $G'$ be a finite index subgroup of a group $G$. Using induction and restriction of
representations from $G'$ to $G$ and vice versa it is easy to see that the group $G'$ has
Kazhdan property $T$ if and only if the group $G$ has property $T$.
\end{reme}

This result was generalized by Kazhdan in the case of Lie groups
and their lattices.

\begin{thq}[Kazhdan]
\label{lattice}
A lattice in a Lie group has property $T$ if and only if the
Lie group has property $T$.
\end{thq}

\begin{reme}
\label{auto:example}
The groups $\Z$ and $\SL_2(\Z)$ do not have Kazhdan property $T$ ---
the unitary dual of the group $\Z$ is isomorphic to the circle $\R/\Z$
and the Felt topology on it coincides with the usual topology on this set,
therefore there are no isolated points in the unitary dual of $\Z$. The group
$\SL_2(\Z)$ does not have property $T$ because it contains a free group as a
finite index subgroup.
\end{reme}
\begin{co}
If $G$ has property $T$ then any quotient also have $T$, in particular a group $G$
can not have property $T$ is it has infinite abelian quotient.
\end{co}

Proving that a group $G$ has a property $T$ directly using the above
definition requires detailed knowledge of the unitary dual $G^*$.
In~\cite{Kazhdan} Kazhdan avoided this problem for simple Lie groups.

\begin{thq}[Kazhdan]
Let $G$ be simple Lie group over $\R$, then $G$ has property $T$ if
the rank of $G$ is at least $2$.
\end{thq}
\begin{reme}
The situation with simple Lie groups of rank $1$ is more complicated - some
groups for example $F_{4,-22}$ have property $T$ and others like $\SL_2(\R)$
do not.
\end{reme}

An obvious obstruction for a Lie group $G$ to have property $T$ is to have $\R$
as a quotient (because the group $\R$ like the group $\Z$ does not have $T$).
It is interesting that for Lie groups with big semi simple part this
is the only obstruction. In order to show this we need to state
the following result proved by Wang in~\cite{Wang} as a consequence of the
so called Maunter Phenomena.

\begin{thq}[Wang]
\label{wang1}
Let $G$ be a Lie group over $\R$. Let $S$ denotes the reductive part and $N$
denotes the unipotent part of the group $G$. If the group $N$ coincides with
$[N,S]$ then there exist finitely many $S$ orbits $O_i$, $i=1,\dots, k$, in $N$,
such that any
element $n$ in $N$ can be written as
$n = o_1 o_2 \dots o_k$, where $o_i \in O_i$ or $o_i=1$.
\end{thq}
\begin{proof}
We will just sketch the proof. For details we refer the reader
to the original paper by Wang~\cite{Wang}.

The proof of this theorem is by induction of the nilpotency class of
$N$. The base case is when the unipotency class is $0$ and the group $N$ is trivial.
In this case there is nothing to prove.

By induction every element in $N/\Cen(N)$ can be written as a product of finitely many
$S$ orbits in $N/\Cen(N)$. We can lift these orbit to orbits in $N$, which
implies that every element in $N$ is a product of representatives of these orbits
up to an element in the center of $N$.

Therefore in order to prove that the same is true for $N$, we
only need to show that $\Cen(N)$ can be written as a product of finitely many
$S$ orbits.  The abelian group $\Cen(N)$ decomposes as a direct product of
finite number of subgroups $N_i$, which are simple $S$ modules.
If the action of $S$ on
$N_i$ is not trivial, then $N_i\setminus{0}$ consists of a single $S$ orbit.
Clearly every element in the product of nontrivial $N_i$-es can be
written as a product of elements in some finite number of $S$ orbits.

The condition $N = [N,S]$ implies that $N_k \subset [N,N]$
for any $k$ such that $N_k$ is isomorphic to the trivial $S$ module. Using that induction
hypothesis we can express every element in $N/\Cen(N)$ as a product of
finite number of $S$ orbits and therefore the same is true for the commutator of
any two elements in $N$.
This shows that we can find $S$ orbits $O_j$ such that every element in
$\Cen(N)$ can be written as a product of elements $o_j \in O_j$,
which completes the proof of the induction step.
\end{proof}



We can use this result to find a necessary and sufficient condition for a group with big nilpotent radical
to have Kazhdan Property $T$.

\begin{thq}[Wang]
\label{wang2}
Let $G$ be a Lie group over $G$. Let $S$ denotes the reductive part and $N$
denotes the unipotent part of the group $G$. Then the following two conditions are equivalent
\begin{enumerate}
\item
$G$ has property $T$;

\item
$S$ has property $T$ and $N / [S,N]$ is compact.
\end{enumerate}
\end{thq}
\begin{proof}
Clearly $1$ implies $2$, because $S$ is a quotient of $G$ and therefore has $T$,
and the quotient $N / [S,N]$ can be embedded into
$G/[G,G]$ therefore it is compact, because the only abelian groups
which have property $T$
are the compact ones.

We will only sketch why $2$ implies $1$.
Let construct a compact generating set $K$ of the group $G$ consisting of
a compact generating set of the reductive part $S$ together with
some compact subset of $N$ which projects onto $N / [S,N]$.
We need also to add to $K$ some generating set of $[S,N]$, let
$G' = S \ltimes [S,N]$ this group satisfies the conditions of
Theorem~\ref{wang1}. Therefore there exist a finite number of $S$ orbits
$O_i$ such that every element in $[S,N]$ is a product of representatives of
these orbits. Let us include in $K$ one representative from each orbit.
It is easy to see that the set $K$ constructed in this way is a
compact generating set of $G$.

Let us take a representation $(\pi,\mc{H})$ of the group $G$ and suppose that
$v$ is an $\epsilon$ invariant unit vector in $\mc{H}$ for epsilon small
enough, with respect to the generating set of $K$.

Using the fact that $S$ has property $T$ we can pass to a new vector $w$
which is $\epsilon'$ invariant with respect to the set $K$ but such that
the group $S$ acts trivially on $w$, where $\epsilon'$ depends only on $\epsilon$
and the Kazhdan constant of $S$.

The last condition implies that
$w$ is $\epsilon'$ invariant with respect to every element in the orbits $O_k$.
By Theorem~\ref{wang1} every element in the group $G$ can be expressed as
a product of a fixed number of elements which almost preserve the vector $w$,
therefore $w$ is moved a distance  less then $1$ by any element in $G$ if
$\epsilon$ is sufficiently small. This shows that $\mc{H}$ contains
invariant vectors, which proves that the group $G$ has Kazhdan property $T$.
\end{proof}
\end{section}

\ifnewpage

\begin{section}{Property $T$ of $\Aut \Gamma (n,d)$}

In~\cite{LubPak} A. Lubotzky and I. Pak studied the product replacement algorithm on
nilpotent groups. They proved that the working time of is algorithm is logarithmic
using the fact that a special subgroup of the automorphism
group $\Aut \Gamma(n,d)$ has property $T$ provided that $n\geq 3$.
Since this subgroup is of infinite index they asked the question
when the whole automorphism group has property $T$.
Using the results from previous sections we can answer this question.

\begin{thq}
\label{whenT}
The automorphism group $\Aut \Gamma(n,d)$ of the free nilpotent group on
$n$ generators and prepotency class $d$ has Kazhdan property $T$ if and only if
$n \geq 3$
and $d \leq h(n)$, where $h(n)$ is the function defined in
definition~\ref{gen:hes}. In
particular it has property $T$ if $d\leq n(n-1)$ and $n\geq 3$. If conjecture~\ref{gen:conj}
holds that it will have property $T$ for any $d$, provided that $n\geq 3$.
\end{thq}
\begin{proof}
The condition $n \geq 3$ is necessary because the group $\Aut \Gamma (2,d)$ has
$\GL_2(\Z)$ as a quotient and this group does not have property $T$, since it is
virtually free.

The main step in the proof of the theorem is the following lemma
\begin{lm}
The following statements are equivalent
\begin{enumerate}
\item
$\Aut \Gamma(n,d)$ has property $T$;
\item
$\Aut_1 \Gamma(n,d)$ has property $T$;
\item
$\Aut_1 G(n,d)$ has property $T$;
\item
$\Aut_1 G(n,d)$ does not have $\R$ as a quotient and $n \geq 3$;
\item
the abelinization of the Lie algebra $\sll_n + UL(n,d)$  is trivial and $n \geq 3$;
\item
$n \geq 3$ and
$(UL(n,d)/[UL(n,d),UL(n,d)])^{\sll} =0$.
\end{enumerate}
\end{lm}

\begin{proof}
Many of the equivalences among the above statements are easy to prove. Later we will use
that 1 is equivalent to 6 in order to prove Theorem~\ref{whenT}.

$1$ is equivalent to $2$, because $\Aut_1 \Gamma(n,d)$ is a subgroup
of index $2$ in $\Aut \Gamma(n,d)$ and we can use remark~\ref{finiteindex};

$2$ is equivalent to $3$, by Kazhdan Theorem~\ref{lattice}
because $\Aut_1 \Gamma(n,d)$ is a lattice in
$\Aut_1 G(n,d)$;

$3$ and $4$ are equivalent by the Maunter phenomenon (Theorem~\ref{wang2}),
because the group $\Aut_1 G(n,d)$ is connected, simply
connected and can not have any non trivial compact quotients;

$4$ and $5$ are the same by Lie algebra arguments;

$5$ and $6$ because the abelinization of the Lie algebra
$\sll_n + UL(n,d)$ is equal to $UL(n,d)/[UL(n,d),UL(n,d)]^{\sll}$.
\end{proof}

By definition $h(n)$ is the smallest degree of homogeneous element
in the factor space
$UL[V]/[UL[V],UL[V]]^{\sll}$. By construction the algebra $UL(n,d)$ is
a quotient of $UL[V]$ containing all homogenous components of degree less then $d$.
Therefore the space $UL(n,d)/[UL(n,d),UL(n,d)]^{\sll}$ is trivial if and only if
$d\leq h(n)$, which completes the proof of the theorem.

\end{proof}

\end{section}

\ifnewpage

\begin{section}{The center of $\Aut \Gamma(n,d)$}

In this section we give a proof of a result by Formanek\cite{Formanek}, which describes
the center of the automorphism group of $\Gamma(n,d)$.
We will use this result latter in Chapter 4, where we will
describe the automorphism tower of the group $\Gamma(n,d)$.

\begin{thq}[Formanek]
The center of the group $\Aut \Gamma(n,d)$ is not trivial
if and only if
$2n | d-1$.
\end{thq}
\begin{proof}
Every element in the center of the group $\Aut \Gamma$ should lie
in the center of the group $\IAut \Gamma$. We have seen that
$$
\Cen(\IAut \Gamma) = \Hom (\Gamma/[\Gamma,\Gamma], \Cen \Gamma).
$$
An element $\phi \in \Cen(\IAut \Gamma)$ lie in the center of the group $\Aut \Gamma$
iff it is preserved by the action of $\GL_n(\Z)$ on this space.
Since the action of $\GL_n(\Z)$ in these spaces is polynomial of degree
$1$ and $d$ respectively, a necessary condition for existence of
invariant elements is that $2n | d-1$.

Therefore we only need to show that if $2n | d-1$ there exists  $\GL_n(\Z)$
invariant elements in the $\Hom (\Gamma/[\Gamma,\Gamma], \Cen \Gamma)$.
If $d > 4n$ or $n \geq 3$
this is true because that are $\GL_n(\Z)$ invariant elements in
the homogeneous component of the free Lie algebra of degree $2kn$ for $k\geq 2$
or $n > 3$, and the homogeneous component of the free Lie algebra can
be embedded in the above space of homomorphisms.
In the cases $n=2$, $d=5$ and $n=3$, $d=7$
we need to verify that in the homogeneous component of the free Lie algebra
of degree $d$ there is a $\sll$ submodule isomorphic to the simple
module  corresponding to the partition $[2k+1,2k^{n-1}]$.
\end{proof}

\end{section}

\end{chapter}
\newpage

\begin{chapter}{Derivation Tower of Free Nilpotent Lie algebras}
\thispagestyle{empty}
	\nocite{DyerForm1}

In this chapter we study the derivation tower of free nilpotent
Lie algebras.  It is known that in the case of Lie groups there is
a close connection between the automorphism group of the Lie group
and the derivation algebra of its Lie algebra. Using this
connection is natural first one to study the derivation tower
of the Lie algebra of a group and after that to try to lift the
results to the automorphism tower of the Lie group.

It turns out (see Theorem~\ref{fulltower}) that the derivation
tower of free nilpotent Lie algebras is very short --- it has
height at most 3 and similar result holds for the automorphism
tower of the free unipotent group.

Unfortunately the situation with the free nilpotent group is
different. One explanation for this difference is that in the ring
$\Z$ there are just $2$ invertible elements and therefore the
automorphism group of a lattice in a Lie group is much smaller
then the automorphism group of the Lie group. A very illustrative
example is $\Z^n \subset \R^n$.  Their automorphism groups are
$\GL(\Z)$ and $\GL(\R)$ respectively. If one wants to use some Lie
group to study the group $\GL(\Z)$ it is better to use $\SL^{\pm
1}(\R)$ not the whole group $\GL(\R)$.

Having in mind this discrepancy, we define the algebra of
restricted derivations (see definition~\ref{der0}) and construct
the tower of restricted derivations. Studying this tower, we found
that the behavior of the tower is different depending whether the
first algebra in this tower $\D_0 L[V,d]$ has trivial center or
not. The main results in this chapter are
Theorems~\ref{dertower:nocenter} and~\ref{dertower:center}, which
state that the tower of restricted derivations of the free
nilpotent Lie algebra stabilizes after finitely many steps.

In the next chapter 4 we will lift these result to the
automorphism tower of the free nilpotent group.

\begin{section}{Derivation tower of Lie algebras}

If we have a Lie algebra $L$ we can construct its derivation
algebra $\D L$, which comes with a natural map $\ad: L \to \D L$.

\begin{df}
A Lie algebra $L$ is called complete if it has trivial center and
no outer derivations. This is equivalent to the map $\ad$ being an
isomorphism.
\end{df}

Iterating the procedure of constructing the derivation algebra, we
can define the derivation tower of the Lie algebra $L$.

\begin{df}
Let $L$ be a (finite dimensional) Lie algebra, by a derivation
tower of the algebra $L$ we mean, the sequence
$$
\D^0 L \to \D^1 L \to D^2 L \to \dots \to \D^n L \to \dots,
$$
where $\D^0 L =L$ and for any $i$ we have $\D^{i+1}L = \D (\D^i L)$.
The maps $\ad_i : \D^{i}L \to \D^{i+1}L$ comes form the inner derivations
of the algebra $\D^{i}L$. In general this maps are neither
injective nor surjective.
\end{df}

\begin{reme}
If the algebra $\D^i L$, for some $i$, has a trivial center then
the map $\ad_i$ is an injection. By Remark~\ref{der:centralizer})
we have that the algebra $\D^{i+1} L$ also has a trivial center
and therefore by induction all maps $\ad_k$ are injections for
every  $k\geq i$.
\end{reme}

\begin{df}
We say that the derivation tower of the algebra $L$ stabilizes at
level $k$ (or after $k$ steps) if the map $\ad_k$ is an
isomorphism. We also say that the tower weakly stabilizes at level
$k$, iff $\D^k L \simeq \D^{k+1} L$.
\end{df}

\begin{ex}
\label{tow:example} Derivation tower of the algebra $\sll_n$. The
algebra $\sll$ has trivial center and has no outer derivations.
Therefore the derivation tower of $\sll_n$ stabilizes at the level
$0$ and the tower is
$$
\sll_n \to \sll_n \to \dots \to \sll_n \to \dots.
$$

\end{ex}

\begin{ex}
\label{tow:example1}
Derivation tower of $\R$. The derivation
algebra of $\R$ is also $\R$, but the map $\ad$ is trivial because
$\R$ is abelian Lie algebra. This shows that the derivation tower
of $\R$ is
$$
\R \to \R \to \dots \to \R \to \dots,
$$
and all maps $\ad_k$ are trivial. Therefore the tower does not stabilize,
 but it stabilizes weakly at level $0$.
\end{ex}

\end{section}

\ifnewpage

\begin{section}{Derivation Algebras}

First let us prove several lemmas about derivation algebras, their
centers and centralizers of certain ideals. We will use these
lemmas in the next section to describe the derivation algebra of
$\D L[V,d]$.

\begin{lm}
\label{der:centralize}
Let $L$ be a Lie algebra. Then the set of all inner derivations
$\ad L$, is an ideal in the algebra $\D\,L$ and its centralizer is
isomorphic to
$$
\Cen_{\D\,L}(\ad(L)) = \Hom(L/[L,L], \Cen(L)).
$$
\end{lm}
\begin{proof}
Let $d\in \D\,L$ be a derivation of $L$, which commutes with all
inner derivations. Then
$$
0=[d,\ad x](y) = d([x,y]) - [x,d(y)],
$$
but we have that
$$
d([x,y]) = [d(x),y] + [x,d(y)].
$$
From these equalities it follows that $[d(x),y]=0$, for all $x,y
\in L$. Therefore $d(x)\in \Cen(L)$ for all $x \in L$. From the
Leibnitz rule, it follows that $d([x,y])=0$, i.e., $d_{|[L,L]} =
0$. This shows that we can view $d$ as a map from $L/[L,L]$ to
$\Cen(L)$. This is an isomorphism, because every element in
$\Hom(L/[L,L], \Cen(L))$ can be extended to a derivation of $L$,
by letting it act trivially on the commutator subalgebra, and the
direct computation verifies that this defines a derivation which
commutes with all inner derivations.
\end{proof}

\begin{reme}
\label{der:centralizer} If the algebra $L$ has a trivial center
then the centralizer of the ideal of inner derivations is trivial,
in particular we have that the algebra $\D\,L$ also has a trivial
center.
\end{reme}

\begin{lm}
\label{der:acttrivally} Let $L$ be a Lie algebra. We can consider
the inner derivations as a map $\ad:L \to \D\,L$. Let $D$ be a
derivation of $\D\,L$ (i.e.\ an element in $\D(\D\,L)$), which
acts trivially on $\ad(x)$ for all $x\in L$. Then $(D(d))(x)\in
\Cen(L)$ for all $x\in g$ and $d\in \D\,L$, i.e., $\im D \subset
\Cen_{\D\,L}(\ad(L))$.
\end{lm}
\begin{proof}
Let us note that $[d,\ad(x)] = \ad(d(x))$ for all $d\in \D\,L$.
Applying the derivation $D$ to both sides yields
$$
0=D(\ad(d(x))) = D [d,\ad(x)] = [D(d),\ad(x)]  + [d, D(\ad(x))]
= \ad (D(d)(x)),
$$
therefore $D(d)(x)\in \Cen{g}$ for all $x$. From
lemma~\ref{der:centralize} it follows that $\im D \subset
\Cen_{\D\,L}(\ad(L))$.
\end{proof}
\begin{co}
Let $L$ be a Lie algebra with trivial center. The adjoint action
gives a map $\ad:L \to \D\,L$. Let $D$ be a derivation of $\D\,L$
(i.e.\ an element in $\D(\D\,L)$), which acts trivially on
$\ad(x)$ for all $x\in L$. Then $D=0$. This shows that any
derivation of $\D L$ is determined by its restriction to $\ad(L)$.
\end{co}

Lemma~\ref{der:acttrivally} can be generalized to derivation
acting trivially on ideals in Lie algebras.

\begin{lm}
\label{der:Cenofideal} Let $D\in \D L$ be a derivation which acts
trivially on the ideal $I$, i.e.,\ $D(i)=0$ for all $i\in I$. Then
$\im D \in \Cen_{L}(I)$.
\end{lm}
\begin{proof}
Let $g\in L$ and $i\in I$. Since $I$ is an ideal we have $[g,i]\in
I$, if we apply the derivation $D$ to both sides we get
$[D(g),i]+[g,D(i)]=0$, which gives that $[D(g),i]=0$, i.e., $D(g)
\in \Cen_{L}(I)$ for any $g\in L$.
\end{proof}

\end{section}

\ifnewpage

\begin{section}{Chevalley Theorem}

In~\cite{Cheva} Chevalley stated that if $\mr{g}$ is a finite
dimensional Lie algebra, over a field of characteristic $0$, with
trivial center, then its derivation tower stabilizes after
finitely many steps. This result was generalized by Schenkman to the
case of positive characteristic see~\cite{Schenkman}. However, this
result can not be generalized to arbitrary finite dimensional Lie
algebra, because there are examples (like
Example~\ref{tow:example1}) of Lie algebras with non stabilizing
derivation towers.

Chevalley's  note contains only a sketch of the proof, in this
section we give a detailed proof of this result, because in the
next sections we will make several constructions similar to the
ones in this proof.

The following technical lemma is a key point in the proof.

\begin{lm}
\label{ch:tech}
Let $\mr{g}$ be a Lie algebra with trivial center, which has
the following decomposition
$$
\mr{g}= \mr{s} + \mr{a} + \mr{n},
$$
where $\mr{s}$ is a semi simple Lie algebra, $\mr{a}$ is
abelian and $\mr{n}$ is nilpotent, such that the following
properties are satisfied:
\begin{list}{}{}
\item $\mr{a}$ commutes with $\mr{s}$;
\item $\mr{a}$ acts diagonally on $\mr{n}$;
\item $\mr{n}$ is an ideal in $\mr{g}$;
\item $\mr{n}$ is generated by 
$[\mr{s},\mr{n}]$ and $[\mr{a},\mr{n}]$  as a Lie subalgebra.
\end{list}
Let  $\D\, \mr{g}$ denotes the derivation algebra of $\mr{g}$. Since
$\mr{g}$ has no center we can consider $\mr{g}$ as an ideal of
$\D\, \mr{g}$. Let $d$ be a derivation from $\mr{g}$ to
$\D\, \mr{g}$, i.e.\ a linear map which satisfies the equality
$$
d[x,y]= [d(x),y] + [x, d(y)]
$$
for any $x,y \in \mr{g}$. Then $\im d \subset \mr{g}$, i.e., $d\in
\D\,\mr{g}$.
\end{lm}
\begin{proof}
We have that $d[x,y]\in \mr{g}$ for any $x,y \in \mr{g}$, because
$\mr{g}$ is an ideal in $\D\, \mr{g}$. Therefore,
$d[\mr{g},\mr{g}]\in \mr{g}$.
From the properties of the decomposition of $\mr{g}$
it follows that $[\mr{g},\mr{g}] = \mr{s}+ \mr{n}$. Therefore
we only need to show that $d(a) \subset \mr{g}$.

The Lie algebra $\mr{a}$ acts diagonally on $\mr{g}$ and
$\D\,\mr{g}$. This action defines a $\mr{a}^*$ gradation on both
$\mr{g}$ and $\D\, \mr{g}$. Without loss of generality, we
may assume that the derivation $d$ is homogeneous with respect to
the gradation of $\D\, \mr{g}$.

Suppose that the degree $\lambda \in \mr{a}^*$ of $d$ is different
from zero. Then $d(\mr{a}) \subset V$ for some nontrivial $\mr{a}$
module $V$, such that $[a,v] = \lambda(a) v$ for any $a\in \mr{a}$
and $v\in V$. But if $\lambda \not = 0$ then $[\mr{a},V] = V$,
which implies that $V\subset\mr{g}$, since $\mr{a}\subset \mr{g}$
and $\mr{g}$ is an ideal in $\D\, \mr{g}$. Thus we have shown that
$d(\mr{a})\subset \mr{g}$, which implies that $d(\mr{g})\subset
\mr{g}$.

Now suppose that the degree of $d$ is zero. Let $v$ be an
eigenvector of $\mr{a}$ in $\D\,\mr{g}$ corresponding to the
eigenvalue $\lambda \in \mr{a}^*$. By the homogeneity of $d$ it
follows that $d(v)$ is also an eigenvector for $\mr{a}$,
corresponding to the same eigenvalue. Thus,
$$
d[a,v] = d(\lambda(a) v) = \lambda(a) d(v)
$$
and
$$
d[a,v] = [d(a),v] + [a,d(v)] = [d(a),v]+  \lambda(a) d(v).
$$
Therefore, $[d(a),v] = 0$, for any $a\in\mr{a}$ and any
eigenvector $v$. This implies that $[d(a), x]=0$ for any $x\in
\D\,\mr{g}$ because $\mr{a}$ acts diagonally on $\D\,\mr{g}$. From
Remark~\ref{der:centralizer} it follows that the center of the algebra $\D\,\mr{g}$ is
trivial, which implies that $d(a)=0$ for all $a\in\mr{a}$. Thus we
have $d(\mr{a})=0 \subset \mr{g}$, i.e.\ $d(\mr{g})\subset
\mr{g}$, which completes  the proof of the lemma.
\end{proof}

\begin{co}
\label{ch:complete}
Let $\mr{g}$ be a Lie algebra which satisfies the conditions
of the lemma, then its derivation algebra
is complete, i.e., $\D (\D \mr{g}) = \D \mr{g}$.
\end{co}

\begin{thq}
\label{ch:main}
Let $\mr{g_0}$ be a finite dimensional Lie algebra over a
field of characteristic $0$, with a trivial center. Then
the derivation tower $\D^n(\mr{g_0})$ of $\mr{g_0}$ stabilizes
after finitely many steps.
\end{thq}
\begin{proof}
The algebra $\mr{g_0}$ has the following decomposition
$$
\mr{g_0}= \mr{s} + \mr{a} + \mr{n_0},
$$
where $\mr{s}$ is a semi simple Lie algebra, $\mr{a}$ is
abelian and $\mr{n_0}$ is nilpotent, such that the following
properties are satisfied:
\begin{list}{}{}
\item $\mr{a}$ commutes with $\mr{s}$;
\item $\mr{a}$ acts diagonally on $\mr{n_0}$;
\item $\mr{n_0}$ is the maximal nilpotent ideal of $\mr{g_0}$.
\end{list}
Let $\mr{n}$ be the subalgebra of $\mr{n_0}$ generated by the
elements $[\mr{s},\mr{n_0}]$ and $[\mr{a},\mr{n_0}]$. Let $\mr{g}$
denote the subalgebra $\mr{g}= \mr{s} + \mr{a} + \mr{n}$.

\begin{lm}
The subalgebra $\mr{g}$ is an ideal in $\mr{g_0}$.
\end{lm}
\begin{proof}
It suffices to show that $[[a,n_1],n_2]\in \mr{n}$ for any $a\in
\mr{a}$ and $n_1,n_2 \in \mr{n_0}$, and that $[[s,n_1],n_2]\in
\mr{n}$ for any $s\in \mr{s}$ and $n_1,n_2 \in \mr{n_0}$. Let us
prove the first one. Without loss of generality, we may assume
that $n_2$ is an eigenvector for $\mr{a}$. If the corresponding
eigenvalue is zero, then $[\mr{a},n_2]=0$, therefore
$[[a,n_1],n_2] = [a,[n_1,n_2]] \in [\mr{a},\mr{n_0}] \subset
\mr{n}$. If the corresponding eigenvalue is not zero then
$n_2=[a_1,n_2]$ for some $a_2\in \mr{a}$, which gives that
$[[a,n_1],n_2] = [[a,n_1],[a_1,n_2]] \in \mr{n}$. The proof that
$[[s,n_1],n_2]\in \mr{n}$ for any $s\in \mr{s}$ and  $n_1,n_2 \in
\mr{n_0}$ is similar.

\end{proof}

\begin{lm}
The centralizer of $\mr{g}$ in $\mr{g_0}$ is trivial.
\end{lm}
\begin{proof}
Let $J = \{ x\in \mr{g_0} \mid [x,\mr{g}]=0 \}$ be the
centralizer of $\mr{g}$. It is easy to see that
$J\subset \mr{n_0}$ and $J$ is an ideal in $\mr{g_0}$.
If $J$ is non trivial then
$J \cap \Cen{\mr{n_0}} \not = \{0\}$,
because in the nilpotent Lie algebra every nontrivial
ideal has nontrivial intersection with the center.
It is easy to see that
$J \cap \Cen{\mr{n_0}} \subset \Cen{\mr{g_0}} = \{0\}$,
which is a contradiction.
Therefore, $J= \{0\}$, which proves the lemma.
\end{proof}

The above two lemmas allows us to identify the algebra $\mr{g_0}$
with some sub algebra of $\D\,\mr{g}$. Let $D\in \D(\mr{g_0})$.
Then the restriction of $D$ to $\mr{g}$, is a derivation from
$\mr{g}$ to $\mr{g_0}\subset \D\,\mr{g}$. By Lemma \ref{ch:tech}
we have that $D_{|\mr{g}}=D' \in \D\,\mr{g}$. We can consider the
$D'$ as a derivation from $\mr{g_0}$ to $\D\,\mr{g}$. By
Lemma~\ref{der:centralize} we have that $D-D'=0$ (in
$\D\,\mr{g}$), which allows us to identify $\D(\mr{g_0})$ as a sub
algebra of $\D\,\mr{g}$. Now by simple induction we have
$$
\mr{g} \subset \mr{g_0}
\subset \D^1(\mr{g_0})
\subset \cdots \subset \D^n (\mr{g_0}) \subset \cdots
\subset \D\,\mr{g}.
$$
Since the algebra $\D\,\mr{g}$ is finite dimensional, the above
sequence has to stabilize at some point. Therefore, the
derivation tower of $\mr{g_0}$ stabilizes.
\end{proof}
\begin{reme}
This proof also gives us a bound on the height of the derivation
tower --- the height is less than $\dim \D\,\mr{g} - \dim
\mr{g_0}$.
\end{reme}

\end{section}

\ifnewpage

\begin{section}{Derivation tower of free nilpotent Lie algebras}

Let us describe the derivation tower of the free nilpotent Lie
algebra $L[V,d]$. We can not use  directly the Chevalley's Theorem
because the algebra $L[V,d]$ has nontrivial center. However we can
apply Chevalley's result to its derivation algebra which does
not have center. In chapter 1, we showed that its derivation
algebra has a decomposition
$$
\D L[V,d ] = \gll_n + UL[V,d] = \sll_n + \R.1 + UL[V,d]
$$
This decomposition is of the same form as the one needed for
Lemma~\ref{ch:tech}. Using the fact that $1$ acts as a degree
derivation on $UL$, which is positively graded, we have
$[1,UL]=UL$. Therefore the conditions in the technical lemma~\ref{ch:tech} are
satisfied, which gives us that the derivation algebra of $\D
L[V,d]$ is complete. This shows that the derivation tower of the
algebra $L[V,d]$ stabilizes at most at the third level.

Let us describe the derivation algebra of $\D L[V,d]$ in detail.
Now let us assume that $d \geq 2$. Without loss of generality we
may assume that any outer derivation $d$ acts trivially on $1$,
$\sll_n$ and $\ad(V)$ (see next section). Therefore $d$ preserves
the grading and the $\sll_n$ module structure on the $UL$. Using
the fact that $d$ acts trivially on the space $\ad(V)$, it follows
that $d$ kills every element of the form $\ad(f)$ for some $f\in
L[V,d]$. The last condition implies that $\im d \subset \Cen UL$.
This argument gives that if $d\geq 3$ there is an isomorphism
between the space of outer derivations of the algebra $\D L[V,d]$
and the space
$$
\Hom_\sll \left(U_{\gen} \cap UL^{(d-1)}, UL^{(d-1)} \right).
$$
In particular the set of outer derivations is not trivial because
$U_{\gen} \cap UL^{(d-1)}$ is not empty (it contains at least the
$\sll$ module generated by $\AD_{x_1^n}$).

In the case $d=2$ situation is a bit more complicated because
$\ad(V) \subset U_{\gen} \cap UL^{(d-1)}$. If $n \geq 3$ then
$U_{\gen} \cap UL^{(1)} = UL^{(1)}$ is a sum of two simple $\sll$
modules, therefore the space of outer derivations is one
dimensional. In the case $n=2$ we have that $UL^{(1)} = \ad(V)$,
therefore there are no outer derivations.

Putting these arguments together we can describe the derivation tower of the
algebra $L[V,d]$.

\begin{thq}
\label{fulltower}
a) If the nilpotency class $d$ is at least $3$ or
the number of generators is at least $3$
then the derivation tower of the
free nilpotent Lie algebra $L[V,d]$ is
$$
L[V,d] \to \gll + UL[V,d] \to \gll + UL[V,d] +
\Hom_\sll \left(U_{\gen}^+ \cap UL^{(d-1)}, UL^{(d-1)} \right),
$$
and all other terms are equal to the $\D^2 L[V,d]$, i.e., the tower
stabilizes at the second level.

b) If the nilpotency class is $d=2$ and $n = 2$ then the tower
stabilizes at the first level.
\end{thq}

\end{section}

\ifnewpage

\begin{section}{Algebra of restricted derivations}

Our goal is to use the description of the derivation tower of the
algebra $L[V,d]$ to obtain information about the automorphism
tower of the free nilpotent group. In chapter 2 we saw the `Lie
algebra' corresponding to the free nilpotent group is $L[V,d]$ and
the `Lie algebra' corresponding to its automorphism group is
$$
\sll_n + UL[V,d],
$$
which is not the whole derivation algebra of $L[V,d]$, but a
subalgebra of codimension $1$. One reason for that is the fact the
in the ring $\Z$ there are just two invertible elements, therefore
the action of any automorphism on any characteristic quotient (as
the abelinization of $\Gamma(n,d)$) should have `determinant' plus
or minus $1$. This implies that `derivations' which corresponds
to such automorphism should has trace zero.

This argument suggest the instead of taking the whole derivation
algebra we should take the subalgebra of all `traceless'
derivations.

\begin{df}
\label{der0}
Let $L$ be a Lie algebra and let $\D_0 L$ denote the
subalgebra of $\D L$ consisting of 'totally traceless'
derivations, i.e,
$$
\D_0 L = \{ D \in \D L \mid \tr D_I =0 \mbox{ and }
\tr D_{L/I} =0
\mbox{ for any characteristic ideal }I\in L \}.
$$
We will call $\D_0 L$ the algebra of restricted derivations.
\end{df}
\begin{reme}
Using the fact the commutator of any two linear maps has zero trace
we can see that $[\D L, \D L] \subset \D_0 L$. This implies that
$\D_0 L$ is a subalgebra $\D L$. We will call this subalgebra the
algebra of restricted derivations of $L$.
\end{reme}

\begin{reme}
In general we do not have that $\ad(\mr{g}) \subset \D_0 \mr{g}$.
For example take $\mr{g}$ to be the two dimensional Lie algebra
with basis $a,b$ and commutator relation $[a,b]=b$. In this case
we have that $\tr \ad(a) = 1$, therefore $\ad (a) \not \in \D_0
\mr{g}$.

In general any Lie algebra $\mr{g}$ can be written uniquely as
$$
\mr{g} = \mr{s} + \mr{a} + \mr{n},
$$
where $\mr{s}$ is semi simple, $\mr{a}$ is abelian and the action
of $\ad(\mr{a})$ on $\mr{g}$ is diagonalizable, and $\mr{n}$ is
nilpotent. We have that $\ad(\mr{g}) \subset \D_0{\mr{g}}$, if and
only if $\mr{a}$ is trivial or $\mr{a}$ acts trivially on
$\mr{n}$.

However, if there exists a Lie group $G$ over $\R$ and Zariski
dense lattice $\Gamma$ in $G$ such that the Lie algebra of $G$ is
$\mr{g}$ then the above condition is satisfied and we have
$\ad(\mr{g}) \subset \D_0 \mr{g}$.
\end{reme}

The next lemma allows us to define the tower of restricted
derivations

\begin{lm}
Let $\mr{g}$ be a finite dimensional Lie algebra, then
$$
\ad (\D_0 \mr{g}) \subset \D_0(\D_0\mr{g}).
$$
\end{lm}
\begin{proof}
Let us decompose $\D \mr{g}$ as
$$
\D \mr{g} = \mr{s} + \mr{a} + \mr{n},
$$
and consider its action on $\mr{g}$. The action of $\mr{a}$ is
diagonalizable and defines an $\mr{a}^*$ grading on $\mr{g}$, which
decomposes into homogeneous components as $\mr{g} = \sum V^{\mu}$.

Let $\mr{n_0}$ be the subalgebra of $\mr{n}$ generated by
$[\mr{n}, \mr{a}]$. The direct computation verifies that
$\mr{n_0}$ is an ideal in $\D \mr{g}$ and its homogeneous
component of degree $0$ acts nilpotently on any $V^\mu$.

Therefore the space $\tilde {\mr{g}} = \mr{g} / \mr{n_0}(\mr{g})$,
has nonzero homogeneous components $\tilde V^\mu$ (under the
action of $\mr{a}$)  for any $\mu$ such that $V^\mu$ is not zero.
This allows us to construct an ideal $I$ in $\mr{g}$ which is
invariant under all derivation such that $\tr_{\mr{g}/I}{\ad(a)}
\not = 0$ for any $a\in \mr{a}$, which acts nontrivially on
$\mr{g}$. Therefore $\mr{a} \not \subset D_0 \mr{g}$ and we have
that
$$
\D_0 \mr{g} = \mr{s} + \mr{n},
$$
and by the previous remark we have that $\ad(D_0 \mr{g}) \subset \D_0( D_0 \mr{g})$,
which finishes the proof.
\end{proof}

\begin{df}
Restricted derivation tower. Let $L$ be an Lie algebra
such that $\ad(L) \subset \D_0 L$. Then we
can define the restricted derivation tower of $L$ as follows
$$
L = \D^0_0 L\to \D^1_0 L \to \D^2_0 L \to \dots \to \D^n_0 L \to \dots
$$
where $\D^{n+1}_0 L = \D_0( \D^n_0 L)$ and the maps come from the
inner derivations. We say the  tower stabilizes at the level $k$
if $\ad_k$ is an isomorphism from $\D_0^k L \to D_0^{k+1} L$.
\end{df}

\begin{reme}
Let $\mr{g}$ be an algebra with trivial center. Then $\mr{g}
\subset \D_0 \mr{g}$ if $\mr{a}=0$ or if $\mr{a}$ acts trivially
on $\mr{n}$ (here we use the decomposition from
Theorem~\ref{ch:main}). In this case it is easy to check that the
proof about the stabilization of the derivation tower carries on to
the tower of restricted derivations without any modifications.
\end{reme}

\begin{ex}
Restricted derivation tower of $\R^n$.
The derivation algebra of $\R^n$ is $\gll_n$, which can be written as
$\sll_n + \R.\id$. The subalgebra of restricted derivations consist only of
$\sll_n$, because $\id$ acts on $\R^n$ as an operator with trave $n \not =0$.
This gives that $\D_0 \R^n = \sll_n$. All other algebras in the tower of restricted
derivations of $\R^n$ are equal to $\sll_n$, because the algebra $\sll_n$ is complete
and coincide with its derived algebra. Therefore the tower is
$$
\R^n \to \sll_n \to \sll_n \to \dots \to \sll_n \to \cdots
$$
and stabilizes at the first level.
\end{ex}

Finally let us define the operator $\Hom^0$ which is similar of
the $\D^0$
\begin{df}
Let $V$, $W$ and $U$ are linear spaces such that we have a
projection $\pi: W \to U$ and an inclusion $i: U \to V$. Then
there is a natural projection form $\Hom(V,W)$ to $\End(U)$. Let
us define the space $\Hom^0(V,W)$, consisting of all elements in
$\Hom(V,W)$ whose projection in $\End(U)$ has zero trace for all invariant
$U$-es which are naturally a subspaces of $V$ and factor factor spaces of $W$.
\end{df}

\end{section}

\ifnewpage

\begin{section}{Distinguished ideals in $UL$}

Before describing the tower of restricted derivations of the algebra $L[V,d]$
we need to define several invariant ideals in the algebras $UL[V]$ and $UL[V,d]$.
The Lie algebra $\sll$ is simple and it acts naturally on $UL$,
therefore we can split $UL$ as a direct sum
$$
UL[V] = UL^{\triv}[V] + \sum_{\lambda\not= \triv} UL^\lambda[V],
$$
where the sum is over all nontrivial partitions containing less
than $n$ parts. Here, $UL^{\triv}$ is the maximal trivial $\sll$
submodule of $UL$ and $UL^\lambda$ is the maximal submodule, which
can be written as a sum of simple $\sll$ submodules corresponding
to the partition $\lambda$. Let us denote $UL^+ = \sum
UL^\lambda$, this is the maximal submodule, which can be written
as a sum of nontrivial $\sll$ submodules.

\begin{df}
\label{dist:ideals}
Let us define the following submodules of $UL[V]$:
\begin{list}{}{}
\item
$UL_\infty$ as the ideal in $UL$ generated by
$[UL^+,UL^+]$;

\item
$UL^+_{(0)} = UL^+$ and $UL^+_{(k+1)} = [UL^+_{(k)},UL^{\triv}]$;

\item
$UL^\lambda_{(0)} = UL^\lambda$ and $UL^\lambda_{(k+1)}
= [UL^\lambda_{(k)},UL^{\triv}]$;

\item
$UL_k = UL_\infty + UL^+_{(k)}$;

\item
$UL^\lambda_k = UL_{k+1} + UL^\lambda_{(k)}$.

\end{list}
\end{df}

The submodules $UL^{\lambda}$ are preserved under all automorphisms
of the algebra $\sll+ UL$. Therefore all the ideals defined above
are also invariant.
\begin{lm}
\label{seq}
The submodules $UL_k$ and $UL^\lambda_k$ form a descending
sequence of characteristic ideals in the algebra $\D_0 L[V]$, i.e.,
$$
UL_0 \supset UL_1 \supset \dots \supset UL_k
\supset \dots \supset UL_\infty
$$
$$
UL^\lambda_0 \supset UL^\lambda_1 \supset \dots \supset UL^\lambda_k
\supset \dots \supset UL^\lambda_\infty.
$$
\end{lm}
\begin{proof}
$UL_k$ is an ideal in $UL$ because, we have the inclusions
$[UL_{(k)}, UL^+] \subset [UL^+,UL^+]$ and $[UL_{(k)}, UL^{\triv}]
\subset UL_{(k+1)} \subset UL_{(k)} \subset UL_k$. These inclusions
prove that for every $k$,  $UL_k$ is an ideal in $UL$. The proof
that $UL_k^\lambda$ are ideals is similar.
\end{proof}
\begin{df}
Let us denote by $U_{\gen,k} = UL_k / UL_{k+1}$ and
$U_{\gen,k}^\lambda = UL_k^\lambda / UL_{k+1}$. We use the
subscript `$\gen$' because a set $S$ generates the algebra $UL$ if
and only if its image in $U_{\gen,0}$ generates that space. We
also want to define
$$
U_{\gen}^{\triv} = UL^{\triv} /\big([UL,UL] \cap UL^{\triv}\big).
$$

These modules will play an important role in the description of
the derivation tower of $L[V,d]$. Also it is convenient to take
the image of $\ad(V)$ away from $U_{\gen,0}$ and write $U_{\gen,0}
= \ad(V) + U_{\gen,0}^+$.
\end{df}


Similarly, we can define
$UL_k[V,d]$, $UL_k^\lambda[V,d]$ and so on,
as the submodules (ideals) of $UL[V,d]$ constructed form the algebra $UL[V,d]$
in the same way as the other modules were constructed from the algebra $UL$.
The submodules $UL_k[V,d]$ and $UL^\lambda_k[V,d]$ again form a descending
sequence of characteristic ideals in the algebra $\D_0 L[V,d]$.

\begin{thq}
\label{ULstabilization}
The sequences  $UL_k[V,d]$ and
$UL^\lambda_k[V,d]$ of characteristic ideals in $UL[V,d]$
stabilize after at most $1+ d/h(n)$ steps. In particular if
Conjecture~\ref{gen:conj} is true (or $d\leq h(n)$), they
stabilize at the second term, i.e., $UL_1[V,d] = UL_2[V,d]$.
\end{thq}
\begin{proof}
In the factor space $UL^{\triv} /\big([UL,UL] \cap
UL^{\triv}\big)$ there are no homogeneous components of degree
less then $h(n)$. By induction, it follows that in
$$
\big(UL^+_{(k)} + [UL^+,UL^+]\big)/\big( UL^+_{k+1} + [UL^+,UL^+]\big)
$$
there are no homogeneous component of degree less then $kh(n)$.
Therefore if $k > d/h(n)$ then $UL^+_{(k)} \subset [UL^+,UL^+]$,
which implies that $UL_k = UL_{k+1} = \dots = UL_{\infty}$.
\end{proof}

\begin{reme}
\label{gen:co} Theorem~\ref{gen:maincom} gives that if $d\leq n(n-1)$ then
$U_{\gen,0}$ is generated by the images of the elements $\delta_k$
for $k=1,\dots,d-1$ and $\delta_-$ as an $\sll$ module.
In that case we have $UL_1[V,d] = UL_\infty[V,d]$ and all other
modules $U_{\gen,k}[V,d]$ are trivial for $k\geq 1$. The Conjecture~\ref{gen:conj}
implies that all modules $U_{\gen,k}$ for $k \geq 1$ are trivial.
\end{reme}

\begin{df}
We need to define other modules, which will take part in the
description of the tower of restricted derivations of $L[V,d]$.
Let us denote by
$$
CUL_k^\lambda[V,d] = UL_{k}^\lambda[V,d] \cap \Cen(UL[V,d]),
$$
and
$$
CU_{\gen,k}^\lambda[V,d] = CUL_k^\lambda[V,d]/CUL_{k+1}^\lambda[V,d]
$$
the intersections of the
modules $UL_k^\lambda$ and $U_{\gen,k}^\lambda$ with the center of
the algebra $UL[V,d]$. These modules are ideals in the algebras
$UL$ and $\sll+ UL$ because they lie in the center of $UL$.
\end{df}

Later when we construct the derivation tower of $UL$, we will see that the
spaces $\Hom_\sll(U_{\gen,k}^\lambda,CUL_k^\lambda)$ play an important role. It is
important to notice that there is a natural projection from that space to
$\End_\sll(CU_{\gen,k}^\lambda)$,
because $CU_{\gen,k}^\lambda$ is naturally a subspace of the domain and
factor space of the image.

\begin{reme}
\label{HomasAlgebra}
We can also put an associative algebra
structure on the spaces
$\Hom_\sll(U_{\gen,k}^\lambda,CUL_k^\lambda)$, by defining the
product of two maps $fg$ to be their composition thought the space
$CU_{\gen,k}^\lambda$, i.e.,
$$
(fg)(x) = f(i(\pi(g(x))),
$$
where $\pi$ is the projection from $CUL_k^\lambda$ to
$CU_{\gen,k}^\lambda$ and $i$ is the inclusion of
$CU_{\gen,k}^\lambda$ into $U_{\gen,k}^\lambda$. Notice that the
associative algebra $\Hom_\sll(U_{\gen,k}^\lambda,CUL_k^\lambda)$ does not
have unit unless $U_{\gen,k}^\lambda=CUL_k^\lambda$.
\end{reme}

\end{section}

\ifnewpage

\begin{section}{Normalizer of the ideal of inner derivations}
\label{sec:normalize}

Let us introduce a notation for the normalizer of the ideal of
inner derivations in $UL[V,d]$. We will show that it coincides
with the ideal of inner derivations if $n\geq 3$. We  will need
this ideal in the description of the second derivation algebra of
$L[V,d]$ (see section~\ref{sec:second}).

\begin{df}
\label{nor:centralizers}
Let $\widetilde{UL}[V,d]$ be the nilpotent part of
the algebra of outer derivations of $L[V,d]$, i.e.\
$\widetilde{UL}[V,d] = UL[V,d]/I_L$. Let
$\widetilde{C} [V,d]$ be the center of the algebra
$\widetilde{UL}[V,d]$ and let
$C[V,d] = \pi^{-1}(\widetilde{C}[V,d])$ be its pre-image
in $UL$, where $\pi$ is the natural projection
$\pi: UL[V,d] \to \widetilde{UL}[V,d]$.
We have the inclusion
$$
\Cen(UL[V,d]) + I_L \subset C[V,d].
$$
Let $S[V,d]$ be a  $\sll$-submodule in $C[V,d]$ which is
the complement to the module $\Cen(UL[V,d]) + I_L$.
Similarly let $C_0[V,d]$ and $S_0[V,d]$ be the corresponding
modules, if we use the subalgebra $UL_0[V,d]$
instead of the whole algebra $UL[V,d]$.
\end{df}

\begin{lm}
\label{nor:normadx2}
Let $M$ denote the $\sll$-module
generated by the derivation $\AD{x_1^2}$.
Then,
$$
\Cen_{\widetilde{UL}[V,d]}(M +I_L) =
\{ \delta \in UL[V,d] \mid [\delta,u]\in I_L \mbox{ for all }u\in M\}=
$$
$$
\{ \delta \in UL[V,d] \mid [\delta,u]=0 \mbox{ for all }u\in M\}=
UL[V,d]^{(d-2)}+UL[V,d]^{(d-1)} + I_L.
$$
\end{lm}
\begin{proof}
In order to prove the lemma it is enough to show that
$$
\Cen_{\widetilde{UL}[V]}(M +I_L) = I_L.
$$

Suppose that $\AD_f \in \Cen_{\widetilde{UL}[V]}(M +I_L)$ for
some associative polynomial $f$. Let $k$ be the smallest number such
that $f$ can be written as a sum of products of no more than $k$
Lie polynomials (let us denote this number by $\mathrm{L}\deg f$).
By Lemma~\ref{com:com}, we have $[\AD_f,\AD_{x^2}] = \AD_h$,
where
$$
h=x^2f-fx^2 + \AD_f(x^2) + \AD_{x^2}(f).
$$
The assumption $\AD_f \in \Cen_{\widetilde{UL}[V,\infty]}(M +I_L)$ implies that
$h$ is a Lie polynomial. But it is easy to see that
$$
\mathrm{L}\deg (x^2f-fx^2) = k+1,
\quad
\mathrm{L}\deg \AD_f(x^2) = 2,\mbox{ and }
\quad
\mathrm{L}\deg \AD_{x^2}(f)=k.
$$
This shows that if $k\geq 2$ then $\mathrm{L}\deg h = k+1 >1$ which is
a contradiction. Therefore, $k=1$, i.e.\ $f$ is a Lie polynomial and
$\AD_f\in I_L$.

Suppose that $\DEL{y}{g} \in \Cen_{\widetilde{UL}[V,d]}(M +I_L)$ for
some Lie polynomial $g$, which does not depend on $y$. Then, direct
computation shows that
$$
[\DEL{y}{g},\AD_{x^2}] = \DEL{y}{h}
$$
where $h= \AD_{x^2}(g) - [x,[x,g]]$ also does not depend on $x$. But
this derivation can not be in $I_L$ unless $h=0$, which is
impossible since $\deg g >1$.

By Theorem~\ref{com:split} the Lie algebra $UL[V]$ is generated as an $\sll$-module
by the elements $\AD_f$ and $\DEL{y}{g}$, we have shown that
$$
\Cen_{\widetilde{UL}[V,\infty]}(M +I_L) = I_L.
$$
\end{proof}
\begin{lm}
\label{nor:normnoass}
Let $n\geq 3$ and let $M$ denote the $\sll$-module
generated by the derivation $\DEL{x_1}{[x_2,x_3]}$.
Then,
$$
\Cen_{\widetilde{UL}[V,d]}(M +I_L) =
\{ u \in UL[V,d] \mid [\delta,u]\in I_L \mbox{ for all }u\in M\}=
$$
$$
\{ u \in UL[V,d] \mid [\delta,u]=0 \mbox{ for all }u\in M\}=
\Cen(UL[V,d]) = UL[V,d]^{(d-1)} + I_L.
$$
\end{lm}
\begin{proof}
Similar to the proof of the previous lemma.
\end{proof}

\begin{thq}
\label{nor:descriptionofS}
The inclusion in Definition~\ref{nor:centralizers}
is in fact an equality if $n\geq 3$, i.e.\
$S[V,d] = S_0[V,d]=0$.
In the case $n=2$ we have that the two spaces are equal
$$
S[V,d] = S_0[V,d] =
UL[V,d]^{(d-2)}/(I_L\cap UL[V,d]^{(d-2)})=
(UL^{(d-2)}+I_L)/I_L.
$$
\end{thq}

\end{section}

\ifnewpage

\begin{section}{On Second Derivation Algebra}
\label{sec:second}

In this section we prove Theorems~\ref{2nd:description} and~\ref{2nd:withcenter},
which describe the derivation algebra of
$\D_0 L[V,d]$ and its algebra of restricted derivations
In the next two sections we will use
these theorems as
a base for an induction to describe the
tower of restricted derivations of the algebra $L[V,d]$

Let us consider the algebras $\D U$ and $\D_0 U$
of derivations of
$U= \D_0 L[V,d]$, where $d \geq 2$.

\begin{lm}
\label{2nd:center}
The center of the Lie algebra $U$ is
$$
\Cen(U) = (\Cen(UL[V,d]))^{\sll} = (UL[V,d]^{(d-1)})^{\sll}.
$$
In particular $\Cen(U)=0$ if $n\not|d-1$.
\end{lm}
\begin{proof}
The first equality holds because $U=\sll + UL[V,d]$ and
$\sll$ is a simple Lie algebra. The
second equality holds
because $\Cen_{UL[V]}(\ad(V)) = 0$ (which follows from
Lemma~\ref{com:com}).
\end{proof}

The algebra of the inner derivations of $U$ is isomorphic to
$U/\Cen(U)$, so if we want to describe the derivation algebra of
$U$ we need only to describe the algebra of outer derivations. We
start with any derivation $D$ and we try to modify $D$ using inner
derivation to make $D$ act trivially on as big subalgebra of $U$
as possible.

\begin{df}
Let $\sim$ be the equivalence relation on $\D U$ defined as
follows: $a \sim b $ if and only if the derivation $a-b$ is an
inner derivation of $U$.
\end{df}

\begin{lm}
Any derivation is equivalent to one which kills $\sll$.
\end{lm}
\begin{proof}
The algebra $\sll$ is simple and it is a standard fact in theory
of semi-simple algebras that $H^1(\mr{g},K) = 0$ for any
simple Lie algebra $\mr{g}$ and any $\mr{g}$-module
$K$.
\end{proof}

\begin{co}
\label{2nd:trivialsl}
The Lie algebra $\sll$ acts trivially on the algebra of outer
derivations $\Out U$.
Another way to prove that, is to use the fact that
the inner derivations form an ideal, which contains $\sll$.
Therefore this ideal contains all non trivial $\sll$ submodules
of $\D U$, i.e., the space $\Out U$ is a trivial $\sll$ module.
Therefore every element from this space acts the algebra $\sll$.
\end{co}

First we will describe the space of derivations which act trivially on
all inner derivations $I_L$. This is equivalent to saying that 
we are interested only in derivations which act trivially on the set
$\ad(V)$.

\begin{lm}
\label{2nd:trivialadx} Let $D$ be a derivation of $U$ such that
$D(\ad(x)) \in \Cen(UL)$ for all $x\in V$. Then for any $d \in U$
and any $y\in V$, we have that $(D d)(y) \in \Cen (L[V,d])$, i.e.,
$\im D\in \Cen(UL[V,d])$. The last inclusion it immediately
implies that $D$ acts trivially on $[UL,UL]$.
\end{lm}
\begin{proof}
The algebra $L[V,d]$ is generated by $V$. Therefore, any
such derivation $d$ acts trivially on $\ad(L[V,d])= I_L$.
Finally we can apply Lemma~\ref{der:acttrivally} and obtain
$$
\im D \in \Cen_U(I_L) = \Cen(UL[V,d]).
$$
\end{proof}

\begin{co}
\label{2nd:triviallyonV}
The set of derivations of $U$, which
send $\ad (V)$ to $\Cen(UL[v,d])$ and act trivially on $\sll$,
can be embedded into the set
$$
\Hom_{\sll}\big(U_{\gen,0}[V,d], \Cen(UL[V,d])\big)
+ \Hom_{\sll}\big( U_{\gen}^{\triv}[V,d], \Cen(UL[V,d])^{\sll} \big)
.
$$
Every element in this set can be extended to a derivation of $\D_0
L[V,d]$, by letting it act trivially on $\sll$ and $[UL,UL]$.
Therefore the above embedding is an isomorphism. This is also an
isomorphism of Lie algebras (the Lie algebra structure of
$\Hom_\sll(\dots)$ comes form the associative algebra structure on
this space described in Remark~\ref{HomasAlgebra}).
\end{co}

Our next step is to see how the derivations of $U$ can act on
the space $\ad(V)$.

\begin{lm}
\label{2nd:adVpossible}
For any derivation $D$ of the algebra $UL$ we have
$$
D(\ad(V)) \subset UL^{(1)} + C[V,d].
$$
\end{lm}
\begin{proof}
Without loss of generality we may assume that $D$ is
homogeneous derivation of positive degree
(if the degree is zero we have that
$D(\ad(V)) \subset UL^{(1)} $). Since the
space $I_L$ is an ideal we have that
$D([I_L,I_L]) \subset I_L$.
Let $\delta \in UL$ be a derivation.
Then for any $x\in V$ we have that
$$
D([\ad(x),\delta]) = [D(\ad(x)),\delta] + [\ad(x),D(\delta)],
$$
but $[\ad(x),D(\delta)]\in I_L$ and
$$
D([\ad(x),\delta]) = D(\ad(\delta(x))) \in I_L.
$$
Therefore, we have that $[D(\ad(x)),\delta] \in I_L$, which
is equivalent to $D(\ad(x))\in C[V,d]$.
\end{proof}

\begin{lm}
\label{2nd:adVinner}
For any derivation $D$ of the algebra $UL$, such that
$D(\ad(V))\in [I_L,I_L]$, there exists an element
$f\in UL$ such that $(D-\ad(f))(\ad(V))=0$.
\end{lm}
\begin{proof}
If $D$ is a derivation such that $D(\ad(V))\in [I_L,I_L]$, then
$D$ defines a map $d : V \to [L[V],L[V]]$ by
$\ad(d(x)) = D(\ad(x))$. We can extend this map to a derivation
$f$ of $L[V]$. It is easy to see that $f\in UL$ and
$(D-\ad(f))(\ad(V))=0$.
\end{proof}

\begin{lm}
\label{2nd:adVouter}
Let $f\in \Hom_{\sll}(V,S[V,d])$ be an $\sll$-invariant
linear map. Then there exists a derivation $D$ of
$\D_0 L[V,d]$ such that $D(\ad(x)) = f(x)$ for all $x\in V$.
\end{lm}
\begin{proof}
The map $f$ can be extended to a derivation
from $I_L$ to $UL[V,d]$ by the Leibnitz rule
(because $I_L$ is the `free' Lie algebra generated
by $V$). Let us denote this extension with $D$.
Note that $D([I_L,I_L])\subset I_L$.

Now we have to define $D(\delta)$ for any $\delta\in UL$.
If such an extension exists we would have
$$
D(\ad(\delta(x))) = D([\ad(x),\delta]) =
[D(\ad(x)),\delta] + [\ad(x),D(\delta)].
$$
Therefore
$$
\ad D(\delta)(x) = [f(x),\delta] - D(\ad(\delta(x))).
$$

By the definition of the module $S$ we have that
$ [f(x),\delta] \in I_L$. Also we have that
$D(\ad(\delta(x)))\in I_L$, because
$\ad(\delta(x))\in [I_L,I_L]$. Therefore,
the right side
is an element in $I_L$, which allows us
to define $D(\delta)(x)$ by the above equality.
Thus, we have defined $D(\delta)$. It is easy to
check that $D$ is a derivation of $UL$. Finally
setting $D(\sll)=0$, defines $D$ as a derivation
of $\D_0 L[V,d]$.
\end{proof}
\begin{reme}
Note that the statement of this lemma is trivial unless
$n=2$ and $2|d$, because if $n>2$ by,
Corollary~\ref{nor:descriptionofS}, $S=0$ and if $n=2$ and
$2\not\!|\,\, d$, then $\Hom_{\sll}(V,S) = 0$.
Also using the fact that the space $S$ is so
small there is no need to check that the map
defined in the proof above is a derivation because it
sends `almost' every element to $0$.
\end{reme}

Now we can describe the derivation algebra of $U$.

\begin{thq}
\label{2nd:trivcenter}
If $\Cen(U)=0$ then the derivation algebra of $U$ can be written as
$$
\D U = \gll +  UL + \sum_{\lambda} W^\lambda_0 + T,
$$
where
$$
W^{\lambda}_0[V,d] = \Hom_{\sll}(U_{\gen,0}^{\lambda}, CUL^{\lambda}_0[V,d]),
$$
for $\lambda \not = [1]$ and
$$
W^{[1]}_0[V,d] = \Hom_{\sll}(U_{\gen,0}^{[1]}, CUL^{[1]}_0[V,d])/(UL^{(d-2)})^\sll.
$$
Finally we have that
$$
T[V,d] \simeq \Hom_{\sll} \big(\ad (V), S[V,d] \big)=
$$
$$
\left\{
\begin{array}{ll}
0 & \mbox{if } n\geq 3 \mbox{ or }d \leq 3 \\
\Hom_{\sll}(\ad(V), UL[V,d]^{(d-2)}+I_L/I_L)
 & \mbox{if } n=2  \mbox{ and } d\geq 4\\
\end{array}
\right. .
$$
The Lie algebra structure on this space is given by the natural
action of $\gll$ on $UL$, and its trivial action on $W^\lambda_0$
and $T$; also $W^\lambda_0$ and $T$ act naturally $UL$ in
particular they act trivially on $[UL,UL]$.
\end{thq}
\begin{reme}
The space $T$ is trivial, because if $n\geq 3$ then $S$ is trivial.
In the case $n=2$, the condition $\Cen(U)=0$ implies that $2 | d$
therefore the space $\Hom_\sll(\ad(V),UL^{(d-2)})$ is trivial because
the source has degree $1$; the target has even degree and
all $\sll$ invariant maps are of even degree.
\end{reme}
\begin{proof}
The algebra of inner derivations of the $U$ is the part $\sll + UL$,
because the center of $U$ is trivial.

By Lemma~\ref{2nd:trivialadx} and
Corollary~\ref{2nd:triviallyonV} we have that the space $W$ of
outer
derivations (which acts `almost' trivially on $\ad(V)$) is
isomorphic to
$$
\Hom_{\sll}\left(U_{\gen,0}[V,d],\Cen(UL[V,d])\right),
$$
since the part coming from $U^{\triv}_{\gen}$ is trivial because
$\Cen(U) =0$. Lemmas~\ref{2nd:adVpossible},~\ref{2nd:adVinner}
and~\ref{2nd:adVouter} describe the set of derivations which act
`non-trivially' in $\ad(V)$ modulo the inner derivations. We need
to factor out some part of the space $W^{[1]}_0$ because
$(UL^{(d-2)})^\sll$ is part of the algebra of inner derivations and
also part of $\Hom_\sll(V,\Cen(UL))$.
\end{proof}

Now we can describe the second algebra in the
restricted derivation tower of $L[V,d]$.
\begin{thq}
\label{2nd:description}
The algebra $\D_0^2 L[V,d]$ can be written as
$$
\D_0^2 L[V,d] = \sll + UL[V,d] + \sum_\lambda R^\lambda_0[V,d] + T
$$
provided that $\Cen(\D_0L[V,d])=0$. Here we have denoted
$$
R^{\lambda}_0[V,d] = \Hom_{\sll}^0(U_{\gen,0}^{\lambda}, CUL^{\lambda}_0[V,d]),
$$
for $\lambda \not = [1]$ and
$$
R^{[1]}_0[V,d] = \Hom_{\sll}^0(U_{\gen,0}^{[1]}, CUL^{[1]}_0[V,d])/(UL^{(d-2)})^\sll.
$$
\end{thq}
\begin{proof}
Theorem~\ref{2nd:trivcenter} 
gives a
description of the full derivation algebra. Here we
are interested only in restricted derivations, which
correspond to putting $\Hom^0$ instead of $\Hom$.
\end{proof}

\begin{reme}
\label{faithful2}
The difference between the reductive part of $\D_0^2 L$ and $\D_0 L$
is just a sum of copies of $\sll_l$, for some $l$-es. It acts faithfully on the
nilpotent part and all
nontrivial modules for some $\sll_l$ are isomorphic to the standard
or to its dual.
Therefore the reductive part of $\D_0^2 L$ is a semi simple Lie algebra which
acts faithfully on the nilpotent part.
\end{reme}

\begin{thq}
\label{2nd:small-stabilize}
The restricted derivation tower of the free
nilpotent Lie algebra of class $d \leq n$ terminates at
the second level,
i.e.\
$\D^2_0 L[V,d] = \D^1_0 L[V,d]$.
\end{thq}
\begin{proof}
Theorem~\ref{2nd:description} describes the algebra
$\D^2_0 L[V,d]$.
In this case we have
$$
R_0^\lambda = \Hom(U_{\gen,0}^{\lambda}, CUL^{\lambda}_0[V,d])
$$
but this space is trivial unless $|\lambda| = d-1$, because we do
not have $\sll$ invariant maps except when $n$ divides the
difference of the degree. If $|\lambda| = d-1$, then by Theorem~\ref{gen:smalldeg},
$U_{\gen,0}^{\lambda}$ is trivial unless for $\lambda = [d-1]$. In the last
case we have that both $U_{\gen,0}^{\lambda}$ and
$CUL^{\lambda}_0[V,d]$ coincide and are simple $\sll$ modules
therefore $\Hom_\sll$ is one dimensional, but $\Hom^0_{\sll}$
is trivial.
\end{proof}

In the case when the center of $U$ is not trivial the situation is
slightly more complicated:

\begin{thq}
\label{2nd:withcenter}
If $\Cen(U) \not = 0$  then the derivation
algebra of $U$ decomposes
$$
\D U = \gll +  UL/\Cen(U) + \sum_{\lambda} W^\lambda_0 + T
$$
and the algebra of restricted derivations is
$$
\D_0^2 L[V,d] = \sll + UL[V,d]/\Cen(\D_0 L[V,d]) + \sum_\lambda R^\lambda_0[V,d] + T
$$
where $W^\lambda_0$ and $R^\lambda_0$ for $\lambda \not = \triv$
are defined as in Theorems~\ref{2nd:trivcenter}
and~\ref{2nd:description}, and
$$
W^{\triv}_0 = \Hom_{\sll}(U_{\gen}^{\triv}, \Cen U)
\quad \mbox{and} \quad
R^{\triv}_0 = \Hom_{\sll}^0(U_{\gen}^{\triv}, \Cen U).
$$
Also note that $W^{\triv}_0$ and $R^{\triv}_0$ split as a direct
summands of $\D \D_0 L[V,d]$ and $\D_0^2 L[V,d]$ respectively.
\end{thq}
\begin{proof}
The proof of this theorem is similar to the one of
Theorem~\ref{2nd:trivcenter}. The only thing that needs to be
checked is that $W^{\triv}_0$ acts trivially on the rest of the
algebra $\D U$.
This is true  because $W^{\triv}_0(U) \subset \Cen (U)$ and
everything else acts trivially on the center of $U$.
Therefore it splits as a direct summand.
\end{proof}

\begin{lm}
The centralizer of the ideal $\ad(UL_0)$ in $\D_0^2 L[V,d]$
is isomorphic to $W_0^{\triv}[V,d] + \Cen(UL[V,d])/\Cen(\D_0 L)$.
\end{lm}
\begin{proof}
The centralizer of $\ad(UL_0)$ in
$UL[V,d]/\Cen(\D_0 L)$ is
$$
\Cen(UL[V,d])/\Cen(\D_0 L),
$$
and for $\lambda\not=\triv$, the space $W_0^{\lambda}$ acts
faithfully on $U_{\gen}^{\lambda}\subset UL_0$. Therefore the
centralizer of $\ad(UL_0)$ in $\D_0^2 L[V,d]$ is contained in
$$
W_0^{\triv}[V,d] + \Cen(UL[V,d])/\Cen(\D_0 L) + T[V,d].
$$
It is easy to see that it is only
$W_{0}^{\triv}[V,d] + \Cen(UL[V,d])/\Cen(\D_0 L)$, because $T$ act faithfully on
$\ad(V)\subset UL_0$.
\end{proof}

\begin{lm}
The center of the algebra $\D_0^2 L[V,d]$ coincides with the
center of $R^{\triv}_0$.
\end{lm}
\begin{proof}
By the previous lemma we have that the center of $\D_0^2 L[V,d]$
lies inside $R^{\triv}_0$, because $\sll$ acts nontrivially on
$\Cen(UL)/\Cen(U)$. Using the algebra $\D_0^2 L[V,d]$ splits as a
direct sum, it is easy to see that its center is exactly $\Cen
R^{\triv}_0$.
\end{proof}
\begin{reme}
The analog of the Remark~\ref{faithful2} holds also in the case when the center of
$U$ is not trivial.
\end{reme}
\end{section}

\ifnewpage

\begin{section}[Derivation Tower. Trivial Center]{Derivation tower of free nilpotent Lie algebras.
Case with trivial center}
\label{sec:dertower}

The description of the second derivation algebra can be
generalized to a description of all algebras in the restricted
derivation tower. As we saw in the previous section there is a big
difference in the structure of the algebra $\D_0^2 L[V,d]$
depending on the triviality of the center of $\D_0 L[V,d]$.

Let us first consider the case when $n \not | d-1$ and the algebra
$\D_0 L[V,d]$ has a trivial center. It is easy to see that this
implies that all algebras in the restricted derivation tower have
trivial centers, and all the maps $\D_0^k L[V,d] \to \D_0^{k+1}
L[V,d]$ are embeddings.

The next theorem gives the description of all algebras in the
restricted derivation tower and its proof is by induction where
Theorem~\ref{2nd:trivcenter} serves as a base case.

\begin{thq}
\label{dertower:nocenter}
The derivation algebra of $\D_0^k L[V,d]$, for $k\geq 2$ is
$$
\D \D_0^k L[V,d] = \D_0^k L[V,d] + \sum_\lambda W_{k-1}^\lambda,
$$
where
$$
W_{k}^\lambda =  \Hom_\sll(U_{\gen,k}^{\lambda}, CUL_{k}^\lambda).
$$
The algebra of restricted derivations is obtained by replacing $W_{k-1}^\lambda$ with
its subspace $R_{k-1}^\lambda$, consisting of maps with zero trace on
$CU_{\gen,k-1}^\lambda$.
\end{thq}
\begin{proof}
Proof is by induction on $k$. Suppose that the theorem is true for
some $k$ (in case $k=2$ we use Theorem~\ref{2nd:description}).
In order to make
the induction step we have to describe the
algebras $\D \D_0^k L[V,d]$ and $\D_0 \D_0^k L[V,d]$.

Let us
describe the derivation algebra
$\D \D_0^k L[V,d]$. First we want to see how the outer derivations
can act on the algebra $\D_0^k L[V,d]$.

\begin{lm}
\label{nil:kill1}
Let $D\in \D \D_0^k L[V,d]$ be a derivation of
$\D_0^k L[V,d]$ such that $D(\ad(V)) = \Cen(UL)$ and $D(\sll)=0$. Then
$D([UL^+,UL^+])=0$.
\end{lm}
\begin{proof}
By Lemma~\ref{der:Cenofideal},
$D(f)$ commutes with $\ad(V)$ for any
$f\in UL$. Note that the space of elements in
$\D^k_0 L[V,d]$ which commute with $\ad(V)$ is
$$
\Cen_{\D^k_0 L[V,d]}(\ad(V)) =
\Cen(UL) + \sum_{i=0}^{k-2} \sum_{\lambda} R_i^\lambda[V,d].
$$
All spaces $R_{i}^{\lambda}[V,d]$ are trivial $\sll$-modules.
Therefore $D(U_{\mu}) \subset \Cen(UL)$ for all simple modules
$U_{\mu}$ corresponding to the partition $\mu \not = \triv$. The
space $\Cen(UL)$ is an abelian Lie algebra.  Therefore $D$ acts
trivially on the commutator algebra of the algebra generated by
$UL^+$.
\end{proof}

\begin{lm}
\label{nil:posdegree}
Let $D\in \D \D_0^k L[V,d]$ be a derivation of
$\D_0^k L[V,d]$ such that $D(\ad(V)) =0$ and $D(\sll)=0$. Then
$$
D(x)\in \Cen(UL) + \sum_{i=0}^{k-2} \sum_{\lambda} P_i^\lambda[V,d],
$$
for all $x\in UL[V,d]$, where $P_i^\lambda$ is the subspace of
$R_i^\lambda$ consisting of all maps which are trivial on
$CU_{\gen,i}^\lambda$.
\end{lm}

\begin{proof}
Assume the contrary --- then there exists $g\in \Cen(UL)_\lambda$
such that $D(x)(g)\not=0$. Therefore, $[D(x),g] = D(x)(g)\not=0$.
But this is impossible because $D(g)=0$ and $[x,g]=0$.
\end{proof}

\begin{lm}
\label{nil:killULk}
Let $D\in \D \D_0^k L[V,d]$ be a derivation of
$\D_0^k L[V,d]$ such that $D(\ad(V)) =0$ and $D(\sll)=0$. Then
$D(UL_{k}[V,d])=0$.
\end{lm}
\begin{proof}
Let $u$ be an element in
$u \in UL_{(i)}^\lambda$ for some
$i\geq k$, which can be written as
$u=[\tilde u, u_0]$, where
$\tilde u \in UL^{\lambda}_{(i-1)}$
and $u_0\in U_{\triv}$. We have
$$
D(u) = [D(\tilde u), u_0] + [\tilde u, D(u_0)].
$$
Using $\sll$ invariance we have
$D(\tilde u)\in \Cen(UL)$ and, therefore
$D(\tilde u)$ commutes with $u_0$.
Also
$D(u_0) \in \sum_{j\leq k, \lambda} R_{j}^{\lambda}$ and by
the induction hypothesis acts
trivially on $UL_{k-1}^0[V,d]$.
Since $\tilde u \in UL_{k-1}[V,d]$, we have that $D(u)=0$.

Finally, we notice that any element in
$U^{\lambda}_{(i)}$ is a sum of the elements considered
above, which completes the proof.
\end{proof}

\begin{lm}
\label{nil:killfinal}
Let $D\in \D \D_0^k L[V,d]$ be a derivation of
$\D_0^k L[V,d]$ such that $D(\ad(V)) =0$ and $D(\sll)=0$. Then
$
D(UL^{\lambda}_{k-1}) \subset CU^{\lambda}_{k-1}
$,
for any partition $\lambda \not = \triv$.
\end{lm}
\begin{proof}
Assume the contrary, ie., that $D(x) \in CUL^\lambda_{i-1}
\setminus CUL^\lambda_{i}$, for some $x\in UL^{\lambda}_{k-1}$ and
some $i < k$. Then there exists some $f\in R^\lambda_{i-1}$ such
that $[f,D(x)]\in CUL^\lambda \setminus 0$. By the choice of $f$
we have $[f,x] =0$, which implies that $[D(f),x ] \in CUL^\lambda
\setminus 0$, but this is impossible because all $R^\mu_j$ act
trivially on $x$. This is a contradiction which finishes the proof
of the lemma. 
\end{proof}

\begin{lm}
Let $f \in
\sum_{\lambda} \Hom_{\sll}\big(
U_{\gen,k-1}^{\lambda},
CUL^{\lambda}_{k-1} \big)$
be a $\sll$ invariant linear map. Then there exists
a derivation $D \in \D \D_0^k L[V,d]$ such that
the restriction of $D$ to $\sum U_{\gen,k-1}^{\lambda}$
coincides with $f$.
\end{lm}
\begin{proof}
Let us construct a derivation $D$ which acts trivially on $\sll$,
$UL_k$, $T$ and $R_{i}^{\lambda}[V,d]$ for any partition $\lambda$
and all $i\leq k$. This allows us to lift the map $f$ to a map $D$
form $UL_{k-1}$ to $CUL^+_{k-1}$. Let us extend $D$ to an $\sll$
invariant map from $UL_0$ to $CLU^+_{k-1}$. In order to define $D$
as a derivation of $\D_0^k L[v,d]$ we only need to say how it acts
on $UL^{\triv}$ and verify that it satisfies the Leibnitz rule.

Now we define
$D: UL^\triv \to \sum R^\mu_j$ by
$$
D(u)(v) = D[u,v],
$$
for any $u\in UL^{\triv}$ and $v\in UL_0$.

In order to verify that $D$ is well defined, we need to verify that for
any $u\in UL_k \cap UL^{\triv}$, we have that $D(u)=0$. Any such $u$
has to lie in $[UL^+,UL^+]$ and therefore $[u,v] \in [UL^+,UL^+] \subset UL_k$,
This shows that $D(u)=0$, which gives that $D$ is well defined.

Finally we have to check that $D$ satisfies the Leibnitz rule,
which can be verified using the construction of $D$.
\end{proof}

\begin{reme}
Lemmas~\ref{2nd:adVpossible} and~\ref{2nd:adVinner}
also apply in this case, with one exception --- since $UL[V,d]$
is no longer closed under derivations, we have to work with
$UL_0[V,d]$, which is still closed. Therefore, we need to substitute
$S[V,d]$ with $S_0[V,d]$, but Corollary~\ref{nor:descriptionofS}
gives that these two spaces are the same.
Thus, for any $D \in \D \D_0^k L[V,d]$ there exists an inner derivation
$\ad(d)$ and a scalar $c$ such that
$D -d$ acts on $\ad(V)$ as a multiplication by $c$.
\end{reme}

\begin{lm}
\label{nil:trivialonUL}
Let $D\in \D \D_0^k L[V,d]$ be a derivation of
$\D_0^k L[V,d]$ such that $D(UL_0)=0$ and $D(\sll)=0$. Then
$D=0$.
\end{lm}
\begin{proof}
The centralizer of the space $UL_0$ in $\D_0^k L[V,d]$
is $\Cen UL$. Lemma~\ref{der:Cenofideal} gives us that $\im D
\subset \Cen UL$. The last space does not have any $\sll$
invariant elements. Therefore $D$ acts trivially on any trivial
$\sll$ module in $\D_0^k L[V,d]$. By construction the factor
algebra $\D_0^k L[V,d] /UL_0$ is $\sll$ invariant, therefore $D$
is trivial.
\end{proof}

We can combine the results of the previous lemmas to obtain a
description of the derivation algebra of $\D_0^k L[V,d]$.

\begin{thq}
\label{nil:detDoknocenter}
The derivation algebra of $\D_0^k L[V,d]$ can be written as
$$
\textstyle
\D \D_0^k L[V,d] = \D_0^k L[V,d] + \C .d +
\sum_\lambda W_{k-1}^\lambda[V,d] ,
$$
where
$$
W_{k-1}^{\lambda }[V,d] \simeq
\Hom_{\sll}\big( \Hom_\sll(U_{\gen,k-1}^{\lambda}, CUL_{k-1}^\lambda) \big).
$$
\end{thq}
\begin{reme}
Note that
above we showed that the set of derivations from
$\sll+UL_0$ to $\D_0^k L[V,d]$ is isomorphic to
$\D \D_0^k L[V,d]$. We are going to use this fact
in chapter 4.
\end{reme}

Now we can finish the proof of
Theorem~\ref{dertower:nocenter}.
We only have to notice that
in order to pass from $\D\D_0^k L[V,d]$ to
$ \D_0^{k+1} L[V,d]$, we need to substitute all
$\Hom$ operators
with $\Hom^0$'s. After applying this to
$W_{k-1}^\lambda[V,d]$,
we get $R_{k-1}^\lambda[V,d]$, which
completes the proof of the theorem.
\end{proof}

The next corollaries are an immediate consequence of the
above Theorem and Theorem~\ref{ULstabilization}.

\begin{co}
If $n\not | d-1$, then the restricted derivation tower of $L[V,d]$
stabilizes after at most $3+ d/h(n)$ steps.
\end{co}
\begin{co}
If $n<d<n(n-1)$ (or if Conjecture~\ref{gen:conj} holds) then the restricted
derivation tower of $L[V,d]$ stabilizes at the third level if
$n \not | d-1$.
\end{co}
\begin{reme}
\label{faithfulk:nocenter}
The difference between the reductive part of $\D_0^{k+1} L$ and $\D_k L$
is just a sum of copies of $\sll_l$, for some $l$-es. It acts faithfully on the
nilpotent part and all
nontrivial modules for some $\sll_l$ are isomorphic to the standard
or to its dual.
Therefore the reductive part of $\D_0^k L$ is a semi simple Lie algebra
which decomposes as a sum of copies of $\sll_l$ for some $l$-es.
It also acts faithfully on the nilpotent part and almost all nontrivial
modules which appear in are either isomorphic to the standard one or its dual.
\end{reme}

\end{section}

\ifnewpage

\begin{section}[Derivation Tower. Nontrivial Center]{Derivation tower of free nilpotent Lie algebras.
Case with non trivial center}

Now let us consider the case $n | d-1$. This case is similar to
the case when the center is trivial, but there are two important
things to consider. First, the centers of some algebras in the
derivation tower are not trivial, and the map between two algebras
in the tower is not always an embedding. Also all algebras in the tower
split as a direct sum of two subalgebras.

\begin{thq}
\label{dertower:center}
For $k > 2$, the $k$-th algebra
in the restricted derivation tower of $L[V,d]$ can be written as a
direct sum $\D^{k}_0 L[V,d] = DL^k[V,d] \oplus \tilde R_{k}^{\triv}[V,d]$.

The algebras $DL^k[V,d]$ and $\tilde R_{k}^{\triv}[V,d]$ can be
constructed recursively by
$$
DL^{k+1}[V,d] = DL^{k}[V,d] + \sum_{\lambda\not = \triv} R_{k-1}^{\lambda}[V,d]
$$
where the spaces $R_{k-1}^{\lambda}$ of outer derivations for
$\lambda \not = \triv$ are
as in Theorem~\ref{dertower:nocenter}, and
$$
\tilde R_{k}^{\triv}= \D_0  \tilde R_{k-1}^{\triv} +
\Hom\big(
DL^{k-1}/[DL^{k-1},DL^{k-1}],
\Cen(\tilde R_{k-1}^{\triv})
\big).
$$
\end{thq}
\begin{proof}
As the proof of Theorem~\ref{dertower:nocenter}, this theorem is proved on induction.
The base case $k=2$ is given by Theorem~\ref{2nd:withcenter}.

The proof is essentially the same as the proof in the case with trivial center,
the main differences are:

a) In order to prove Lemma~\ref{nil:kill1} we need to use the fact that
$\tilde R_{i}^{\triv}[V,d]$ commutes with $UL_0[V,d]$, which is part in the induction
hypothesis.

b) Lemma~\ref{nil:trivialonUL} does not hold in this case, but it can
be generalized to

\begin{lm}
Let $D$ be a derivation of $\D_0^k L[V,d]$ such that $D(UL_0)=0$ and $D(\sll)=0$,
then $\im D \subset \tilde R_{k}^{\triv}$.
\end{lm}
\begin{proof}
As in Lemma~\ref{nil:trivialonUL} we want to apply Lemma~\ref{der:centralize},
therefore we
need to describe the centralizers of $UL_0$ in
$\D_0^k L[V,d]$. We have that
$$
\Cen_{\D_0^k L[V,d]} (UL_0) = \tilde R_{k}^{\triv} + \Cen (UL)/\Cen(UL)^{\sll},
$$
as before $\Cen (UL)/\Cen(UL)^{\sll}$ does not contain $\sll$ invariant vectors,
and $\sll$ acts trivially on $\D_0^k L[V,d]/UL_0$. Therefore the image of $D$
has to lie in the subalgebra $\tilde R_{k}^{\triv}$.
\end{proof}

\begin{lm}
\label{nil:derin0} The space of derivations from $\D_0^k L[V,d]$
to $\tilde R_{k}^{\triv}$ is
$$
\tilde R_{k+1}^{\triv}= \D \tilde R_{k}^{\triv} +
\Hom\big(DL^{k}/[DL^{k},DL^{k}],\Cen(\tilde R_{k}^{\triv})\big).
$$
\end{lm}
\begin{proof}
Let $L$ be a Lie algebra such that $L=L_1 \oplus L_2$, where
$L_1$ and $L_2$ are ideals. Let $D$ be a derivation of $L$
with image in $L_1$; using the fact that $L_1$ and $L_2$ commute, we can see
that $D$ kills $[L_2,L_2]$, and $D(L_2) \subset \Cen L_1$.
This gives that the space of derivations form $L$ to $L_1$ is isomorphic to
$$
\Hom (L_2/[L_2,L_2],\Cen(L_1)) + \D L_1.
$$
In order to prove the lemma we need to apply the above argument to
$L=\D_0^k L[V,d]$, $L_1=\tilde R_{k}^{\triv}$
 and $L_2= DL^{k}$ .
\end{proof}

As in Theorem~\ref{nil:detDoknocenter}, we can describe
the derivation algebra of $\D_0^k L[V,d]$.

\begin{thq}
\label{nil:detDokcenter}
The derivation algebra of $\D_0^k L[V,d]$ can be written as
$$
\textstyle
\D \D_0^k L[V,d] = \left(
DL^k [V,d] + \C .d + \sum_{\lambda \not = \triv} \widetilde W_{k-1}^{\lambda }[V,d]
\right)
\oplus \tilde W_{k+1}^{\triv }[V,d],
$$
where $\tilde W_{k+1,\lambda }[V,d]$
is defined in Theorem~\ref{nil:detDoknocenter} for $\lambda\not=\triv$,
and in Lemma~\ref{nil:derin0} for $\lambda=\triv$.
\end{thq}

To finish the proof of
Theorem~\ref{dertower:center},
we only need to substitute all
$\Hom$ and $\D$ operators with $\Hom^0$'s and $\D_0$'s, respectively.
\end{proof}

The sequence $DL^k$ stabilizes after finitely many steps because the sequence
$UL_k$ stabilizes. In order to show that the restricted derivation tower
of $L[V,d]$ stabilizes in this case we need to show that the sequence
$\tilde R_k^{\triv}$ stabilizes.

\begin{lm}
\label{bigdim}
If $\tilde R_2^{\triv} \not = 0$ then for all $k\geq 2$
$$
\dim DL^{k}/[DL^{k},DL^{k}] \geq 2.
$$
\end{lm}
\begin{proof}
If $n\geq 3$ is follows from Theorem~\ref{gen:maincom}. Because in
$\Cen UL$ there are many $\sll$ (more than 2) modules of type
$[n+s,s^{n-1}]$ for $s=(d-1)/n-1$, because in the free Lie algebra
$L[V]$ there are many modules of type $[n+s+1,s^{n-1}]$ (here we
are using the $d> h(n)\geq n(n-1)$). The same argument works in
the case $n=2$ if $d\geq7$. Therefore the only case which is left
is $n=2$ and $d=5$. In this case the direct computation shows that
$\dim DL^{k}/[DL^{k},DL^{k}] =2$ for all $k\geq 2$.
\end{proof}

\begin{lm}
\label{nil:Wsequence}
The sequence
$$
\tilde R_2^{\triv} \to
\tilde R_3^{\triv} \to
\tilde R_4^{\triv} \to
\tilde R_5^{\triv} \to
\cdots
$$
stabilizes at most at the sixth level, i.e., $\tilde
R_{6}^{\triv}= \tilde R_{7}^{\triv}$.
\end{lm}
\begin{proof}
The behavior of this sequence is determined by
the dimensions of the spaces
$U_{\gen}^{\triv}$, $\Cen(UL)^{\sll}$ and $\Cen(UL)^{\sll} \cap [UL,UL]$.
Let us define $p,q,r$ as follows
$$
p = \dim (UL / (UL^+ + \Cen(UL) + [UL,UL]), \quad
q = \dim (\Cen(UL)^{\sll} \cap [UL,UL]) \mbox{ and}
$$
$$
r = \dim (\Cen (UL)^{\sll} / \Cen(UL)^{\sll} \cap [UL,UL]).
$$
Let also denote $l_i= \dim(DL^i/[DL^i,DL^i])$, by Lemma~\ref{bigdim} we have
that $l_i\geq 2$ for all $i$ unless $\tilde R_2^{\triv} =0$.
Let us consider the different possibilities for this dimensions:
\begin{enumerate}
\item
$p=r=0$ or $r=q=0$.
In this case all spaces $\tilde R_i^{\triv}$, for $i\geq 2$ are trivial and
the sequence stabilizes.

\item $r=0$, $p=q=1$. In this case we case we have $\tilde
R_2^{\triv} = \R$ is abelian, this implies that $\tilde
R_3^{\triv} = \R^{l_2}$ is also abelian and $\tilde R_4^{\triv} =
\sll_{l_2} + \R^{l_2l_3}$ and finally $\tilde R_i^{\triv} =
\sll_{l_2} + \R^{l_2l_3} + \sll_{l_3}$ for $i\geq 5$.

\item $r=0$ and $p.q > 1$ In this case we case we have $\tilde
R_2^{\triv} = \R^{pq}$ is abelian, therefore $\tilde R_3^{\triv} =
\sll_{pq} + \R^{pql_2}$ and finally $\tilde R_i^{\triv} =
\sll_{pq} + \R^{pql_2} + \sll_{l_2}$ for $i\geq 4$.

\item
$r=1$ and $p=q=0$.
In this case all spaces $\tilde R_i^{\triv}$, for $i\geq 2$ are trivial and
the sequence stabilizes.

\item $r=1$ and $p=0$ or $q=0$( say $q=0$). $\tilde R_2^{\triv} = \R^{p}$ is
abelian is in cases 2 and 3, the sequence $\tilde R_i^{\triv}$
stabilizes after 5 or 6 steps depending on whether $p=1$ or not.

\item $r=1$ and $p=q=1$. In this case $\tilde R_2^{\triv}$ is
three dimensional nilpotent of class 2 Lie algebra. Therefore we
have
$\tilde R_3^{\triv} = \sll_2 + \R^2 \oplus \R^{l_2}$.
Similarly to case 3 the sequences $\tilde R_i^{\triv}$ stabilizes
after 5 steps.

\item
$r=1$ and $p=1$ or $q=1$.
$\tilde R_2^{\triv} = \R + \R^p + \R^p$, where one copy of $\R^p$ lies in the center and
commutation with elements in $\R$ map one copy of $\R^p$ to the other. After that we have
$\tilde R_3^{\triv} = \R + \R^p + \sll_p + \R^{pl_2}$,
$\tilde R_4^{\triv} = \R^p + \sll_p + \R^{pl_2} + \sll_{l_2} + \R^{l_3}$ and
$\tilde R_5^{\triv} = \R^p + \sll_p + \R^{pl_2} + \sll_{l_2} + \sll_{l_3} + \R^{l_3l_4} $ and finally
$$
\tilde R_6^{\triv} = \R^p + \sll_p + \R^{pl_2} + \sll_{l_2} + \sll_{l_3} + \R^{l_3l_4} + \sll_{l_4}
$$
for all $i \geq 5$

\item
$r=1$ and $p,q>1$.
We have that
$\tilde R_2^{\triv} = \R^p + \R^q + \R^{pq}$, where $\R^{pq}$ is in the center and
the commutator of $\R^p$ and $\R^q$ is $\R^{pq}$. By the construction of
$\tilde R^{\triv}_i$ we have
$$
\tilde R^{\triv}_3 = \sll_p + \sll_q + \R^p + \R^q + \R^{pql_2}
$$
and for $i\geq 4$ we have
$$
\tilde R^{\triv}_i = \sll_p + \sll_q + \R^p + \R^q + \R^{pql_2} + \sll_{l_2}.
$$

\item
$r\geq 2$ and $p=q=0$.
$\tilde R_i^{\triv} = \sll_r$ for all $i\geq 2$.

\item
$r\geq 2$ and $p=0$ or $q=0$.
$\tilde R_2^{\triv} = \sll_r + \R^{r(p+q)}$. If $p+q=1$ then all
$\tilde R_i^{\triv}$ are equal to $\tilde R_2^{\triv}$, otherwise we have
$\tilde R_i^{\triv} = \sll_r + \R^{r(p+q)} + \sll_{p+q}$ for all $i\geq 3$.

\item
$r\geq 2$ and $p=q=1$.
$\tilde R_2^{\triv} = \sll_r + \R^{r} + \R^r + \R$, where $\R$ is in the center.
$\tilde R_3^{\triv} = \sll_r + \R^{r} + \R^r + \R^{l_2}$, and
$\tilde R_4^{\triv} = \sll_r + \R^{r} + \R^r + \sll_{l_2}+\R^{l_2l_3}$, and finally
$$
\tilde R_i^{\triv} = \sll_r + \R^{r} + \R^r + \sll_{l_2}+\R^{l_2l_3} + \sll_{l_3},
$$
for all $i \leq 5$.

\item
$r\geq 2$ and $p=1$ or $q=1$.
$\tilde R_2^{\triv} = \sll_r + \R^{r} + \R^{rp} + \R^p$, where $\R^p$ is in the center.
$\tilde R_3^{\triv} = \sll_r + \R^{r} + \R^{rp} + \sll_p + \R^{pl_2}$, and finally
$$
\tilde R_i^{\triv} = \sll_r + \R^{r} + \R^{rp} + \sll_p + \sll_{l_2}+\R^{pl_2},
$$
for all $i \leq 4$.

\item
$r\geq 2$ and $p,q>1$.
$\tilde R_2^{\triv} = \sll_r + \R^{rq} + \R^{rp} + \R^{pq}$, where $\R^{pq}$ is in the center.
$\tilde R_3^{\triv} = \sll_r + \R^{rq} + \R^{rp} + \sll_p + \sll_q + \R^{pql_2}$, and finally
$$
\tilde R_i^{\triv} = \sll_r + \R^{rq} + \R^{rp} + \sll_p + \sll_q + \sll_{l_2}+\R^{pql_2},
$$
for all $i \leq 4$.

\end{enumerate}

Since there are no other possibilities for the integers $p,q$ and $r$,
we have shown that the sequence $\tilde R^{\triv}_i$ stabilizes after at most
$6$ steps, which completes the proof.
\end{proof}

\begin{reme}
Actually several of these cases are impossible and the sequence $\tilde R_i^{\triv}$
stabilizes in less then $6$ steps.
\end{reme}

\begin{reme}
We used that $l_i \geq 2$ only in cases 2, 5, 6, 7 and 11. If we
had that $l_i=1$ for all $i\geq 2$, we would have that the
sequence $R^{\triv}_i$ does not stabilize, but it will at least
weakly stabilize.

This together with $h(2)=4$ suggests that it may be possible to
construct relatively free Lie algebra on $2$ generators in some
nilpotent variety such which has not stabilizing tower of
restricted derivations.
\end{reme}

The next corollary is an immediate consequence of the theorem.

\begin{co}
The restricted derivation tower of $L[V,d]$
stabilizes after at most $\max\{3+ d/h(n),6\}$ steps.
\end{co}
\begin{co}
If Conjecture~\ref{gen:conj} holds, then the restricted
derivation tower of $L[V,d]$ stabilizes at the third level if
$d \geq  n+1$ and $n | d-1$.
\end{co}

\begin{reme}
\label{faithfulk:center}
As in the case when the center of $\D_0 L$ is trivial (
see Remark~\ref{faithfulk:nocenter}),
we have some description of
the reductive part of the algebras $DL^k$
The reductive part of $DL^k$ is a semi simple Lie algebra
which decomposes as a sum of copies of $\sll_l$ for some $l$-es.
It also acts faithfully on the nilpotent part and almost all nontrivial
modules which appear in are either isomorphic to the standard one or its dual.
Actually the same is the for the reductive part of $\tilde R^{\triv}_k$
but to show this one need to go over all cases in Lemma~\ref{nil:Wsequence}
and verify this statement.
\end{reme}

\end{section}

\end{chapter}
\newpage

\begin{chapter}{Automorphism Tower of Free Nilpotent Groups}
\thispagestyle{empty}
	\nocite{Weiss,Kass4,Kass5,Hochschild,Andre}

In this chapter we study the automorphism tower of the free
nilpotent group $\Gamma(n,d)$. The main tool in obtaining
the description of this tower is to embed the groups
$\Aut^k \Gamma(n,d)$ into a Lie group over (the reals)
$G^k$ and use methods form Lie algebra to study the
automorphism group of $G^k$. After that we use a
rigidity result proved in Grunewald and Platonov in~\cite{GrunPla},
to show that we can embed the next group in
automorphism tower into $\Aut G^k$.

Theorems~\ref{auttower:nocenter} and~\ref{grtower:center}
show that there is a close connection between the
groups in the automorphism tower of $\Gamma(n,d)$, for $n\geq 3$,
and the algebras in the tower of restricted derivations
of the free Lie algebra $L[V,d]$. This allows us to
derive the main result in this dissertation that the
automorphism tower of the free nilpotent group stabilizes
after finitely many steps (at least in the case of more
than $2$ generators).

The two generated case is slightly mote complicated, because
we can not embed the groups $\Aut^k \Gamma$ for $k\geq 2$ into
Lie groups.
The reason for this is that
the group $\SL_2(\Z)$ has nontrivial cocycles.
We go around this problem in section~\ref{sec:2gen} by finding a
characteristic subgroup $A\Gamma^k$ in $\Aut^k \Gamma$, which
can be embedded in a Lie group $G^k$ and first study the
automorphism group of this subgroup.

\begin{section}{Automorphism tower}

For any group $G$ we have a natural map from
$\ad: G \to \Aut G$, which sends every element $g$
to the inner derivation of the group $G$, given by
conjugation with the element $g$.

\begin{df}
A group $G$ is called compete if the map $\ad$ is an isomorphism,
i.e., if every automorphism of $G$ is inner and the center of $G$
is trivial.
\end{df}

Iterating the construction of the automorphism group
allows us to define the automorphism tower of the group $G$.

\begin{df}
Let $G$ be a group, by the automorphism tower of $G$ we mean the sequence
$$
\Aut^0 G \to \Aut^1 G \to \Aut^2 G \to
\dots \to \Aut^k G \to \cdots,
$$
where $\Aut^0 G = G$ and $\Aut^{k+1} G = \Aut (\Aut^{k} G)$.
The homomorphisms $\ad_k : \Aut^{k} G \to \Aut^{k+1} G$
comes form the adjoint action of $\Aut^k G$ on itself.
\end{df}
\begin{reme}
If the group $G$ has a trivial center then the map $\ad_0$
 is injection and it can easily be seen that all groups $\Aut^k G$
have trivial centers. Therefore $\Aut^k G$ form an increasing
sequence of groups, which explains why $\Aut^k$ is called a tower.
\end{reme}

\begin{df}
\label{auttower:stab}
We say that the automorphism tower of the group $G$ stabilizes at
level $k$ (or after $k$ steps) if the map $\ad_k$ is an
isomorphism. We also say that the tower weakly stabilizes at level
$k$, if $\Aut^k G \simeq \Aut^{k+1} G$.
\end{df}

A classical result of Wielandt~\cite{Wielandt}, asserts that the
automorphism tower stabilizes for any finite group $\Gamma$ with
trivial center. Similar result for infinite groups with center
does not hold in general.

Baumslag conjectured that for any finitely generated
nilpotent group $\Gamma$, the automorphism tower of $\Gamma$
stabilizes after finitely many
steps. (Since this conjecture is not published, we will refer to the strongest
version of the conjecture, i.e., that the tower
stabilizes in the sense of Definition~\ref{auttower:stab}.)

In~\cite{DyerForm1}, Formanek and Dyer proved that the automorphism tower of
the free nilpotent group of class $2$ stabilizes at the second level,
when the number of generators is different than
$3$. For the three generator case, the tower stabilizes at the third level.
Similar results are known for other nilpotent groups, like the group of
$n\times n$ upper triangular matrices, and in all these case the
automorphism tower is very short.

In the rest of this chapter we study the automorphism tower of the free
nilpotent groups $\Gamma(n,d)$.
\end{section}

\ifnewpage

\begin{section}{Rigidity results}

Our approach in studying the automorphism tower of the free
nilpotent group $\Gamma(n,d)$ is to embed each group
$\Aut^k \Gamma$ as a lattice into a Lie group over reals $G^k$.

We say  that the pair $(G,\Gamma)$ of Lie group and a lattice in it,
is rigid if any automorphism of the lattice $\Gamma$ comes from
an automorphism of the Lie group $G$. (In section~\ref{sec:2gen} we will
see that not all groups in the automorphism tower of $\Gamma(n,d)$
can be embedded as a rigid lattices in a Lie group.)

The first rigidity results were proved my Margulis~\cite{Margulis2}
and Mostow~\cite{Mostow} for lattices in higher rank simple
Lie groups. These results were latter generalized to some
simple rank $1$ groups.

All groups in the automorphism tower of the free
nilpotent group $\Gamma(n,d)$ have big nilpotent
part and they are very far from semi simple group,
for which the classical rigidity results hold.

In a resent paper, F. Grunewald and  Platonov~\cite{GrunPla}
proved a rigidity result for a groups with large
nilpotent radical. Here we will state one of the main
results from their paper, which we will use in the
next sections. This rigidity theorem
is valid only for groups with strong unipotent radical.

\begin{df}
Let $G$ be a Lie group. We can write $G=U.H$, where $U$ is
the unipotent radical of $G$ and $H$ is a reductive group,
which has trivial intersection with $U$.
We say that $G$ has a strong unipotent radical
if the adjoint action of $H$ on $U$ is faithful.
\end{df}

For groups with strong unipotent radical we
have the following rigidity result.

\begin{thq}[Grunewald, Platonov]
\label{GP}
Let $G$ be an algebraic  group with strong unipotent
radical and let $\Gamma$ be a lattice in $G$. Then for any
automorphism $\Phi$ of the lattice $\Gamma$ there exists an
automorphism $\phi$ of $G$ and a map $v :\Gamma \to \Cen(U)$ such
that
$$
\Phi(\gamma) = v(\gamma) \phi(\gamma),
$$
for any $\gamma \in \Gamma$. Also we have that
$v$ is trivial on $\Gamma \cap U$.
\end{thq}
\begin{proof}
Here we give a heuristic explanation why such
result holds. For a complete proof we refer the readers
to the original paper (see~\cite{GrunPla}).

Let denote $\IGamma = \Gamma \cap U$. This is a nilpotent group. We
can reconstruct $U$ starting from $\IGamma$ first by taking the
Malcev completion of the nilpotent group $I\Gamma$, and then by
completing the resulting arithmetic group over $\Q$ to a Lie group
over $\R$. This allows us to lift any automorphism $\Phi$ of the
group $\IGamma$ to an automorphism of $\tilde \phi$ of $U$.

If group $G$ have strong unipotent radical then $G/\Cen(U)$ acts
faithfully on $U$ therefore we can embed $G/\Cen(U)$ into $\Aut
U$. The image of this embedding coincides with the Zariski
closure of $\Gamma/\Cen(\IGamma)$ in $\Aut U$, where the last
embedding comes from the adjoint action of $\Gamma$ on $\IGamma$
and lifting the automorphisms for $\IGamma$ to automorphisms of
$U$.

Let $\Phi$ is an automorphism of the group $\Gamma$, we can
restrict $\Phi$ to an automorphism $\bar \Phi$ of $\IGamma$ and
then lift $\bar \Phi$ to an automorphism $\bar \phi$ of $U$. This
automorphism is an element of $\Aut U$ and it acts on $\Aut U$ by
conjugation as an automorphism $\tilde \phi$.

By construction $\tilde \phi$ preserves the image of
$\Gamma/\Cen(\IGamma)$ in $\Aut U$, therefore it preserves
$G/\Cen(U)$. Thus we have obtained an automorphism $\tilde \phi$
of the group $G/\Cen(G)$ which we can restrict to the reductive
part of the group $G$. This together with $\bar \phi$,
allows us to construct an
automorphism $\phi$ of the group $G$. By construction we have that
for any $\gamma \in \IGamma$ we have that $\phi(\gamma) =
\Phi(\gamma)$. Also this construction gives that the adjoint
actions on $I\Gamma$ of $\phi(\gamma)$ and $\Phi(\gamma)$ on
$I\Gamma$ coincide for any $\gamma \in \Gamma$ . This shows that
$$
\Phi(\gamma) = v(\gamma) \phi(\gamma),
$$
for some map $v: \Gamma \to \Cen U$, which is trivial on $I\Gamma$.
\end{proof}
\begin{reme}
$\phi$ and $v$ are not uniquely defined.
If $\phi$ acts trivially on $G/U$ then
the map $v$ is a cocycle and it is defined up to a
cohomotopy. Therefore we say that if $G$ has strong unipotent
radical every automorphism of the lattice $\Gamma$ comes from
an automorphism of the Lie group up to a cocycle from
the reductive part $\Gamma/(\Gamma \cap U)$ of $\Gamma$ to the
center of the group $U$. This gives that
$\Aut \Gamma / (\Aut \Gamma \cap \Aut G)$
can be embedded into $H^1(\Gamma/(\Gamma \cap U, \Cen(U))$.
\end{reme}

If the reductive part of $G$ is a semi simple group, which is a
product of higher rank simple Lie groups, then we have that all
cohomology groups are trivial, which gives that any automorphism
of $\Gamma$ comes from an automorphism of the Lie group $G$.
Similar results hold if the semi simple part have $\SL_2(\Z)$ as a
quotient and $\Cen(U)$ decomposes as $\SL_2(\Z)$ module as a sum of
trivial modules and standard ones (see~\cite{Borel}).

\end{section}
\ifnewpage

\begin{section}{Defining analogs of $UL_k$ at group level}

In the previous chapter we showed that the ideals $UL_k$ play in
important role in the description of the tower of restricted
derivations of the free nilpotent Lie algebra. In this section we
will define analogs of this ideals inside the group
$\Aut \Gamma$.

Let us observe that the semi simple part of the group $\Aut \Gamma$ is
$\GL_n(\Z)$, therefore we need to construct the analogs of the
ideals $UL_k$ using the action of the group $\GL_n(\Z)$, not the
action of the Lie algebra $\sll_n$.

\begin{df}
Let $\widetilde{UL}^+$ (respectively $\widetilde{UL}^{\triv}$) denotes the maximal
submodule of $UL$ without $\GL_n(\Z)$ invariant vectors
(consisting entirely of invariant vectors). Let us construct the ideals
$\widetilde{UL}_k$ of $UL$ in the same way we constructed the ideals $UL_k$
in definition~\ref{dist:ideals}, but starting form $\widetilde{UL}^+$
instead $UL^+$.
\end{df}

Now can define the group analogs of the ideals $UL_k$ which will help us
to construct the automorphism tower of the group $\Gamma(n,d)$.

\begin{df}
Let $UG_k(n,d)$ be the Lie subgroup of $UG(n,d)$ obtained by
exponentiating the ideal $\widetilde{UL}_k$.
\end{df}

\begin{lm}
For any $k$, $UG_k(n,d)$ is normal subgroup of $G^1(n,d)$,
which is preserved by all automorphisms of this group.
\end{lm}
\begin{proof}
The Lie algebras corresponding to $UG_k$ by construction are
ideals in the Lie algebra $\D_0 L(n,d)$ corresponding to $G^1$.
Therefore $UG_k$ are normal subgroups. These subgroups are
invariant under all automorphisms of $UG_k$, because the ideals
$UL_k$ are characteristic ideals in $UL$ and $\D_0 L$ by
Lemma~\ref{seq}.
\end{proof}

\begin{lm}
$UG_k$ form a decreasing sequence of normal subgroups of $UG$ (and
$G^1$), which stabilizes after at most $1+d/h(n)$ steps.
\end{lm}
\begin{proof}
The stabilization of this sequence of subgroups follows
immediately from the stabilization their Lie algebras $UL_k$ which
is given by Theorem~\ref{ULstabilization}.
\end{proof}

\begin{lm}
The factor group
$UG_0/UG_\infty$ is abelian.
\end{lm}
\begin{proof}
The Lie algebra corresponding to this factor group is
$UL_0/UL_{\infty}$ which is abelian by construction.
\end{proof}

\begin{df}
Let us define the analogs of the spaces $U_{\gen,k}$ using the Lie
group $UG$. Let us denote $UG_{\gen,k}(n,d) =
UG_k(n,d)/UG_{k+1}(n,d)$ and the corresponding subgroups in the
center of the $UG$:
$$
CUG_k(n,d) = UG_k(n,d) \cap \Cen(UG(n,d)) \quad \mbox{and}
$$
$$
CUG_{\gen,k}(n,d) = CUG_k(n,d)/CUG_{k+1}(n,d).
$$
\end{df}
\begin{reme}
All groups $UG_{\gen,k}$, $CUG_{k}$ and  $CUG_{\gen,k}$ are
abelian and are naturally isomorphic to their Lie algebras.
\end{reme}

Finally we need to define the discrete analogs of all these groups using the
group $\IAut \Gamma(n,d)$.

\begin{df}
Let us define the following groups using $\IAut \Gamma(n,d)$:
\begin{itemize}
\item
$\IAut \Gamma_k = \IAut \Gamma \cap  UG_k$;
\item
$\CIAut \Gamma_k = \IAut \Gamma \cap  CUG_k = \IAut \Gamma_k \cap \Cen(\IAut \Gamma)$;
\item
$\IAut \Gamma_{\gen,k} = \IAut \Gamma_k / \IAut \Gamma_{k+1}$;
\item
$\CIAut \Gamma_{\gen,k} = \CIAut \Gamma_k / \CIAut \Gamma_{k+1}$.
\end{itemize}
\end{df}

\begin{reme}
The groups $\IAut \Gamma_k$ and $\CIAut \Gamma_k$ are normal
subgroups of $\Aut \Gamma$ and $\IAut \Gamma$, which are invariant
under all automorphisms. They also form two decreasing sequences
of normal subgroups of length at most $1+ d/h(n)$.
\end{reme}

\begin{lm}
The groups $C\IAut \Gamma_k$, $\IAut \Gamma_{\gen,k}$ and $C\IAut \Gamma_{\gen,k}$
are free abelian and there is a natural action of $\GL_n(\Z)$ on them.
\end{lm}
\begin{proof}
The groups are free abelian groups because they are lattices in
the corresponding Lie group which are abelian (without compact
factors).

The group $\Aut \Gamma$ acts by conjugation on $\IAut \Gamma_k$.
Its subgroup $\IAut \Gamma$ acts trivially modulo $\IAut
\Gamma_{k+1}$ (because $[UL,UL_k] \subset UL_{k+1}$ by
construction), therefore we can define action of the factor group
$\Aut \Gamma/\IAut \Gamma$ on $\IAut \Gamma_k/\IAut \Gamma_{k+1}$.
This defines the action of $\GL(\Z)$ on $\IAut \Gamma_{\gen,k}$,
the action on the other groups are defined in similar way.
\end{proof}

\end{section}
\ifnewpage

\begin{section}{Automorphism tower case with trivial center}

In this section we will describe the automorphism tower of the
free nilpotent group $\Gamma(n,d)$ in the case $n\geq 2$ and $\Cen
\Aut \Gamma =0$ (by a result of Formanek~\cite{Formanek} this is
equivalent to $2n\not | d-1$).
In this and in the next two section we will write $\Gamma$ on 
instead of $\Gamma(n,d)$ to prevent to notations from becomming too
complicated.


Let us start by first describing the second group in the automorphism tower.

The group $\Aut \Gamma$ is a lattice in the Lie group $G^1$. Since
the semi simple part of $\Aut \Gamma$ is $\GL_n(\Z)$ which do not
have any non trivial cocycles if $n\geq 3$ (the case $n=2$ is
slightly more complicated --- see the proof of Theorem~\ref{outd:nocenter})
by
the rigidity results we can embed $\Aut (\Aut \Gamma)$ into $\Aut
G^1 $. Therefore out first goal is to describe the automorphism
group $\Aut G^1$.

\begin{thq}
\label{outlie:nocenter}
The group of outer automorphisms of $G^1(n,d)$ is isomorphic to
$$
(\R^*/\{\pm 1\}) \ltimes WG_0,
$$
where
$$
WG_0 \subset \Hom_{\GL(\Z)}(UG_{\gen,0}, CUG_0)/
\big( UG^{(d-2)}/\Cen(UG)\big)^{\GL(\Z)}
$$
is the subset consisting of all map whose projection in
$\End(CUG_{\gen,0})$ is invertible, here the projection from
$\Hom(UG_{\gen,0}, CUG_0)$ to $\End(CUG_{\gen,0})$ is defined
using the natural embedding of $CUG_{\gen,0}$ in $UG_{\gen,0}$ and
the natural projection from $CUG_{0}$ to $CUG_{\gen,0}$.
The group structure on the set $WG_0$ comes from the composition of the
maps through $CUG_{\gen,0}$.

The group $\big( UG^{(d-2)}/\Cen(UG)\big)^{\GL(\Z)}$
is the group of all inner derivations which lie in \\
$\Hom(UG_{\gen,0}, CUG_0)$. Finally $\R^*/\{\pm 1\}$ comes from $\GL(\R) / \SL^{\pm 1} (\R)$

\end{thq}
\begin{proof}
The Lie algebra of the group $\Aut G^1$ is a subalgebra of the derivation algebra of
$G^1$. Our first goal is to describe this subalgebra.

By Theorem~\ref{2nd:trivcenter} the space of outer derivations of the Lie algebra of the group
$G^1$ is isomorphic to
$$
Out = \R + \Hom_\sll (UL_{\gen,0}, CUL_0)/(UL^{(d-2)})^{\sll},
$$
because if $n\geq 3$ the spaces $T$ is trivial. But the group
$G^1$ has two connected components, therefore we need to take the
subalgebra which preserved by the action of the outer automorphism,
which comes from the conjugation with an element in the connected
component of $G^1$ which does not contain the identity.
The subalgebra of $Out$ preserved by this automorphism
consists of all maps which preserve the action of $\GL(\Z)$.

Now we only need to exponentiate this Lie algebra in order to
construct the Lie group. Exponentiating $\R$ gives us the
extension of the semi simple part form $\SL_n^{\pm 1}(\R)$ to
$\GL_n(\R)$, which corresponds to the factor $\R^*/\{\pm 1\}$ in the
outer automorphism group. Exponentiation of the other part
$$
\Hom_{\GL(\Z)} (\widetilde{UL}_{\gen,0}, \widetilde{CUL}_0)/(UL^{(d-2)})^{\GL(\Z)}
$$
gives the group $WG_0$.
\end{proof}

We can use the description of $\Aut G^1$, to obtain a description
of $\Aut (\Aut \Gamma))$.

\begin{thq}
\label{outd:nocenter}
The group of outer automorphisms of $\Aut \Gamma$ can be embedded into
$$
\Hom_{\GL(\Z)}(\IAut \Gamma_{\gen,0} , \CIAut \Gamma_0) /
\big(\IAut \Gamma^{(d-2)}/\Cen(\IAut \Gamma)\big)^{\GL(\Z)}
$$
The image of this embedding consists of all maps, whose projection
$\phi$ in  the space $\End_{\GL(\Z)}(\CIAut_{\gen,0})$ is invertible and it is
possible to lift $\phi$ to an automorphism $\tilde \phi$ of $G^1$ which
preserves the group $\Aut \Gamma$.
\end{thq}
\begin{proof}
First we want to prove that $\Aut (\Aut \Gamma)$ can be embedded
into the Lie group $\Aut G_1$. In order to prove that we need to
verify that the pair $(G^1,\Aut \Gamma)$ satisfies the condition of
Theorem~\ref{GP}.

The group $G_1$ has a strong unipotent radical, because
$SL^{\pm1}_n(\R)$ acts faithfully on $UG$ (since the action of
$\sll$ on $UL$ is not trivial). Therefore we can lift any
automorphism of $\Aut \Gamma$ to an automorphism of $G^1$ up to a
cocycle $\gamma$ from $\GL_n(\Z)$ to $\Cen(UG)$.

If we have $n\geq 3$ then every cocycle of $\GL_n(\Z)$ is
homologically trivial.
In the case $n=2$, we have that $\Cen(UL) =
\sum_k V^{[2k+1]}$ as $\SL_2(\Z)$ module because it lies in the
homogeneous component of odd degree. But the group $\GL_2(\Z)$
does not have homologically non trivial cocycles in the module $U$
unless $U$ contains submodules corresponding to a partition $\lambda=[2k]$
for $k \geq 2$ (see~\cite{Borel}). This
implies that $\gamma$ is homologically trivial, therefore we
may assume that $\gamma$ is trivial and
that $\Aut(\Aut \Gamma)$ is a subgroup
of $\Aut (G^1(n,d))$. In order to describe the image it is enough
to find all outer automorphisms of $G^1$ which preserve $\Aut
\Gamma$.

Let $\phi \in WG_0$ be an outer automorphism, in order to preserve
$\Aut \Gamma$ its image has to be in $\Cen(\IAut \Gamma)$
therefore $\phi$ is an element in
$$
\Hom_{\GL(\Z)}(\IAut \Gamma_{\gen,0} , \CIAut \Gamma_0) /
\big(\IAut \Gamma^{(d-2)}/\Cen(\IAut \Gamma)\big)^{\GL(\Z)}.
$$
For an element $\phi$ to preserve the lattice $\Aut \Gamma$ it is necessary
that its projection in $\End_{\GL(\Z)} (\CIAut_{\gen,0})$ is
invertible. A sufficient condition for $\phi$ to preserve this
lattice is that $\phi$ lies in some unipotent subgroup of $WG_0$.
Therefore using that $\SL(\Z)$ is generated by the elementary
matrices  we can see that if for any partition $\lambda$ we have
$\det \phi \End_{\GL(\Z)} (\CIAut^{\lambda}_{\gen,0}) =1$. Then
$\phi$ can be lifted to an isomorphism of $\Aut \Gamma$.
\end{proof}

Let $G^2$ be the Lie group corresponding to $\Aut^2 \Gamma$, i.e., we define
$G^2$ as the Zariski closure of $\Aut^2 \Gamma$ in the Lie group $\Aut G^1$.
The next lemma describes the Lie algebra of the group $G^2$. From this result
we can notice that the Lie algebra of $G^2$ is almost the same as the
second algebra in the restricted derivation tower of $L[V,d]$
which justifies the study of this tower.

\begin{reme}
\label{Zcl:nocenter}
The Lie algebra of the group $G^2$ is
$$
\sll + UL + \Hom^0_{\GL(\Z)}(\widetilde{UL}_{\gen,0},\widetilde{CUL}_0).
$$
\end{reme}
\begin{proof}
The Lie algebra contains the subalgebra of $W_0$ generated by all
nilpotent derivations, therefore it contains
$R_0^\lambda \cap W_0$.
This algebra is also a sub algebra of
$\D_0 {\mbox{Lie}}G_1$,
which implies that the Lie algebra of $G^2$ is exactly
$$
\sll + UL + \Hom^0_{\GL(\Z)}(\widetilde{UL}_{\gen,0},\widetilde{CUL}_0).
$$
By Remark~\ref{faithful2} the semi simple part of this Lie algebra
is a sum of copies of $\sll_l$ for some $l$-es and it acts
faithfully on the nilpotent part. Also almost all nontrivial modules which
appears in the nilpotent part are either isomorphic to the standard modules
or their duals.
\end{proof}

\begin{reme}
We do not have a good description of the connected components of
$G^2$, in order to obtain such description we need a detailed knowledge of how the
group $\Aut \Gamma$ sits inside $\Aut G$. For example it can be
shown that $G^2$ has two connected components if $d\leq n(n-1)$
unless $n=3, d=2$, when it has $4$ connected components (for $d=2$
this follows from a result of Formanek and Dyer~\cite{DyerForm2}, about
the description of the automorphism tower of free nilpotent groups
of class $2$).
\end{reme}

In section~\ref{sec:dertower} we used the description of the second derivation algebra
of $L[V,d]$ as a base of induction to describe all algebras in the
tower of restricted derivation of $L[V,d]$. We can do the same for the case
of discrete group $\Gamma$ an its automorphism tower.

We can use the description of the second automorphism group of
$\Gamma(n,d)$ as a base for induction which give us the structure
of every group in the automorphism tower of the group $\Gamma$.

\begin{thq}
\label{auttower:nocenter}
a) The group of outer automorphisms of $G^k$ is isomorphic to
$$
\big(\R^*/\{\pm 1\} \big) \ltimes WG_{k-1},
$$
where
$$
WG_{k-1} \subset \Hom_{\GL(\Z)}(UG_{\gen,k-1}, CUG_{k-1})
$$
is the subset consisting of all map whose projection in
$\End(CUG_{\gen,k-1})$ is invertible, where the multiplication on this
space is given by a composition trough $CUG_{\gen,k-1}$.

b) The group of outer automorphisms of $\Aut^k \Gamma$ can be
embedded in
$$
\Hom_{\GL(\Z)}(\IAut \Gamma_{\gen,k-1} , \CIAut \Gamma_{k-1}).
$$
The image of this embedding consists of all maps $\phi$, whose projection
of  in the space $\End_{\GL(\Z)}(C\IAut_{\gen,k-1})$ is invertible and the
lifting of $\phi$ to an automorphism $\tilde \phi$ of
$\SL^{\pm 1}(\R) \ltimes UG_0$
preserves the projection of $\Aut \Gamma$ into this Lie group.

c) The automorphism group of $\Aut^k \Gamma$ can be embedded into
$\Aut G^k$. The Lie algebra of the Zariski closure $G^{k+1}$ of
$\Aut^k \Gamma$ is isomorphic to
$$
\sll + UL + \sum_{i=0}^{k-1} \Hom^0_{\GL(\Z)}(\widetilde{UL}_{\gen,i},\widetilde{CUL}_i).
$$
\end{thq}
\begin{proof}
a) The proof is by induction. Let us first describe the
automorphism group of $G^k$. Using the lemma from the proof of
Theorem~\ref{dertower:nocenter} we can see that the algebra of
outer derivations which are preserved by the action of the outer
automorphism of order $2$ obtained by conjugation by elements in
the connected component of $G^k$ corresponding to $\SL_n^{-1}(\R)$
is
$$
\R + \Hom_{\GL(\Z)}(\widetilde{UL}_{\gen,k-1}, \widetilde{CUL}_{k-1}).
$$
From Lemma~\ref{nil:posdegree} it follows that this algebra is preserved by
all other outer automorphisms of the Lie algebra of $G_k$ coming
from the other connected components. Therefore it coincides with
the Lie algebra of the group of outer automorphisms of $G_k$.

We can obtain the group $WG_k$ by exponentiating this Lie algebra
as we did in Theorem~\ref{outlie:nocenter} for the group $\Aut G^1$.

b) In order to prove that $\Aut^{k+1} \Gamma$ can be embedded in
$\Aut G_k$ we need to verify that the pair $(G^k,\Aut^k \Gamma)$
satisfies the condition in Theorem~\ref{GP}. The group $G^k$ has
strong unipotent radical, because the semi simple part of its Lie
algebra acts faithfully on the nilpotent part by
Remark~\ref{faithfulk:nocenter}.

Now we need to verify that there are no homologically nontrivial
cocycles from the semi simple part of the group $\Aut^k \Gamma$ to
the center $V$ of the unipotent part of the group $G^k$. All
simple factors $S$ of $\Aut^k \Gamma$ are isomorphic to
$\SL_l(\Z)$ or $\GL_l(\Z)$ for different $l$-es. Unless $S$ is original copy of
$\SL_n(\Z)$ then under the action of $S$, the space  $V$
decomposes as a sum of modules which are either trivial or
isomorphic to the standard module or its dual, therefore every
cocycle from $S$ to $V$ is homologically trivial. If $S$ is the
original copy of $\GL_n(\Z)$ then as we saw in the proof of
Theorem~\ref{outd:nocenter} there are no nontrivial cocycles.

This shows that $\Aut (\Aut \Gamma) \subset \Aut G^k$, and in
order to describe the group of outer automorphisms of
$\Aut^k \Gamma$ we only need to see which elements in $WG_{k-1}$
preserve the lattice $\Aut^k \Gamma$. The argument is the same as
the argument in the case for the second automorphism group.

c) The Lie algebra of $G^{k+1}$ contains the subalgebra of
$W_{k-1}$ generated by all nilpotent derivations, therefore it
contains $R_k^\lambda \cap W_k$. This algebra is also a sub
algebra of $\D_0 {\mbox{Lie }}G^k$, which implies that the Lie
algebra of $G^{k+1}$ is exactly
$$
\sll + UL + \sum_{i=0}^{k-1} \Hom^0_{\GL(\Z)}(\widetilde{UL}_{\gen,i},\widetilde{CUL}_i).
$$
Notice that by Remark~\ref{faithfulk:nocenter} the semi simple part of
this Lie algebra acts faithfully on the nilpotent part.
\end{proof}

An immediate corollary of this theorem is the following result
about the stabilization of the automorphism tower of the free nilpotent groups.

\begin{thq}
\label{auttower:st}
If $2n \not | d-1 $ then the automorphism tower of the free nilpotent
group $\Gamma(n,d)$ stabilizes after at most $2 + d /h(n)$
steps, in particular it has finite height.
\end{thq}

\begin{reme}
The height of the automorphism tower of $\Gamma$ is usually
the same as the length of the sequence $\widetilde{UL}_k$
(more precisely the maximal $k$ such that the space
$\Hom_{\GL(\Z)}(\widetilde{UL}_{\gen,k}, \widetilde{CUL}_k\big)$
is not trivial).
The height of the
automorphism tower may be smaller then the length of this sequence
if and only if for the maximal  $k$ such that $\widetilde{UL}_{\gen,k}$ is not
trivial we have that
$$
\widetilde{UL}_{\gen,k}^\lambda = \widetilde{CUL}_{\gen,k}^\lambda
= \widetilde{CUL}_{k}^\lambda
$$
and this is a simple $\sll$ module corresponding to the partition $\lambda$,
 for all $\lambda$ such that
$$
\Hom_{\GL(\Z)}(\widetilde{UL}_{\gen,k}^\lambda, \widetilde{CUL}_{k}^\lambda) \not = 0.
$$
In this case it is possible that the height of automorphism tower
is the length of the sequence $UL_k$ minus $1$.

In particular if $d \leq 2n$ the automorphism tower stabilizes at
the second or the third
level. If the tower does not stabile at the second level then the
group of outer automorphisms $\Aut^2 \Gamma /\Aut^1 \Gamma$ is
cyclic group of order $2$. It can be shown that unless $d=2$ and $n=2$
the tower stabilizes at the second level.
Similarly
if $2n < d \leq h(n)$ and $2n\not | d-1$, then the automorphism tower stabilizes
at the third level.
\end{reme}

\begin{reme}
This result is generalization of the result by Formanek and Dyer for
$d=2$, however in the case $d=2$ their result is more precise
because it says that the tower stabilizes after $2$ steps unless
$n=2$ and our result gives that it stabilizes after at most $3$
steps and that the quotient $\Aut^2 \Gamma /\Aut \Gamma$ is
finite.
\end{reme}
\end{section}

\ifnewpage

\begin{section}{Automorphism tower case with center}

In the previous section we obtained a description of
the automorphism tower of the group $\Gamma$ in the case when
$\Cen(\Aut \Gamma)$ is trivial. Now we want to consider the case
when the center of the first group in the tower is not trivial.

There is a close analogy between this case and the tower of
restricted derivations of $L[V,d]$ if the first algebra
has a center, for example all groups in the automorphism tower of
$\Gamma$ splits as a direct product of two groups.

In this section
we will assume that $n\geq 3$ which allows us to embed every group
in the automorphism tower in some Lie group. The case of $2$
generated groups will be considered in the next section.

Let us start with the second group in the tower. As in the case without center
we first describe the automorphism group of $G^1$ and show that
$\Aut^2 \Gamma$ can be embedded into this group.

\begin{thq}
\label{2ndlevel:center}
a)
The automorphism group of $G^1$ splits as
$$
\Aut G^1 = (\R^*/\{\pm 1\})\ltimes (\widetilde{NG}^2 \times \widetilde{CG}^2),
$$
where the group
$\widetilde{NG}^2$ is the
extension of the image of the group of inner automorphisms by
$$
WG_0 \subset \Hom_{\GL_n(\Z)}(UG_{\gen,0}, CUG_0),
$$
where $WG_0$ is the subset consisting of all maps whose projection in
$\End(CUG_{\gen,0})$ is invertible and
$$
CG^2 \subset \Hom_{\GL_n(\Z)}(UG_{\gen}^{\triv}, \Cen (G^1)).
$$

b)
The automorphism group of $\Aut \Gamma$ splits as a direct
product $\Aut^2 \Gamma =  N\Gamma^2 \times C\Gamma^2$, where
$N\Gamma^2$ and $C\Gamma^2$ are discrete subgroups of
$\widetilde{NG}^2$ and $\widetilde{CG}^2$. The group $N\Gamma^2$ is an extension of
$\Aut \Gamma / \Cen(\Aut \Gamma)$ by the subgroup of
$$
\Hom_{\GL(\Z)}(\IAut \Gamma_{\gen,0} , \CIAut \Gamma_0)
$$
consisting of all maps $\phi$, whose projection in
$\End_{\GL(\Z)}(\CIAut_{\gen,0})$ can be lifted to
automorphism of $\tilde \phi$ which preserves $\Aut \Gamma$. The group
$C\Gamma^2$ is
$$
C\Gamma^2 \subset \Hom_{\GL_n(\Z)}(\IAut\Gamma_{\gen}^{\triv}, \Cen (\Aut \Gamma))
$$

c)
The Zariski closure $G^2$ of the  group $\Aut^2 \Gamma$ in
$\Aut G^1$ splits as a direct product $G^2 = NG^2 \times CG^2$, and the
Lie algebra of $G^2$ is isomorphic to
$$
\sll + UL/\Cen(UL)^{\GL(\Z)} + \Hom^0_{\GL(\Z)}(\widetilde{UL}_{\gen,0},\widetilde{CUL}_0)
\oplus \Hom^0_{\GL(\Z)}(\widetilde{UL}_{\gen}^{\triv},\Cen(UL)^{\GL(\Z)}).
$$
\end{thq}
\begin{proof}
a) The situation is similar to  Theorems~\ref{outlie:nocenter}.
The automorphism group of $G^1$ splits because the its Lie algebra is
a direct sum of two components (see Theorem~\ref{2nd:withcenter}).

b) The automorphism group of $\Aut \Gamma$ can be embedded into $\Aut G^1$,
because the group $G^1$ has strong unipotent radical and $\GL_n$ goes not have
any homologically non trivial  cocycles (see Theorem~\ref{outd:nocenter}).
The group splits as a direct product because automorphisms which acts trivially on
$\Cen(\Aut \Gamma)$ commutes with automorphism which acts trivially on
$\Aut \Gamma /\Cen(\Aut \Gamma)$.

c) The group $G^2$ splits as a direct product because the group $\Aut^2 \Gamma$
does. We can obtain the description of its Lie algebra from
Theorems~\ref{outlie:nocenter} and~\ref{2nd:withcenter}.

\end{proof}

The group $\Aut^2 \Gamma$ is a very good model for all group in the automorphism tower.
We use the previous theorem as a base for induction which will describe all
groups in the automorphism tower of $\Gamma(n,d)$.
\begin{thq}
\label{grtower:center}
a)
The automorphism group of $G^k$ splits as a direct product
$$
\Aut G^k = \big(\R^*/\{\pm 1\}\big)
\ltimes (\widetilde{NG}^{k+1} \times \widetilde{CG}^{k+1}),
$$
where the group
$\R^*/\{\pm 1\}\ltimes \widetilde{NG}^{k+1}$ is the automorphism group of  $NG^k$,
also $\widetilde{NG}^{k+1}$
is an extension of $NG^k$ by
$$
WG_{k-1} = \Hom_{\GL_n(\Z)}(UG_{\gen,k-1}, CUG_{k-1}),
$$
where $WG_{k-1}$ is the subset consisting of all maps whose projection in
$\End(CUG_{\gen,k-1})$ is invertible and
$$
\widetilde{CG}^{k+1} = \Aut CG^k \ltimes \Hom \big(NG^k/[NG^k,NG^k], \Cen(CG^k)\big).
$$

b)
The automorphism group of $\Aut^k \Gamma$ splits as a direct
product $\Aut^{k+1} \Gamma =  N\Gamma^{k+1} \times C\Gamma^{k+1}$, where
$N\Gamma^{k+1}$ and $C\Gamma^{k+1}$ are discrete subgroups of
$\widetilde{NG}^{k+1}$ and $\widetilde{CG}^{k+1}$.
The group $N\Gamma^{k+1}$ is an extension of
$N\Gamma^k$ by the subgroup of
$$
\Hom_{\GL(\Z)}(\IAut \Gamma_{\gen,k-1} , \CIAut \Gamma_{k-1})
$$
consisting of all maps $\phi$, whose projection in
$\End_{\GL(\Z)}(\CIAut_{\gen,0})$ can be lifted to
automorphism of $\tilde \phi$ which preserves $\Aut \Gamma$.
The group
$C\Gamma^{k+1}$ is isomorphic to
$$
\Aut C\Gamma^k \ltimes \Hom(N\Gamma^k/[N\Gamma^k,N\Gamma^k], \Cen(C\Gamma^k))
$$

c)
The Zariski closure $G^k$ of the  group $\Aut^{k+1} \Gamma$ in
$\Aut G^k$ splits as a direct product $G^{k+1} = NG^{k+1} \times CG^{k+1}$, and the
Lie algebra of $NG^{k+1}$ is isomorphic to
$$
\sll + UL/\Cen(UL)^{\GL(\Z)} +
\sum_{i=0}^{k-1} \Hom^0_{\GL(\Z)}(\widetilde{UL}_{\gen,i},\widetilde{CUL}_i).
$$
\end{thq}

\begin{proof}
The proof of this theorem is by induction on $k$, where Theorem~\ref{2ndlevel:center}
serves as a base case. The induction step is the same as in
Theorem~\ref{auttower:nocenter} the case when the center of $\Aut \Gamma$ is
trivial. The only thing that have to be added is the splitting of the group
$Auk^k \Gamma$ which follows from the splitting of the algebras in the
tower of restricted derivations.
\end{proof}

As in the case when the center is trivial we have that the
sequence $N\Gamma^k$ stabilizes after finitely many steps, because
$\IAut \Gamma_{\gen,k}$ are trivial for $k \geq 1 + d/h(n)$.
However this is not enough to claim that the automorphism tower
of $\Gamma(n,d)$ stabilizes.

The stabilization of the split part $C\Gamma^k$ does not follow
the stabilization of $\tilde R^{\triv}_k$, since we do not have that
$\tilde R^{\triv}_k$ is the Lie algebra of $CG^k$.
Before starting the proof of the stabilization of $C\Gamma^k$ we need
a technical Lemma similar to Lemma~\ref{bigdim}

\begin{lm}
\label{bigabelian} The abelinization of the group $N\Gamma^k$ is
$$
N\Gamma^k/[N\Gamma^k,N\Gamma^k] = \Z^{l_k} \times \F_2^{m_k},
$$
where $l_k \geq 2$ and $m_k \geq 1$, unless $C\Gamma^2$ is trivial.
The cyclic factors comes from the copies of $\GL_l(\Z)$ which
appear in the semi simple part of the group $\Aut^k \Gamma$.
\end{lm}

\begin{proof}
Form the construction of the group $NG^k$ it follows that their
abelinization is a sum of several copies of $\Z$ and cyclic groups of
order $2$, therefore we have that the above isomorphism for some
numbers $l_k$ and $m_k$. Lemma~\ref{bigdim} gives that $l_k \geq
2$. $m_k\geq 1$ because $N\Gamma^k$ has $\GL_n(\Z)$ as a quotient
therefore its abelinization contains at least one element of order
$2$.
\end{proof}

Now we are ready to proof the stabilization of the sequence $C\Gamma^k$.
\begin{lm}
\label{CG:stab}
The sequence
$$
C\Gamma^2 \to C\Gamma^3 \to C\Gamma^4 \to cdots
$$

stabilizes after at most $5$ steps.
\end{lm}
\begin{proof}
The behavior of this sequence is determined by
the dimensions of the ranks of the
free abelian groups
$\IAut \Gamma_{\gen}^{\triv}$, $\Cen(\Aut \Gamma)$ and $\Cen(\Aut \Gamma) \cap [\IAut \Gamma,\IAut \Gamma]$.
Let us define $p,q,r$ as follows
$$
p = \rank (\IAut \Gamma / (\Cen(\IAut \Gamma) [\IAut \Gamma,\Aut \Gamma]), \quad
q = \rank \Cen(\Aut \Gamma  \cap [\IAut \Gamma,\IAut \Gamma])
$$
and
$$
r = \rank (\Cen (\Aut \Gamma) / \Cen(\Aut \Gamma) \cap [\IAut \Gamma,\IAut \Gamma]).
$$
We will use $l_i$ (and $m_i$) to denote the rank (and the dimension of
the torsion part of the group $N\Gamma^k/[N\Gamma^k,N\Gamma^k]$. Denote $s_i=l_i + p_i$.

Let us consider the different possibilities for this dimensions:
\begin{enumerate}
\item
$p=r=0$ or $r=q=0$.
In this case all groups $C\Gamma^i$, for $i\geq 2$ are trivial and
the sequence stabilizes.

\item $r=0$, $p=q=1$. In this case we case we have that $C\Gamma^2 =
\Z$ is abelian; this implies that $C\Gamma^3 = \Z^{l_2}$ is also
abelian and $C\Gamma^4 = \GL_{l_2}(\Z) \ltimes \Z^{l_2l_3}$ and
finally $C\Gamma^i = \GL_{l_2}(\Z) \ltimes \Z^{l_2l_3} \rtimes
\GL_{l_3}(\Z)$ for $i\geq 5$.

\item $r=0$ and $p.q > 1$ In this case we case we have $C\Gamma^2
= \Z^{pq}$ is abelian, therefore $C\Gamma^3 = \GL_{pq}(\Z) \ltimes
\Z^{pql_2}$ and finally $C\Gamma^i = \GL_{pq}(\Z) \ltimes
\Z^{pql_2} \rtimes \GL_{l_2}(\Z)$ for $i\geq 4$.

\item
$r=1$ and $p=q=0$. In this case we have two possibilities. The first one is $C\Gamma^2=\{1\}$,
in which case all groups $C\Gamma^i$ are trivial. The other is
$C\Gamma^2 = \F_2$, which leads to
$C\Gamma^3 = \F^{s_2}$ and $C\Gamma^4 = \SL_{s_2}(\F_2) \ltimes \F_2^{s_2s_3}$ and
$C\Gamma^i = \SL_{s_2}(\F_2) \ltimes \F_2^{s_2s_3} \rtimes \SL_{s_2}(\F_2)$ for all $i \geq 4$.

\item
$r=1$ and $p=0$ or $q=0$ (say $q=0$).
This case again splits into two subcases  --- in the first one
we have $C\Gamma^2 = \Z^{p}$, therefore the sequence $C\Gamma^i$ stabilizes as in cases
2 and 3. The other possibility is
$C\Gamma^2 = \F_2 \ltimes \Z^{p}$, which has trivial center.
If $p=1$ this is the infinite
dihydric group and it has non trivial automorphisms, i.e.,
$C\Gamma^i = C\Gamma^2$ for all $i$.
If $p \geq 1$ we have that
$C\Gamma^i = \GL_{p}(\Z) \ltimes Z^{p}$ for all $i\geq 3$.

\item
$r=1$ and $p=q=1$. In this case $C\Gamma^2$ is either the free two generated nilpotent
group of class $2$, or its extension by $\F_2$ in either case we have that
$C\Gamma^3 = \GL_2(\Z) \ltimes (\Z^2 \times Z^{l_2})$,
where $\GL$ acts naturally on $\Z^2$
and acts by multiplication with the determinant on $\Z^{l_2}$,
therefore we have that
$C\Gamma^i = (\GL_2(\Z) \times \GL_{l_2}(\Z))\ltimes (\Z^2 \times Z^{l_2})$,
for all $i > 3$.

\item $r=1$ and $p > 1$ or $q > 1$. In this case $C\Gamma^2$ is
either the nilpotent group of class $2$ having $\Z^{p+q}$ as
abelinization and $\Z^{pq}$ as center, or its extension by $\F_2$.
In either case we have that
$C\Gamma^3 = (\GL_p(\Z) \times \GL_q(\Z)) \ltimes (\Z^p\times \Z^g \times Z^{pql_2})$,
since this
group has trivial center we have that
$C\Gamma^i = (\GL_p(\Z) \times \GL_q(\Z)) \times \GL_{l_2}(\Z)) \ltimes (\Z^q \times \Z^p
\times Z^{pql_2})$, for all $i\geq 3$.

\item
$r\geq 2$ and $p=q=0$.
$C\Gamma^2 = \SL_r(\Z)$ or $C\Gamma^2 = \SL_r(\Z)$. Depending on the parity of $r$ we have that
$\Cen (C\Gamma^2)$ is either trivial or cyclic of order $2$, in the first case we have
$C\Gamma^i$ is  $PSL_r$ or $PGL_r$ for all $r\geq 3$. In the second we have
$C\Gamma^i = PSL_r(\Z) \times \F_2^{s_2}$, and as in case 4 the sequence $C\Gamma^i$ stabilizes
after 5 steps.

\item
$r\geq 2$ and $p=0$ or $q=0$ (say $q=0$).
In this case
$C\Gamma^2 = \SL_r(\Z) \ltimes \Z^{rp}$ or $C\Gamma^2 = \GL_r(\Z)\ltimes \Z^{rp}$.
In both cases we have that
$C\Gamma^3 = (\GL_r(\Z)\times \GL_p(\Z)) \ltimes \Z^{rp}$
and all others $C\Gamma^i$ are equal to $C\Gamma^3$.

\item $r\geq 2$ and $p,q\geq 1$. In this case $C\Gamma^2$ is an
extension of a nilpotent group of class $2$ having $\Z^{r(p+q)}$
as abelinization and $\Z^{pq}$ as center, by either $\SL_r(\Z)$ or
$\GL_z(\Z)$. In both cases we have that $ C\Gamma^3 = (\GL_r(\Z)
\times \GL_p(\Z) \times \GL_q(\Z)) \ltimes (\Z^{rp}\times \Z^{rq}
\times \Z^{pql_2}) $ and all other $C\Gamma^i$ are equal to
$$
C\Gamma^3 = (\GL_r(\Z) \times \GL_p(\Z) \times \GL_q(\Z) \times \GL_{l_2}(\Z))
 \ltimes (\Z^rp\times \Z^rq \times \Z^{pql_2})
$$
\end{enumerate}

Since there are no other possibilities for the integers $p$, $q$ and $r$,
we have shown that the sequence $C\Gamma^k$ stabilizes after at most
$5$ steps, which completes the proof.
\end{proof}

As an immediate corollary of this Lemma and Theorem~\ref{grtower:center} we obtain the
stabilization of the automorphism tower of the free nilpotent group $\Gamma(n,d)$
in the case when $\Cen(\Aut \Gamma)$ is not trivial.

\begin{thq}
\label{grtower:stab}
If $2n |d -1 $, then The automorphism tower of the free nilpotent group
$\Gamma(n,d)$ stabilizes after at most $\max \{5, 2+ d /h(n)\}$ steps, in particular it has
finite height.
\end{thq}

If we have that $d \leq h(n)$ we have that the tower stabilizes after at most $3$ steps.
\end{section}
\ifnewpage

\begin{section}{Two generated case}
\label{sec:2gen}

Finally let us consider the case $n=2$ and $d=4k+1$. In this case
we can not apply directly the results from the previous section
because $\SL_2(\Z)$ has homologically nontrivial cocycles into
$\Cen(UL)$ and therefore we can not embed the whole group $\Aut^2
\Gamma$ into $\Aut G^1$, but we can embed very big part of it.

The automorphism group $\Aut G^1$ can be described in the same way as
in Theorem~\ref{2ndlevel:center} part a). Lets us obtain a similar
description of the second group in the automorphism tower.

\begin{lm}
The group $\Aut^2 \Gamma$ contains a normal subgroup $A\Gamma^2$, which can be
embedded into $\Aut G^1$. Also the factor group $\Aut^2 \Gamma / A\Gamma^2$
is abelian. Also the group $A\Gamma^2$ slits as a direct product
and has description similar to the one Theorem~\ref{2ndlevel:center} b).
\end{lm}
\begin{proof}
Let $A\Gamma^2$ be the subgroup of $\Aut^2 \Gamma$ consisting of all
automorphisms which can be extended to the group $G^1$, i.e., the
cocycle $v:\GL_2(\Z)\to  \Cen(UG)$ is homologically trivial.
This group can be described as in
Theorem~\ref{2ndlevel:center} and therefore splits as a direct product.
The factor group $\Aut^2 \Gamma / A\Gamma^2$ is abelian because
$H^1(\GL_2(\Z), \Cen UG)$ is abelian.
\end{proof}

Let $G^2$ denote the Zariski closure of the group $A\Gamma^2$ in to $\Aut G^1$.
This group has similar description to the one in Theorem~\ref{2ndlevel:center}.

We would like to use the group $A\Gamma^2$ to study the next group in the
automorphism tower of $\Gamma$. Therefore we need to show that this is a
characteristic subgroup of $\Aut^2 \Gamma$.

\begin{lm}
\label{2gen:chgroup2}
The subgroup $A\Gamma^2$ is preserved by all automorphisms of
the group $\Aut^2 \Gamma$.
\end{lm}
\begin{proof}
Let $B$ is the commutator subgroup of $\Gamma$. This is a
subgroup of $A\Gamma^2$, and therefore can be embedded into some
quotient $G$ of $G^2$. Let us fix one such embedding $\phi : B \to G$.

Let $U$ denote the unipotent radical of $G$ and let
$C = \phi^{-1}(\phi B \cap U)$ be
the intersection of the group $B$ with the radical of $G$.
The rigidity Theorem~\ref{GP} implies the the subgroup $C$
does not depend on the choice of the embedding $\phi$, since
$G$ has strong unipotent radical.

Every automorphism of the group $B$ determines a cocycle form
the semi simple part of $B$ into the center of the unipotent radical of $G$.
This gives us a map $\pi$ form $\Aut^2 \Gamma / B$ to
$H^1(B/C, \Cen(U))$, because the conjugation with elements in
$\Aut^2 \Gamma / B$ gives automorphisms of the group $B$.
Applying Theorem~\ref{GP} again gives us  that
$\pi$ does not depend on the choice of $\phi$. Finally we can define
$A\Gamma^2$ as $\pi^{-1}(0)$. Which shows that $A\Gamma^2$ is
preserved by all automorphisms of the group $\Aut^2 \Gamma$.
\end{proof}

Combainning all this results we have shown that
\begin{thq}
\label{2gen:2nd}
The second group $\Aut^2 \Gamma$ in the automorphsim tower of $\Gamma$,
contains a inraviant subgroup $A\Gamma^2$, which is a Zariski dense
lattice int hte Lie group $G^2$. Also  
$A\Gamma^2$ splits as a direct product
of $N\Gamma^2$ and $C\Gamma^2$. The groups  $N\Gamma^2$ and $C\Gamma^2$
can be described in a way similar to Theorem~\ref{2ndlevel:center}.
It can also be shown that the whole group $\Aut^2 \Gamma$ splits as a direct product of
$\tilde N\Gamma^2$ and $C\Gamma^2$ for some abelian extension $\tilde N \Gamma^2$ of
$N\Gamma^2$
\end{thq}
\begin{proof}
The only thing that does nnot follow from the previous lemmas is that the
whole group $\Aut^2\Gamma$ splits as a direct product. This is true because
there are no notrivial cocycles in the center of $\Aut\Gamma$, which implies that
the conjugation with elements in $\Aut^2 \Gamma$ acts trivially on 
the group $C\Gamma^2$. 
\end{proof}

As in the case of more than $2$ generators, description of every group in 
the automorphism tower is similar to the one for the secon group. This
allows us to generalize the previous Theorem to a description of the
automorphism tower of the group $\Gamma$.

\begin{thq}
\label{2gen:tower}
a)
Every group $\Aut^k\Gamma$ in the automorphism tower of
$\Gamma$, contains a characteristic subgroup $A\Gamma^k$
such that all factors $\Aut^k \Gamma / A\Gamma^k$ for
$k \geq 2$ are abelian and isomorphic to each other.

b) Each group $A\Gamma^k$ can be embedded as a lattice in a Lie
group $G^k$ and $A\Gamma^{k+1}$ can be defined as the subgroup of
$\Aut G^k$ consisting of all automorphisms, which preserve the
subgroup $A\Gamma^k$.

c)
Each group $A\Gamma^k$ splits as a direct product of $N\Gamma^k$ and
$C\Gamma^k$ where the group $N\Gamma^k$ has description
similar to the ones in Theorem~\ref{grtower:center}. Also the group
$\Aut^k \Gamma$ splits as a direct product of $\tilde N \Gamma^k$ and
$C\Gamma^k$. The group $C\Gamma^k$ are defined recursivly as follows
$$
C\Gamma^{k+1} = \Aut C\Gamma^k \ltimes 
\Hom(\tilde N\Gamma^k/[\tilde N\Gamma^k,\tilde N\Gamma^k], \Cen(C\Gamma^k))
$$ 
\end{thq}
\begin{proof}
The proof is by induction and uses result of
Lemma~\ref{2gen:chgroup2} and Theorem~\ref{2gen:2nd} as a base case. 
The induction step is
similar to the one in Theorem~\ref{grtower:center}, 
but few additional things have to be shown.
First, we need to show that $A\Gamma^k$ is a invarinat subgroup of
$\Aut^k \Gamma$. The proof of this fact is the same as in the case $k=2$
and we can repeat the argument from Lemma~\ref{2gen:chgroup2}.
All groups $\Aut^k \Gamma / A\Gamma^k$ are isomorphic because
all nontrivial $\GL_2(\Z)$ modules in the center of the nilpotent part
of the group $A\Gamma^k$ comes from $\Cen(\IAut \Gamma)$.

Second, we need to mention that every automorphism of the group $\Aut^k \Gamma$ is
almost determined by its action on $A\Gamma^k$. Using the fact that the adjoint 
action of $\Aut^k \Gamma / \Cen (\Aut^k \Gamma)$ is faithfull we can see that
if $\phi \in \Aut(\Aut^k \Gamma)$ is such that the restriction 
$\phi_{|A\Gamma^k}$ of $\phi$ on $A\Gamma^k$ is the identity, then
$\phi$ acts trivially on $\Aut^k \Gamma / \Cen (\Aut^k \Gamma)$. That is the 
reason why in the
description of the group $C\Gamma^k$ we use the groups $\tilde N\Gamma^k$ not 
$N\Gamma^k$ as in Theorem~\ref{grtower:center}.

Finally we need to mention that the whole group $\Aut^k \Gamma$ splits as a 
direct product because all the nontriviall cocylces of $\GL_2(\Z)$ does not
interesect the center of the group $\Aut^k \Gamma$.
\end{proof}

As in the case of more than $2$ generators, we have that the sequences 
$N\Gamma^k$ and $\tilde N\Gamma^k$ stabilize, since the 
sequence $\tilde{UL}_k$ stabilizes. 
In order to claim that the automorphism towers of the $2$ generated 
nilpotent group stabilizes we need to show that the sequence $C\Gamma^k$
stabilizes.

We have that $C\Gamma^2$ is not trivial if and only if $d=4k+1$ for some
positive integer $k$. If $k\geq 2$ then Lemma~\ref{bigabelian}
holds, because the homogeneus component of degree $4k$ of the 
free Lie algebra on $2$ generators contains enough invarinat $\sll_2$ 
modules. 
In this case and the stabilization of the sequnece 
$C\Gamma^k$ follows from Lemma~\ref{CG:stab}

The case $d=5$ have to be examined separately because there are no $\sll_2$
invarinat modules in the homogeneous component of degree $4$ in the 
free Lie algebra on $2$ generators.
In this case we have that
$\Cen(\Aut\Gamma)= \Z$ and that $\IAut \Gamma^{\triv}_{\gen} = \Z$, therefore
we have that $C\Gamma^2$ is either trivial or cyclic of order $2$ (it is actually a
trivial but to show this we need to examin how $\Aut \Gamma$ sitts inside $G^1$).
If it is trivial then all groups $C\Gamma^k$ are trivial and the sequence stabilizes.

If it is not trivial we have to use the fact that 
$$
\tilde N\Gamma^2/[\tilde N\Gamma^2,\tilde N\Gamma^2] \simeq \Z \times \F_2,
$$
where the $\F_2$ copes from the abealinization of $\GL_2(\Z)$
and $\Z$ comes from the extension $\tilde N \Gamma^2$ over $N\Gamma^2$.
Using this information we can see that
$C\Gamma^3$ is elementary abelaing group of order $4$ and that
$$
C\Gamma^4 \simeq \GL_2(\F_2) \ltimes \Hom(\F_2^2,\F_2^2)
$$
and that all other group are equal to
$$
C\Gamma^k \simeq (\GL_2(\F_2) \times \GL_2(\F_2)) \ltimes \Hom(\F_2^2,\F_2^2)
$$
for $k\geq 5$.

This prove that the sequence $C\Gamma^k$ stabilizes after at most $6$ steps, which 
together with the stabilization of $\tilde N\Gamma^k$ finishes the proof of
the following theorem

\begin{thq}
\label{auttowerrst2gen}
The automorphism tower of the group $\Gamma(2,d)$ stabilizes after at most 
$\max \{ 6, 2 + d/4 \}$ steps, if $d=4k+1$.
\end{thq}

This result together with Theorems~\ref{auttower:st} and~\ref{grtower:stab}, shows that the
automorphim towers of the free nilpotent groups $\Gamma(n,d)$, stabilize
after finetly many steps, for any values of $n$ and $d$.

\end{section}

\end{chapter}
\newpage

\bibliographystyle{amsplain}
\bibliography{thesis}

\providecommand{\bysame}{\leavevmode\hbox to3em{\hrulefill}\thinspace}
\begin{thebibliography}{10}

\bibitem{Andre}
S.~Andreadakis, \emph{On the automorphisms of free groups and free nilpotent
  groups}, Proc. London Math. Soc. (3) \textbf{15} (1965), 239--268.

\bibitem{Borel}
Armand Borel and Nolan~R. Wallach, \emph{Continuous cohomology, discrete
  subgroups, and representations of reductive groups}, Annals of Mathematics
  Studies, vol.~94, Princeton University Press, Princeton, N.J., 1980.

\bibitem{Cheva}
C.~Chevalley, \emph{On groups of automorphism of {L}ie groups}, Proc. Nat.
  Acad. Sci. U. S. A. \textbf{30} (1944), 274--275.

\bibitem{DyerForm2}
Joan~L. Dyer and Edward Formanek, \emph{The automorphism group of a free group
  is complete}, J. London Math. Soc. (2) \textbf{11} (1975), no.~2, 181--190.

\bibitem{DyerForm1}
\bysame, \emph{Automorphism sequences of free nilpotent groups of class two},
  Math. Proc. Cambridge Philos. Soc. \textbf{79} (1976), no.~2, 271--279.

\bibitem{Formanek}
Edward Formanek, \emph{Fixed points and centers of automorphism groups of free
  nilpotent groups}, Comm. Algebra \textbf{30} (2002), no.~2, 1033--1038.

\bibitem{GrunPla}
Fritz Grunewald and Vladimir Platonov, \emph{Rigidity results for groups with
  radical cohomology of finite groups and arithmeticity problems}, Duke Math.
  J. \textbf{100} (1999), no.~2, 321--358.

\bibitem{Hall}
Philip Hall, \emph{The {E}dmonton notes on nilpotent groups}, Mathematics
  Department, Queen Mary College, London, 1969.

\bibitem{Hamkins}
Joel~David Hamkins, \emph{Every group has a terminating transfinite
  automorphism tower}, Proc. Amer. Math. Soc. \textbf{126} (1998), no.~11,
  3223--3226.

\bibitem{Hochschild}
G.~Hochschild, \emph{Automorphism towers of affine algebraic groups}, J.
  Algebra \textbf{22} (1972), 365--373.

\bibitem{Hulse}
J.~A. Hulse, \emph{Automorphism towers of polycyclic groups}, J. Algebra
  \textbf{16} (1970), 347--398.

\bibitem{Kass4}
M.~Kassabov, \emph{Automorphism tower of free nilpotent groups}, submitted to
  Electronic Research Announcements of AMS.

\bibitem{Kass5}
\bysame, \emph{An the automorphism group of free nilpotent groups and property
  $t$}, submitted to Journal of Algebra.

\bibitem{Kazhdan}
D.~A. Ka{\v{z}}dan, \emph{On the connection of the dual space of a group with
  the structure of its closed subgroups}, Funkcional. Anal. i Prilo\v zen.
  \textbf{1} (1967), 71--74.

\bibitem{LubPak}
Alexander Lubotzky and Igor Pak, \emph{The product replacement algorithm and
  {K}azhdan's property ({T})}, J. Amer. Math. Soc. \textbf{14} (2001), no.~2,
  347--363 (electronic).

\bibitem{Margulis2}
G.~A. Margulis, \emph{Discrete subgroups of semisimple {L}ie groups},
  Ergebnisse der Mathematik und ihrer Grenzgebiete (3) [Results in Mathematics
  and Related Areas (3)], vol.~17, Springer-Verlag, Berlin, 1991.

\bibitem{Mostow}
G.~D. Mostow, \emph{Strong rigidity of locally symmetric spaces}, Princeton
  University Press, Princeton, N.J., 1973.

\bibitem{Pettet}
Martin~R. Pettet, \emph{A note on the automorphism tower theorem for finite
  groups}, Proc. Amer. Math. Soc. \textbf{89} (1983), no.~1, 182--183.

\bibitem{Schenkman}
Eugene Schenkman, \emph{A theory of subinvariant {L}ie algebras}, Amer. J.
  Math. \textbf{73} (1951), 453--474.

\bibitem{Schenkman2}
\bysame, \emph{The tower theorem for finite groups}, Proc. Amer. Math. Soc.
  \textbf{22} (1969), 458--459.

\bibitem{Tolstykh}
Vladimir Tolstykh, \emph{The automorphism tower of a free group}, J. London
  Math. Soc. (2) \textbf{61} (2000), no.~2, 423--440.

\bibitem{Wang}
S.~P. Wang, \emph{On the {M}autner phenomenon and groups with property $({\rm
  {t}})$}, Amer. J. Math. \textbf{104} (1982), no.~6, 1191--1210.

\bibitem{Weiss}
Edwin Weiss, \emph{Cohomology of groups}, Pure and Applied Mathematics, Vol.
  34, Academic Press, New York, 1969.

\bibitem{Wielandt}
Helmut Wielandt, \emph{Eine verallgemeinerung der invarianten untergruppen},
  Math. Z. \textbf{45} (1939), 209--244.

\end{thebibliography}


\end{document}							